\documentclass[12pt]{article}
\usepackage[top=1.25in, bottom=1.25in, left=1in, right=1in]{geometry}
\usepackage{setspace}
\usepackage{amsfonts}
\usepackage{amsmath}
\usepackage{amsthm}
\usepackage{amssymb}
\usepackage{url}
\usepackage[colorlinks=true, allcolors=blue]{hyperref}
\usepackage{bm}
\usepackage{bbm}
\usepackage{mathtools}
\usepackage{booktabs}
\usepackage{xcolor}
\usepackage{graphicx}
\usepackage[round,semicolon]{natbib}
\usepackage[page]{appendix}
\DeclareMathOperator{\Diag}{Diag}
\DeclareMathOperator{\tr}{tr}
\DeclareMathOperator{\rank}{rank}
\DeclareMathOperator{\var}{Var}
\DeclareMathOperator*\argmin{\arg\!\min}

\newcommand{\E}{{\mathbb E}}
\newcommand{\Pa}{{\mathbb P}}
\newcommand{\R}{{\mathbb R}}
\newcommand{\N}{{\mathbb N}}
\newcommand{\Sy}{{\mathbb S}}
\newcommand{\Acal}{{\mathcal A}}
\newcommand{\Ccal}{{\mathcal C}}
\newcommand{\Dcal}{{\mathcal D}}
\newcommand{\Ecal}{{\mathcal E}}
\newcommand{\Fcal}{{\mathcal F}}
\newcommand{\Gcal}{{\mathcal G}}
\newcommand{\Hcal}{{\mathcal H}}
\newcommand{\Ical}{{\mathcal I}}
\newcommand{\Tcal}{{\mathcal T}}
\newcommand{\Jcal}{{\mathcal J}}
\newcommand{\Kcal}{{\mathcal K}}
\newcommand{\opn}[1]{\|#1\|_{\rm op}}
\newcommand{\Lcal}{{\mathcal L}}
\newcommand{\Mcal}{{\mathcal M}}
\newcommand{\Ncal}{{\mathcal N}}
\newcommand{\Qcal}{{\mathcal Q}}
\newcommand{\Ucal}{{\mathcal U}}
\newcommand{\Xcal}{{\mathcal X}}
\newcommand{\HH}{{\mathbb H}}
\newcommand{\Ocal}{{\mathcal O}}
\newcommand{\Pcal}{{\mathcal P}}
\newcommand{\Rcal}{{\mathcal R}}
\newcommand{\Scal}{{\mathcal S}}
\newcommand{\Vcal}{{\mathcal V}}

\newcommand{\Zcal}{{\mathcal Z}}
\newcommand{\rmix}{\alpha_{\rm mix}}
\newcommand{\e}{{\rm e}}

\providecommand{\joesubmission}{0}
\newcommand{\figorbox}[2][]{\IfFileExists{#2}{\includegraphics[#1]{#2}}{\typeout{WARNING: figure #2 is missing}\fbox{\parbox[c][0.2\textheight][c]{0.6\textwidth}{\centering\footnotesize figure missing\\\texttt{\detokenize{#2}}}}}}

\newtheorem{proposition}{Proposition}[section]
\newtheorem{lemma}[proposition]{Lemma}
\newtheorem{theorem}[proposition]{Theorem}

\newtheorem{corollary}[proposition]{Corollary}
\newtheorem{remark}[proposition]{Remark}
\newtheorem{assumption}[proposition]{Assumption}

\author{
 Paul\ Schneider\footnote{
 Universit\`a della Svizzera italiana and Swiss Finance Institute, Email: \texttt{paul.schneider@usi.ch}
 }}
\date{\today}

\begin{document}

\title{Covariate-Indexed Inclusion Probabilities on Unbalanced Panels\thanks{I thank Damir Filipovic, without whom this paper would not have happened. I am grateful to the SNF grant 215528 \emph{Large-scale kernel methods in financial economics}.\ifnum\joesubmission=0{} During the preparation of this work I used Claude (Anthropic) to check proofs, write and test code for the simulations and the application, search the literature and edit the text. I reviewed and edited all output and take full responsibility for the content of the paper.\fi}}
\maketitle

%%%%%%%%%%%%%%%%%%%%%%%%%%%%%%%%%%%%%%%%%%%%%%%%%%%%%%%%%%%%%%%%%%%%%%%%
\begin{abstract}
We estimate the probability that an item with given covariates belongs to a set when the cross-section changes over time and items and periods are dependent. We propose two convex models whose fitted probabilities are nonnegative by construction and remain valid whatever the dependence between items, for which we establish consistency, finite-sample bounds and limit distributions for serially dependent panels.
One suits events of moderate or high probability, the other rare events, and half the disagreement between their fits estimates a lower bound on the misspecification error of the worse fit, without knowledge of the true probabilities. On monthly US stock failures, the rare-event model predicts better than the ordinary logit and the classical failure logit. With the eight classical predictors, most of the gain comes from basis functions of the rank angle of each predictor and summary.

\vspace{2ex}

\vspace{1ex}
\end{abstract}
\clearpage

%%%%%%%%%%%%%%%%%%%%%%%%%%%%%%%%%%%%%%%%%%%%%%%%%%%%%%%%%%%%%%%%%%%%%%%%
\section{Introduction}\label{sec_introduction}

Many panel data sets record a set, such as the obligors that default, the policies that produce a claim, the stocks that trade or land in the tails of the return distribution, the species present at surveyed sites, or the items placed in a basket. The cross-section over which the subset is drawn changes from period to period, so the data are unbalanced. A basic quantity attached to such data is the probability that an item with given covariates belongs to the realized subset, the conditional marginal inclusion probability. It is the object estimated by probability-of-default models, occupancy models, first-order inclusion probabilities in survey sampling, and marginal relevance models in information retrieval.

An estimator of this object has to meet three requirements, which we state before any method. (R1) Its output is nonnegative for every covariate value and every composition of the cross-section, the constraint of a probability that a linear index violates, and calibrated, in the sense that fitted probabilities match realized frequencies, over a class of functions flexible enough to follow the data. (R2) It is valid under arbitrary dependence between the items of a cross-section, since the defaults of a month or the species at a site are not independent and their dependence is rarely known. (R3) It supports inference when the periods are serially dependent and the cross-section changes in size and composition from one period to the next. The classical answer is logistic regression with a linear index. Like any generalized linear model, it fixes a link, a function of the probability that the index represents, here the log-odds \citep{mccullaghnelder89}. It meets (R2), since any marginal quasi-likelihood is calibrated to the margins whatever the dependence between items, and (R3), since the asymptotics of a parametric quasi-likelihood under mixing are standard \citep{white94,wooldridge94,PoetscherPrucha1997}, but it fails the flexibility part of (R1). Under the logit link the index is unconstrained, so making it flexible is easy \citep[e.g.,][]{sadhwanigieseckesirignano21}. Under the identity link the index is the probability itself, and under the log link it is minus its logarithm, as in the log-binomial model \citep{wacholder86}. A flexible index must then remain nonnegative at every covariate value. Flexible classes that are nonnegative by construction exist for the identity link (Section~\ref{sec_lit}), but not, to our knowledge, for the log link, and the guarantees of (R3) for the constrained estimator that such a class requires are missing for both. This is the gap the paper fills.

We meet the three requirements with one family of estimators, assembled from two established pieces. The inclusion-probability curve is a sum-of-squares quadratic form in a finite dictionary of functions on the covariate space, nonnegative by construction for every covariate value and every cross-sectional composition, and estimation is a finite-dimensional convex program over the positive semidefinite cone. We adopt this positive semidefinite (psd) model from \citet{marteauferey2020nonparametric}. The loss is the generalized Kullback--Leibler contrast between the observed inclusion indicators and the fitted curve, that is, the Poisson log-likelihood read as a quasi-likelihood. Each requirement is met by one ingredient. The cone delivers nonnegativity, and the first-order conditions of the fit are calibration equations, which is (R1). The upper bound is built in by a second link. Read as the negative logarithm of the probability and fitted by the Bernoulli likelihood, again a convex program over the cone, the sum of squares gives probabilities in $(0,1]$ (Section~\ref{sec_loglink}). Read as the probability itself, it can exceed one, which matters little for inclusion probabilities and can be prevented (Section~\ref{sec_loss}). For rare events the log link is the natural reading. Its calibration equations weight the residuals nearly uniformly, where those of the identity link weight them by the inverse fitted probability, and to first order it is a logistic regression whose index is minus a sum of squares. The quasi-likelihood property, that the unconstrained population minimizer is the conditional mean whatever the rest of the law (Section~\ref{sec_lit}), is (R2). The theory of this paper, consistency, finite-sample bounds and limit distributions for a cone-constrained estimator under strong mixing, is (R3). The two pieces have been combined before for Poisson intensities, the continuum counterpart of our setting (Section~\ref{sec_lit}).

There is one limit that cannot be removed through the choice of the dictionary or the loss. Any model of the own covariate sees the margins only through their average over the cross-sectional compositions in which a covariate value occurs.
% , the composition-averaged margin $\bar p$ of Section~\ref{sec_pbar}.
When the margins depend on the composition of the cross-section, as they can when items are dependent \citep{dasetal07,duffieecknerhorelsaita09,azizpourgieseckeschwenkler18} and do by construction when the number of events per period is fixed on a panel of varying size, the distance between the margins and the composition-averaged margin
% $\bar p$
is a floor under the error of every marginal model, whatever its capacity (Proposition~\ref{prop_target}\ref{prop_target_comp}). The floor is visible in simulation, where under a repulsive law the error of the fitted curve settles at about two percentage points as the sample size grows and is barely moved by doubling the dictionary.
% , from $0.0207$ to $0.0189$ (Table~\ref{tab_e1}).
The floor is relative to the information set of the marginal model. Covariates that summarize the cross-section lower it, and a model of the subset, such as a marked point process for the event times \citep{gieseckeschwenkler18}, is the alternative where no such summary is available.

We make three contributions. First, the estimand and its identification. The population target of the estimator is the generalized Kullback--Leibler projection of the margins onto the sum-of-squares class, which exists for a nondegenerate bounded dictionary and is unique under any joint law of the observed subsets with positive margins whenever the class gap is finite (Proposition~\ref{prop_target}\ref{prop_target_proj} and~\ref{prop_target_comp}), and the gap of this projection splits exactly into the composition part just described, which no dictionary removes, and a dictionary part, which a rich enough dictionary removes under the density conditions of Section~\ref{sec_loss} (Proposition~\ref{prop_target}\ref{prop_target_comp}).

Second, inference under dependence over a cone. Existing guarantees for sum-of-squares and positive semidefinite models are i.i.d.\ excess-risk bounds \citep{marteauferey2020nonparametric,rudiciliberto21}, and self-concordance gives fast rates for regularized empirical risk minimization under independent sampling \citep{marteauferey19}. We give almost-sure consistency of any empirical minimizer over the cone, with existence derived rather than assumed, finite-sample concentration under i.i.d.\ sampling and under geometric strong mixing, and a mean-squared-error bound at the parametric rate on the event of localization (Theorem~\ref{thm_main}). We also give the limit distribution of the fitted curve, for the identity link and, under the mirror image of one condition, for the log link (Theorem~\ref{thm_log}). It is Gaussian with a sandwich long-run covariance at an interior target, its projection onto the tangent cone of the constraint at a boundary target with vanishing multiplier, and Gaussian on the tangent space of the manifold of fixed rank at a boundary target under strict complementarity, the limit that \citet[Theorem~1]{shapiro19} gives for sample semidefinite programs, here for a quasi-likelihood over the cone under strong mixing (Theorem~\ref{thm_clt}). Because the fitted curve is linear in the parameter matrix, every statement transfers to the curve without a delta-method step, and under the log link to the probability with one. The conditions of the two fits fail at opposite ends of the probability range, so the identity link suits events of moderate or high probability and the log link rare events (Section~\ref{sec_theory_log}).

Third, a diagnostic of the model class that needs no knowledge of the margins. The identity and the log links represent different curves of the margins,
% $p$ and $-\log p$, by a sum of squares,
and half the disagreement between the two fits bounds the larger of their two errors on the probability scale from below (Proposition~\ref{prop_scalediag}), and through it the larger of the two dictionary gaps, with a limit distribution that gives it a standard error. It sees the dictionary part of the gap and is blind to its composition part, which both fits share. The application illustrates the estimators and the diagnostic. On monthly failures of US common stocks, a rare event, the log-link fit in a rank-angle dictionary predicts better in log score than the ordinary logit on $121$ characteristics and the failure logit of \citet{campbellhilscherszilagyi08}, and the identity-link fit ranks the names less well and shows the loose calibration at high probabilities where a class gap of the identity link shows.

Section~\ref{sec_estimand} defines the estimand and what any marginal model can recover. Section~\ref{sec_estimator} introduces the two estimators, their population targets and their calibration equations. Section~\ref{sec_theory} contains the statistical theory. Section~\ref{sec_diag} develops the diagnostic. Sections~\ref{sec_sim} and~\ref{sec_app} report the simulation study and the application, and Section~\ref{sec_conclusion} concludes. Appendix~\ref{app_proofs} contains the proofs and the auxiliary lemmas and Appendix~\ref{app_faceband} the bands at a boundary target. The internet Appendix~\ref{app_internet} gives the projection band at a zero multiplier, the identification of the fitted matrix in the application, the simulation in the regime of the application, the missing covariates, the masked dictionary and further results.

\subsection{Related literature}\label{sec_lit}

In this section we set out the antecedents of each component and mark where our contribution begins.
The sum-of-squares parametrization on the positive semidefinite cone is the positive semidefinite model of \citet{marteauferey2020nonparametric}. They establish its representer theorem, the convexity of the resulting problems, and its approximation properties, and they motivate it as we use it, as the nonnegative analogue of a linear model that preserves convexity. The construction descends from the sum-of-squares representation of nonnegative polynomials and the associated moment hierarchies \citep{lasserre10}. Later work extends it to matrix-valued functions \citep{muzellec21psd}, to probability models \citep{rudiciliberto21}, and to squared neural parametrizations \citep{tsuchida24}. In econometrics the same device underlies the constrained polynomial likelihood of \citet{almeidamasinischneider25}. Positive semidefinite quadratic forms read as class probabilities and fitted by maximum likelihood appear in \citet{crammerglobersonuai06}.

The loss used with the identity link is the negative log-likelihood of independent Poisson counts, evaluated at binary outcomes. Consistency for a correctly specified conditional mean, irrespective of the rest of the distribution, is the defining property of quasi-likelihood \citep{wedderburn74} and of pseudo-maximum likelihood in a linear exponential family \citep{gmt84poisson,gmt84theory}, and we use the first term throughout. This property explains the widespread use of the Poisson quasi-likelihood on nonnegative data in economics \citep{silvatenreyro06}. Equivalently it is a Bregman divergence, and the unique minimizer of an expected Bregman divergence is the mean \citep{savage71,banerjeeguowang05,gneitingraftery07}. Estimating marginal means under unspecified within-cluster dependence, with a sandwich variance in place of a correctly specified joint law, is the working-independence generalized estimating equation of \citet{liangzeger86}, and applying it to binary outcomes through a Poisson quasi-likelihood is standard practice in epidemiology \citep{zou04}. Proposition~\ref{prop_target} specializes these results to inclusion indicators, and the sandwich covariance of Proposition~\ref{prop_mds} belongs to the same tradition.

The two ingredients have been combined before, for intensities. \citet{alizadehetal08} and \citet{pappalizadeh14} fit nonnegative splines through a positive semidefinite representation, \citet{flaxman17} estimate a Poisson intensity as the square of a function in a reproducing kernel Hilbert space, the rank-one case of our objective, and \citet{walderbishop17}, \citet{tsuchida24} and \citet{mutnykrause22} use squared or positive semidefinite parametrizations with the same likelihood. Our loss is the finite-ground-set version of theirs. The summands are Bernoulli indicators on a varying ground set, so the fitted quantity is a probability, whose bound by one the sum of squares does not impose, and the index $t$ carries serial dependence.

Consistency and asymptotic normality of extremum estimators under mixing are classical \citep{white94,PoetscherPrucha1997}. At a boundary of the parameter set the limit is the projection law of \citet{gey_94} and \citet{andrews1999b} when the multiplier vanishes \citep[see also][]{chernoff54,selfliang87}, and Gaussian on a linear space under strict complementarity \citep{shapiro89}, for semidefinite programs on the tangent space of the manifold of fixed rank with the curvature term of \citet{shapiro97} \citep[Theorem~1]{shapiro19}. Theorem~\ref{thm_clt} carries these limits to the risks of this paper under strong mixing. Unlike sieve estimation \citep{gallantnychka87,wongshen95}, the dictionary is held fixed, as within each window of the application. Pseudo-self-concordance gives fast rates without global strong convexity \citep{bach10}, and \citet{marteauferey19} and \citet{ost_bac_21} use it for excess-risk and parameter deviation bounds, all under independent sampling and with a vanishing gradient at the target. Section~\ref{sec_theory} extends the localization to a cone whose boundary the target may touch with a nonzero multiplier, to a linear penalty and to strongly mixing data.

Comparing two estimators that share a probability limit under correct specification and not otherwise is the principle of \citet{hausman78}, and \citet{wooldridge91} applies it to quasi-maximum likelihood estimators of a conditional mean, neither of which needs to be efficient. The link discrepancy of Section~\ref{sec_diag} compares two such fits. Each projects onto its own class, so the discrepancy bounds a misspecification error from below rather than testing a null.

%%%%%%%%%%%%%%%%%%%%%%%%%%%%%%%%%%%%%%%%%%%%%%%%%%%%%%%%%%%%%%%%%%%%%%%%
\section{What can be estimated}\label{sec_estimand}
% Question: what is the target, and what can any marginal model recover? Before any loss.

\subsection{Data and estimand}\label{sec_setup}

We observe a training set $\Dcal_T=\{\xi_t:t=0,1,\dots,T-1\}$ of $T$ samples from an unknown distribution. Each sample $\xi_t=(n_t,S_{t+1},\bm z_t)$ consists of a cross-sectional size $n_t\in\N$, free to vary across samples, covariates $\bm z_t=[z_{t,i}:1\le i\le n_t]\in\Zcal^{n_t}$ in a fixed covariate space $\Zcal$, and an observed subset $S_{t+1}$ of the ground set $\{1,\dots,n_t\}$, for instance the items that default, claim, trade, or activate at time $t+1$. The index $t$ may denote time or, more generally, a task. The ground set may vary in both size and composition.

The estimand is the conditional marginal inclusion probability
\begin{equation}\label{ptidef}
    p_{t,i}\coloneqq\Pa\big(i\in S_{t+1}\mid n_t,\bm z_t\big),\qquad 1\le i\le n_t,
\end{equation}
as a function of the item's covariates. The joint law of $S_{t+1}$ is left unrestricted. Items within a cross-section may be dependent, and we impose no parametric form and no restriction on the dependence between inclusion events. The theory concerns the margins~\eqref{ptidef} alone, and the identification result (Proposition~\ref{prop_target}) holds for any joint law with positive margins.

The data points take values in the sample space $\Xi\coloneqq\bigcup_{n\in\N}\big(\{n\}\times 2^{\{1,\dots,n\}}\times\Zcal^n\big)$, the union over all cross-sectional sizes. We write $\xi=(n,S,\bm z)$ for a generic point in $\Xi$, and $y_{t,i}\coloneqq\mathbbm{1}_{\{i\in S_{t+1}\}}$ for the inclusion indicators.

\subsection{What a marginal model can recover}\label{sec_pbar}

A marginal model, in the sense of this paper, fits a function of the item's own covariate to the indicators. Its fitted value at item $i$ of period $t$ is a function of $z_{t,i}$ alone, whereas the margin $p_{t,i}$ is in general a functional of the whole conditioning information $(n_t,\bm z_t)$. Two items with identical covariates in differently composed cross-sections can then have different inclusion probabilities, and when the number of events per period is fixed by construction while the cross-section varies in size, as for the $k$ lowest returns of a month, the margins sum to $k$ and so depend on how many, and in general on which, other items are in the cross-section. This composition dependence is a property of the margins given the conditioning set, not of the joint law. Items driven by a common shock can have margins that depend on their own covariates alone, and independent items can have margins that depend on the composition. What a function of the own covariate can recover is the following average. Throughout, $\{\xi_t\}_{t\ge0}$ is strictly stationary with marginal law $\Pa$, and all population objects refer to this marginal.

Write $\mu_1(A)\coloneqq\E\big[\sum_{i\le n_0}\mathbbm{1}_{\{z_{0,i}\in A\}}\big]$ for the \emph{item intensity} of the covariates, the expected number of items per period with covariate in $A$. It is a measure on $\Zcal$ of total mass $\E[n_0]$, which we assume finite from here on, so that $\mu_1$ is finite. Write further $\mu_p(A)\coloneqq\E\big[\sum_{i\le n_0}p_{0,i}\mathbbm{1}_{\{z_{0,i}\in A\}}\big]\le\mu_1(A)$ for the \emph{intensity of the margins}, the expected number of events per period with covariate in $A$, so that $\mu_p\ll\mu_1$ and the \emph{composition-averaged margin}
\begin{equation}\label{pbardef}
    \bar p\coloneqq\frac{d\mu_p}{d\mu_1}\in[0,1]\quad\mu_1\text{-a.e.}
\end{equation}
is the margin of an item given only its own covariate, averaged over the compositions in which that covariate occurs. It is a function on $\Zcal$, whereas $p_{0,i}$ is a function of all of $(n_0,\bm z_0)$.
Since $\E[y_{0,i}\mid n_0,\bm z_0]=p_{0,i}$, $\mu_p(A)=\E\big[\sum_{i\le n_0}y_{0,i}\mathbbm{1}_{\{z_{0,i}\in A\}}\big]$ is the expected number of included items per period whose covariates fall in $A$, the mean measure of the selected subset read as a point pattern, with total mass $\mu_p(\Zcal)=\E|S_1|$, the expected number of events per period. In the language of marked point processes, with the covariates as points and the indicators as marks, $\mu_1$ is the first moment measure of the ground process and $\bar p(z)$ is the probability of the mark one at location $z$ under the mark kernel \citep{daleyverejones03}.
It depends on the joint law of the subset only through its margins. If the margins depend on the own covariate alone, $p_{0,i}=\pi(z_{0,i})$, then $\mu_p=\pi\,\mu_1$ and $\bar p=\pi$ $\mu_1$-a.e. Otherwise, $\mu_p$ adds up the margins of the items near $z$ over all the compositions in which they occur.
For every Borel $f\ge0$ on $\Zcal$,
\begin{equation}\label{pbarident}
    \E\Big[\sum_{i\le n_0}p_{0,i}f(z_{0,i})\Big]=\int\bar p\,f\,d\mu_1,
\end{equation}
and $f(z_{0,i})=0$ for all $i\le n_0$ $\Pa$-a.s.\ if and only if $f=0$ $\mu_1$-a.e. For $f=\mathbbm{1}_A$ the identity is the definition of $\mu_p$ and of $\bar p$, and it extends to Borel $f\ge0$ by monotone convergence. The second claim is $\mu_1(\{f\ne0\})=\E\big[\sum_{i\le n_0}\mathbbm{1}_{\{f(z_{0,i})\ne0\}}\big]$.

The identity is elementary, and it fixes the scope of every marginal model. An expected loss that is affine in the margin at a fixed fitted value, as the losses of Sections~\ref{sec_estimator} and~\ref{sec_diag} are, is by~\eqref{pbarident} a functional of $(\mu_1,\bar p)$ alone. Whatever the dictionary, and for every loss that adds over the items a term in the item's own covariate, its fitted value and its indicator, the fit is scored against a panel whose margins are $\bar p(z_{t,i})$, and the dependence of $p_{t,i}$ on the rest of the cross-section is invisible to it. The difference between $p_{t,i}$ and $\bar p(z_{t,i})$ is thus a \emph{composition part} of the error that no marginal model removes.

%%%%%%%%%%%%%%%%%%%%%%%%%%%%%%%%%%%%%%%%%%%%%%%%%%%%%%%%%%%%%%%%%%%%%%%%
\section{The estimators}\label{sec_estimator}
To adapt to situations where estimated probabilities might be particularly large or small, and to diagnose the model class, we introduce two estimators. Each is calibrated to the margins on its own scale, and their disagreement measures the link-specific part of the dictionary gap (Section~\ref{sec_diag}). No marginal model, either fit included, sees the composition part (Section~\ref{sec_pbar}).

\subsection{The sum-of-squares index}\label{sec_dict}

A dictionary of $m$ linearly independent Borel functions $\phi_1,\dots,\phi_m:\Zcal\to\R$ is collected in the row-vector-valued feature map $\bm\phi\coloneqq[\phi_1,\dots,\phi_m]$. Any basis functions are admissible, for instance polynomials, trigonometric or spline systems, random features, or the data-driven dictionary of kernel sections at pivoted training covariates of \citet{filipovicschneider24mom}, which lets the data choose the dictionary through a pivoted Cholesky reduction of an RKHS kernel. The candidate indices are the quadratic forms
\begin{equation}\label{elldef}
    \ell_{\bm U}(z)\coloneqq\bm\phi(z)\,\bm U\,\bm\phi(z)^\top,\qquad \bm U\in\Sy_+^m,
\end{equation}
the sum-of-squares parametrization of nonnegative functions. Factoring $\bm U=\bm L\bm L^\top$ gives $\ell_{\bm U}(z)=\|\bm\phi(z)\bm L\|_2^2\ge0$, so the positivity constraint in function space becomes the cone constraint on $\bm U$. With the trace inner product $\langle\bm A,\bm B\rangle\coloneqq\tr(\bm A\bm B)$ on the symmetric matrices, the index is linear in the parameter,
\begin{equation}\label{cdef}
    \ell_{\bm U}(z)=\langle\bm C(z),\bm U\rangle,\qquad
    \bm C(z)\coloneqq\bm\phi(z)^\top\bm\phi(z)\in\Sy_+^m.
\end{equation}
% which later turns statements about the parameter into statements about the curve without a delta-method step.
On a cross-section with feature matrix $\bm\Phi_t\coloneqq\bm\phi(\bm z_t)=[\bm\phi(z_{t,i}):1\le i\le n_t]\in\R^{n_t\times m}$ we write $\bm C_{t,i}\coloneqq\bm C(z_{t,i})$ and
\begin{equation}\label{elltidef}
    \ell_{t,i}\coloneqq\ell_{\bm U}(z_{t,i})=\big[\bm\Phi_t\bm U\bm\Phi_t^\top\big]_{ii},
\end{equation}
suppressing $\bm U$ where no confusion arises. The parametrization has $m(m+1)/2$ free parameters, independent of the cross-sectional sizes and of the total number of observed items.

Panels of firm characteristics have missing coordinates, and in the application only a third of the name-months observe every characteristic. Rather than at an imputed value, a kernel dictionary can be evaluated at the conditional law of the missing coordinates given the observed ones \citep{smolavishwanathanhofmann05,songetal09}, as \citet{smiejaetal18} do for the first layer of a neural network. Write $o=o(z)$ for the observed coordinates of a covariate vector $z$ in $\Zcal\subseteq\R^d$ and $\bar o$ for the missing ones, let $\bm e$ be a dictionary of kernel sections at pivoted training covariates, orthonormal in the reproducing kernel Hilbert space of the kernel, and let $w_1(z),\dots,w_{N_{\rm ref}}(z)\ge0$ be weights that sum to one on complete reference records $x_1,\dots,x_{N_{\rm ref}}$, an estimate of the conditional law of the missing coordinates given $X_o=z_o$. The \emph{masked dictionary}
\begin{equation}\label{maskeddict}
    \bm\phi(z)\coloneqq\sum_{i=1}^{N_{\rm ref}}w_i(z)\,\bm e(z_o,x_{i,\bar o})
\end{equation}
keeps the observed coordinates at their values and averages the basis over completions of the missing ones. A complete $z$ gets $\bm\phi(z)=\bm e(z)$, and the index $\ell_{\bm U}=\bm\phi\bm U\bm\phi^\top$ is nonnegative, linear in $\bm U$ and fitted by the same programs. Appendix~\ref{app_masked} gives the construction, its properties and the conditions under which the weights estimate the conditional law.

Below, we introduce two link functions.
The identity link reads $\ell_{\bm U}$ as the probability (Section~\ref{sec_loss}), the log link as minus its logarithm (Section~\ref{sec_loglink}). We name the two estimators by their links, the identity-link fit and the log-link fit. The superscripts ${\rm Po}$ and ${\rm log}$ on their losses, targets and scores refer to the Poisson loss and the log-link loss. Both fits solve the regularized problem
\begin{equation}\label{ermdef}
    \underset{\bm U\in\Sy_+^m}{\text{minimize }}\ \Rcal_T(\bm U)\coloneqq\frac1T\sum_{t=0}^{T-1}\Lcal(\bm U,\xi_t)+\lambda\tr(\bm G\bm U),
\end{equation}
with the loss $\Lcal=\Lcal^{\rm Po}$ of~\eqref{lossdef} or $\Lcal=\Lcal^{\rm log}$ of~\eqref{lossloglink}, a regularization parameter $\lambda\ge0$, and a weight matrix $\bm G\in\Sy_+^m$. For a kernel-generated dictionary, $\bm G$ is the Gram matrix of the dictionary functions in their reproducing kernel Hilbert space, the penalty is the sum of the squared RKHS norms of the factors $\bm\phi\bm t_k$ of $\bm U=\sum_k\bm t_k\bm t_k^\top$, and the pivoted Cholesky step selects $m$ kernel sections whose span approximates the one in which a representer theorem places the solution \citep{filipovicschneider24mom}, with the residual bounded in Proposition~\ref{prop_masked_app}\ref{prop_masked_err}. The linear and rank-angle dictionaries of the application are the full feature maps of their kernels, so no reduction is made there. The choice $\bm G=\bm I_m$ gives the trace norm, which encourages low-rank solutions. Both problems are finite-dimensional convex programs, representable with semidefinite and exponential cones, and the feasible set may be any closed convex subset of the cone. The fitted curve keeps the letter $\ell$ throughout, $\hat p$ denotes a fitted probability, and $\Lcal$ always denotes an objective.

\subsection{The identity link and the Poisson quasi-likelihood}\label{sec_loss}

The identity link reads the sum of squares as the probability itself, $\hat p=\hat\ell$, and fits it by the \emph{Poisson loss}
\begin{equation}\label{lossdef}
    \Lcal^{\rm Po}(\bm U,\xi_t)\coloneqq\sum_{i=1}^{n_t}\ell_{t,i}-\sum_{i\in S_{t+1}}\log\ell_{t,i},
\end{equation}
with the conventions $\sum_{i\in\emptyset}=0$ and $-\log x\coloneqq+\infty$ for $x\le0$. Up to a term free of $\bm U$, \eqref{lossdef} is the negative log-likelihood of independent Poisson counts $y_{t,i}\sim{\rm Poisson}(\ell_{t,i})$ evaluated at binary outcomes, the Poisson quasi-likelihood \citep{wedderburn74,gmt84poisson} in the sum-of-squares parametrization. In the continuum this pairing is the squared-kernel intensity estimator of \citet{flaxman17} and \citet{tsuchida24}, and Poisson quasi-likelihood for binary outcomes is routine in marginal modeling \citep{liangzeger86,zou04}. Its population target is a projection in the generalized Kullback--Leibler divergence
\begin{equation}\label{dkldef}
    d_{\rm KL}(p\,\|\,\ell)\coloneqq p\log\frac{p}{\ell}-p+\ell\ \ge\ 0,\qquad p\ge0,\quad\ell>0,
\end{equation}
with $0\log0\coloneqq0$, extended by $d_{\rm KL}(p\,\|\,0)\coloneqq+\infty$ for $p>0$, which vanishes if and only if $\ell=p$. It is the Bregman divergence \citep{bregman67} generated by $x\log x-x$. Unlike the divergence between probability vectors it needs no normalization, which is what the margins require, since they sum to the expected subset size and not to one. The linear term $-p+\ell$ makes $\ell=p$, rather than any multiple of $p$, the unique minimizer. We write $\bm U^{\rm Po}_\star$ for a minimizer over $\Sy_+^m$ of the Poisson risk $\bm U\mapsto\E[\Lcal^{\rm Po}(\bm U,\xi_0)]$, where one exists, and $\ell^{\rm Po}_\star\coloneqq\ell_{\bm U^{\rm Po}_\star}$ for its curve. With a penalty the target moves to $\bm U^{\rm Po}_{\star,\lambda}$ of Section~\ref{sec_theory_po}, and $\bm U^{\rm Po}_{\star,0}=\bm U^{\rm Po}_\star$. The next proposition collects the properties of the fit, with the \emph{residual weight} $w_{\rm Po}(\ell)\coloneqq1/\ell$, the ratio of the derivative of the mean map to the working variance of the Poisson contrast.

\begin{proposition}[The identity-link fit]\label{prop_target}
\begin{enumerate}
\item\label{prop_target_convex} \emph{Convexity.} $\bm U\mapsto\Lcal^{\rm Po}(\bm U,\xi_t)$ is convex and lower semicontinuous on $\Sy^m$.
\item\label{prop_target_unc} \emph{Unconstrained target.} Let $S\subseteq\{1,\dots,n\}$ be a random subset with margins $p_i\coloneqq\Pa(i\in S)>0$ under any joint law. The Poisson risk over unconstrained $\bm\ell\in(0,\infty)^n$,
\begin{equation}\label{rhodef}
    \rho^{\rm Po}(\bm\ell)\coloneqq\E\Big[\sum_{i=1}^n\ell_i-\sum_{i\in S}\log\ell_i\Big]=\sum_{i=1}^n\big(\ell_i-p_i\log\ell_i\big),
\end{equation}
is strictly convex with unique minimizer $\ell_i=p_i$, $i\le n$.
\item\label{prop_target_proj} \emph{Projection onto the class.} Let $\{\xi_t\}_{t\ge0}$ be strictly stationary with $n_t\le\bar n$ and $p_{0,i}>0$ $\Pa$-a.s. For every $\bm U\in\Sy_+^m$, as an identity in $[0,\infty]$ and without any moment condition,
\begin{equation}\label{klident}
    \E\big[\Lcal^{\rm Po}(\bm U,\xi_0)\big]=\E\Big[\sum_{i=1}^{n_0}d_{\rm KL}\big(p_{0,i}\,\big\|\,\ell_{\bm U}(z_{0,i})\big)\Big]+c_0^{\rm Po},\qquad c_0^{\rm Po}\coloneqq\E\Big[\sum_{i=1}^{n_0}\big(p_{0,i}-p_{0,i}\log p_{0,i}\big)\Big],
\end{equation}
with $c_0^{\rm Po}\in[0,\bar n]$ free of $\bm U$. The unregularized population target is the $d_{\rm KL}$-projection of the margins onto the class $\{\ell_{\bm U}:\bm U\in\Sy_+^m\}$, with \emph{class gap}
\begin{equation}\label{epsdictdef}
    \epsilon_{\rm class}^{\rm Po}\coloneqq\inf_{\bm U\in\Sy_+^{m}}\E\Big[\sum_{i=1}^{n_0}d_{\rm KL}\big(p_{0,i}\,\big\|\,\ell_{\bm U}(z_{0,i})\big)\Big].
\end{equation}
At an interior minimizer $\bm U^{\rm Po}_\star\succ\bm 0$ at which differentiation under the expectation is permitted,
\begin{equation}\label{orthocond}
    \E\Big[\sum_{i=1}^{n_0}\Big(1-\frac{p_{0,i}}{\ell^{\rm Po}_\star(z_{0,i})}\Big)\,\bm C(z_{0,i})\Big]=\bm 0,
\end{equation}
and the same holds on the boundary of the cone when $\epsilon_{\rm class}^{\rm Po}=0$ is attained, since then $\ell^{\rm Po}_\star=p$ a.s.\ and the integrand vanishes termwise.
\item\label{prop_target_comp} \emph{Composition.} Under the hypotheses of part~\ref{prop_target_proj}, for every $\bm U\in\Sy_+^m$, as an identity in $[0,\infty]$,
\begin{equation}\label{compsplit}
    \E\Big[\sum_{i=1}^{n_0}d_{\rm KL}\big(p_{0,i}\,\big\|\,\ell_{\bm U}(z_{0,i})\big)\Big]
    =\underbrace{\E\Big[\sum_{i=1}^{n_0}d_{\rm KL}\big(p_{0,i}\,\big\|\,\bar p(z_{0,i})\big)\Big]}_{\epsilon_{\rm comp}^{\rm Po}}
    +\int d_{\rm KL}\big(\bar p\,\big\|\,\ell_{\bm U}\big)\,d\mu_1 ,
\end{equation}
so $\epsilon_{\rm class}^{\rm Po}=\epsilon_{\rm comp}^{\rm Po}+\epsilon_{\rm dict}^{\rm Po}$ with the \emph{dictionary gap} $\epsilon_{\rm dict}^{\rm Po}\coloneqq\inf_{\bm U\in\Sy_+^m}\int d_{\rm KL}(\bar p\,\|\,\ell_{\bm U})\,d\mu_1$. The composition part $\epsilon_{\rm comp}^{\rm Po}\in[0,\bar n/\e]$ does not depend on the dictionary and vanishes if and only if $p_{0,i}=\bar p(z_{0,i})$ for all $i\le n_0$ $\Pa$-a.s., that is, if and only if the margins are functions of the own covariate. The population target is the $d_{\rm KL}$-projection of $\bar p$ onto the class in the geometry of $\mu_1$, whatever the joint law of the subsets, and its curve is unique $\mu_1$-a.e.\ whenever $\epsilon^{\rm Po}_{\rm class}<\infty$. If the dictionary is bounded as in (P1) of Section~\ref{sec_theory_po} and $\int\bm C\,d\mu_1\succ\bm 0$, a minimizer exists whenever $\epsilon^{\rm Po}_{\rm class}<\infty$.
\item\label{prop_target_score} \emph{Score and calibration equations.} Wherever $\ell_{t,i}>0$ at the selected items,
\[ \nabla_{\bm U}\Lcal^{\rm Po}(\bm U,\xi_t)=\sum_{i=1}^{n_t}\bm C_{t,i}-\sum_{i\in S_{t+1}}\frac{\bm C_{t,i}}{\ell_{t,i}}, \]
the quasi-score $\sum_{i\le n_t}w_{\rm Po}(\ell_{t,i})(\ell_{t,i}-y_{t,i})\,\bm C_{t,i}$ wherever every $\ell_{t,i}>0$. A minimizer $\hat{\bm U}$ of~\eqref{ermdef} over $\Sy_+^m$ with $\hat\ell_{t,i}>0$ at the selected items satisfies, for every $\bm v$ in its range and every $\bm v'\in\R^m$, including directions $\bm v'$ in its kernel,
\begin{equation}\label{calibeq}
   \frac1T\sum_{t<T}\sum_{i\le n_t}\Big(\frac{y_{t,i}}{\hat\ell_{t,i}}-1\Big)\,\big(\bm\phi(z_{t,i})\bm v\big)\big(\bm\phi(z_{t,i})\bm v'\big)=\lambda\,\bm v^\top\bm G\bm v',
\end{equation}
with $0/0\coloneqq0$. An item with $\hat\ell_{t,i}=0$ has $\bm\phi(z_{t,i})$ in the kernel of $\hat{\bm U}$, so it contributes nothing. Summing~\eqref{calibeq} with $\bm v'=\bm v$ over an eigenbasis of the range weighted by the eigenvalues of $\hat{\bm U}$ gives
\begin{equation}\label{aggcalemp}
   \frac1T\sum_{t<T}\sum_{i\le n_t}\big(y_{t,i}-\hat p_{t,i}\big)=\lambda\tr(\bm G\hat{\bm U}).
\end{equation}
\end{enumerate}
\end{proposition}

Identity~\eqref{compsplit} is the Pythagorean identity of a Bregman divergence for conditional means, see the proof of Theorem~1 of \citet{banerjeeguowang05} and Proposition~1 of \citet{banerjeeetal05jmlr}, here with the item intensity $\mu_1$ as measure and $\bar p$ as the conditional mean. The same argument gives the exact split $\E[\sum_{i\le n_0}(p_{0,i}-\ell(z_{0,i}))^2]=\E[\sum_{i\le n_0}(p_{0,i}-\bar p(z_{0,i}))^2]+\int(\bar p-\ell)^2\,d\mu_1$ for every bounded Borel $\ell$, which is the split used in Section~\ref{sec_sim_rate}.

The composition part of the gap is a dictionary-free floor. It is positive as soon as items with the same covariate carry different margins in differently composed cross-sections, as when the number of events per period is fixed on a panel of varying size. The dictionary part is the projection error of $\bar p$, which a richer dictionary drives to zero. For $\Zcal$ compact, $\bar p$ continuous and bounded away from zero, and a sequence of dictionaries whose spans are dense in $C(\Zcal)$, approximating $\sqrt{\bar p}$ uniformly by $\bm\phi\bm r$ and squaring gives $\ell_{\bm r\bm r^\top}\to\bar p$ uniformly. The approximating matrix has rank one, so interiority of $\bm U^{\rm Po}_\star$ is a property to check rather than assume, hence condition~\eqref{interior}. A richer dictionary relieves a capacity shortage but not a structural gap, which makes the floor testable (Section~\ref{sec_sim_rate}). No Pythagorean identity holds for the projection of $\bar p$ onto the cone, so~\eqref{compsplit} is the only exact split of the gap.

The calibration equations~\eqref{calibeq} are the in-sample face of part~\ref{prop_target_proj} and the in-sample part of the calibration that requirement (R1) asks for. They weight the residuals by $1/\hat p$, so heuristically the fit is calibrated tightly where inclusion is rare and loosely where it is likely, whereas the logistic regression, whose link is canonical \citep{mccullaghnelder89}, has unit weights. By~\eqref{aggcalemp} the fit reproduces the number of events without imposing it, and the penalty pushes it below the realized frequencies on its face, by $\lambda\bm v^\top\bm G\bm v$ per direction. The cone bounds the fitted probability from below but not from above. Since the unconstrained target is the margin, values above one arise only through misspecification or sampling error. The bound can be imposed convexly, by the linear constraints $\langle\bm C(z_k),\bm U\rangle+\eta_k\tr\bm U\le1$ at the centers $z_k$ of a finite covering of the covariate space, with $\eta_k$ bounding the variation of $\bm C$ in operator norm over the $k$th cell. They imply $\ell_{\bm U}\le1$ on $\Zcal$, since $\langle\bm A,\bm U\rangle\le\lambda_{\max}(\bm A)\tr\bm U$ for $\bm U\succeq\bm 0$. This is a tightening in the spirit of \citet{aubinfrankowskiszabo20}, who impose affine shape constraints on a whole set in this way for kernel machines. It can also be imposed by $\bm U\preceq B^{-2}\bm I_m$ under the bounded dictionary of Section~\ref{sec_theory_po}, at the price of a much smaller class, and Theorem~\ref{thm_main} then applies verbatim (part~\ref{thm_main_restr}), while the calibration identities acquire the multipliers of the added constraints. In the application, with an event rate of $0.25$ percent, the largest fitted value of any identity-link fit among the $3.0$ million test name-months is $0.26$.

\subsection{The log link and the Bernoulli likelihood}\label{sec_loglink}

The log link reads the sum of squares as the negative logarithm of the probability, $p_{\bm U}(z)\coloneqq\e^{-\ell_{\bm U}(z)}$, so that
\begin{equation}\label{loglinkmodel}
    0<p_{\bm U}\le1\quad\text{for every }\bm U\in\Sy_+^m,\qquad \ell_{\bm U}=-\log p_{\bm U},
\end{equation}
and fits it by the Bernoulli negative log-likelihood of the indicators, the \emph{log-link loss}
\begin{equation}\label{lossloglink}
    \Lcal^{\rm log}(\bm U,\xi_t)\coloneqq\sum_{i\in S_{t+1}}\ell_{t,i}-\sum_{i\notin S_{t+1}}\log\big(1-\e^{-\ell_{t,i}}\big),
\end{equation}
finite when $\ell_{t,i}>0$ at every non-event. With $\sum_{i\in S_{t+1}}\ell_{t,i}=\tr(\bm\Phi_t^\top\bm\Phi_t\bm U)-\sum_{i\notin S_{t+1}}\ell_{t,i}$ and $x+\log(1-\e^{-x})=\log(\e^x-1)$,
\[
   \Lcal^{\rm log}(\bm U,\xi_t)=\tr(\bm\Phi_t^\top\bm\Phi_t\bm U)-\sum_{i\notin S_{t+1}}\log\big(\e^{\ell_{t,i}}-1\big),
\]
the complementary log-log likelihood of the complement of the event. The model $p=\e^{-\ell}$ for the event is the model $1-p=1-\e^{-\ell}$ for its complement, so one solver serves both. Model~\eqref{loglinkmodel} is log-binomial (relative-risk) regression \citep{wacholder86} with a sum-of-squares index, and the standard difficulty of that model, fitted probabilities above one, the parametrization excludes \citep[on that difficulty see][]{marschnergillett12}. The price is the shape of the class. When the dictionary contains the constant and the products $\phi_k\phi_l$ are linearly independent, an affine function $a+\bm\phi\bm b$ of the dictionary lies in the class only if $b_k=0$ for every nonconstant $\phi_k$, so the log-linear probabilities of a logistic regression with rare events are approximated, not represented. With ${\rm expit}(x)\coloneqq1/(1+\e^{-x})$, $\e^{-\ell}=(1+\e^{-\ell})\,{\rm expit}(-\ell)$, so the log-link probability is the logistic probability with index $-\ell$ times $1+p\in(1,2]$, and for rare events the log-link fit is, to first order in $p$, a logistic regression whose index is minus a sum of squares. The dictionary supplies the curvature that a sum of squares requires, and the constant supplies the level. The probability saturates at one where $\ell$ vanishes, as the logistic probability saturates where the index is large. The log-link loss enters problem~\eqref{ermdef}, and its population target is a projection in the \emph{Bernoulli divergence}
\begin{equation}\label{dbedef}
    d_{\rm Be}(p\,\|\,q)\coloneqq p\log\frac{p}{q}+(1-p)\log\frac{1-p}{1-q}\ \ge\ 0,\qquad p,q\in(0,1),
\end{equation}
extended by $d_{\rm Be}(p\,\|\,1)\coloneqq+\infty$ for $p<1$, the ordinary Kullback--Leibler divergence between the two-point laws $(p,1-p)$ and $(q,1-q)$, which vanishes if and only if $q=p$ and is normalized by construction. The next proposition collects what the log-link fit shares with the identity-link fit. It is convex, calibrated to the margins on its own scale under any joint law, and blind to composition in the same way.

\begin{proposition}[The log-link fit]\label{prop_loglink}
\begin{enumerate}
\item\label{prop_loglink_convex} \emph{Convexity.} $\bm U\mapsto\Lcal^{\rm log}(\bm U,\xi_t)$ is convex on $\{\bm U\in\Sy^m:\ell_{t,i}>0\ \forall i\notin S_{t+1}\}$.
\item\label{prop_loglink_target} \emph{Population target.} Let $S\subseteq\{1,\dots,n\}$ be a random subset with margins $p_i\in(0,1)$ under any joint law. The log-link risk over unconstrained $\bm\ell\in(0,\infty)^n$,
\begin{equation}\label{rhologdef}
    \rho^{\rm log}(\bm\ell)\coloneqq\E\big[\Lcal^{\rm log}\big]=\sum_{i=1}^n\Big(p_i\ell_i-(1-p_i)\log\big(1-\e^{-\ell_i}\big)\Big),
\end{equation}
is strictly convex with unique minimizer $\ell_i=-\log p_i$, $i\le n$.
\item\label{prop_loglink_proj} \emph{Projection and composition.} Under the hypotheses of Proposition~\ref{prop_target}\ref{prop_target_proj} with moreover $p_{0,i}<1$ $\Pa$-a.s., for every $\bm U\in\Sy_+^m$, as an identity in $[0,\infty]$,
\begin{equation}\label{bklident}
    \E\big[\Lcal^{\rm log}(\bm U,\xi_0)\big]=\E\Big[\sum_{i=1}^{n_0}d_{\rm Be}\big(p_{0,i}\,\big\|\,\e^{-\ell_{\bm U}(z_{0,i})}\big)\Big]+c^{\rm log}_0,\qquad c^{\rm log}_0\coloneqq\E\Big[\sum_{i=1}^{n_0}h(p_{0,i})\Big],
\end{equation}
with $h(p)\coloneqq-p\log p-(1-p)\log(1-p)$ the binary entropy and $c^{\rm log}_0\in[0,\bar n\log 2]$ free of $\bm U$. The unregularized population target is the $d_{\rm Be}$-projection of the margins onto the class $\{\e^{-\ell_{\bm U}}:\bm U\in\Sy_+^m\}$, the class gap $\epsilon^{\rm log}_{\rm class}\coloneqq\inf_{\bm U\in\Sy_+^m}\E\big[\sum_{i\le n_0}d_{\rm Be}(p_{0,i}\|\e^{-\ell_{\bm U}(z_{0,i})})\big]$ splits as in~\eqref{compsplit}, with $d_{\rm Be}$ in place of $d_{\rm KL}$, into $\epsilon^{\rm log}_{\rm class}=\epsilon^{\rm log}_{\rm comp}+\epsilon^{\rm log}_{\rm dict}$ with
\begin{equation}\label{compsplitlog}
    \epsilon^{\rm log}_{\rm comp}\coloneqq\E\Big[\sum_{i=1}^{n_0}d_{\rm Be}\big(p_{0,i}\,\big\|\,\bar p(z_{0,i})\big)\Big],\qquad\epsilon^{\rm log}_{\rm dict}\coloneqq\inf_{\bm U\in\Sy_+^m}\int d_{\rm Be}\big(\bar p\,\big\|\,\e^{-\ell_{\bm U}}\big)\,d\mu_1,
\end{equation}
with $\epsilon^{\rm log}_{\rm comp}\in[0,\bar n\log2]$ free of the dictionary and zero if and only if the margins are functions of the own covariate, and the target is that of a panel with margins $\bar p(z_{t,i})$. If $\int\bar p\,\bm C\,d\mu_1\succ\bm 0$, a minimizer exists whenever $\epsilon^{\rm log}_{\rm class}<\infty$.
\item\label{prop_loglink_score} \emph{Score and calibration equations.} With the residual weight $w_{\rm log}(\ell)\coloneqq-1/(1-\e^{-\ell})$,
\[ \nabla_{\bm U}\Lcal^{\rm log}(\bm U,\xi_t)=\sum_{i=1}^{n_t}w_{\rm log}(\ell_{t,i})\big(\e^{-\ell_{t,i}}-y_{t,i}\big)\,\bm C_{t,i}, \]
and a minimizer $\hat{\bm U}$ of the regularized problem over $\Sy_+^m$ in the interior of the domain satisfies, with $\hat p_{t,i}=\e^{-\hat\ell_{t,i}}$, for every $\bm v$ in its range,
\begin{equation}\label{calibeqlog}
   \frac1T\sum_{t<T}\sum_{i\le n_t}\frac{\hat p_{t,i}-y_{t,i}}{1-\hat p_{t,i}}\,\big(\bm\phi(z_{t,i})\bm v\big)^2=\lambda\,\bm v^\top\bm G\bm v .
\end{equation}
At an event with $\ell_{t,i}=0$, where $\hat p_{t,i}=1$, the coefficient $w_{\rm log}(\ell_{t,i})(\e^{-\ell_{t,i}}-y_{t,i})$ and the quotient in~\eqref{calibeqlog} are read as their limits $1$ and $-1$. Equivalently, they are $y_{t,i}-(1-y_{t,i})\e^{-\ell_{t,i}}/(1-\e^{-\ell_{t,i}})$ and $(1-y_{t,i})\hat p_{t,i}/(1-\hat p_{t,i})-y_{t,i}$, which are defined there.
\end{enumerate}
\end{proposition}

Two differences from the identity-link fit matter below. First, the weight. Under the log link $|w_{\rm log}|=1/(1-\hat p)$ lies between one and $1.25$ for probabilities up to $0.2$. Against a common test function $\bm\phi\bm A\bm\phi^\top$ the log-link equations therefore weight the probability residual nearly uniformly and the identity-link equations by $1/\hat p$ (Proposition~\ref{prop_target}\ref{prop_target_score}). The difference comes from the classes. A direction $\bm A$ moves the probability by $\bm\phi\bm A\bm\phi^\top$ under the identity link and by $-\hat p\,\bm\phi\bm A\bm\phi^\top$ under the log link, while on the probability scale both losses weight squared errors by about $1/\hat p$ for rare events. Which link is calibrated where is therefore a question of the class gaps, which Section~\ref{sec_diag} measures. Second, the class. The two fits represent different curves, $p$ and $-\log p$, by a sum of squares in the same dictionary, and a margin curve that is a sum of squares does not in general have a sum-of-squares negative logarithm, so in general at most one of the two classes represents the margins. Both do only in special cases. The one that matters below is constant margins with a constant in the dictionary, where $c$ and $-\log c$ both lie in the class.

%%%%%%%%%%%%%%%%%%%%%%%%%%%%%%%%%%%%%%%%%%%%%%%%%%%%%%%%%%%%%%%%%%%%%%%%
\section{Statistical properties}\label{sec_theory}
% Question: can I trust it and do inference? Proofs and auxiliary lemmas are in Appendix app_proofs.

Section~\ref{sec_theory_po} treats the identity-link fit and Section~\ref{sec_theory_log} the log-link fit. The analysis works in the Euclidean space $(\Sy^m,\langle\cdot,\cdot\rangle)$ of~\eqref{cdef}, with the Frobenius norm $\|\cdot\|_F$ and dimension
\begin{equation}\label{Mdef}
    M\coloneqq\dim\Sy^m=m(m+1)/2 ,
\end{equation}
and states gradients, Hessians, and limit laws in this space without coordinates. We write $\opn{\cdot}$ for the operator norm of a linear map of $\Sy^m$ into itself, $\Ical$ for the identity map, $\bm A\otimes\bm\Delta$ for the rank-one map $\bm E\mapsto\langle\bm\Delta,\bm E\rangle\,\bm A$, and $\bar B(\bm U^{\rm Po}_\star,\delta)\coloneqq\{\bm U\in\Sy^m:\|\bm U-\bm U^{\rm Po}_\star\|_F\le\delta\}$. A linear map $\Xcal$ of $\Sy^m$ into itself is self-adjoint if $\langle\Xcal[\bm A],\bm\Delta\rangle=\langle\bm A,\Xcal[\bm\Delta]\rangle$ for all $\bm A,\bm\Delta$, and $\Xcal\succeq\nu\Ical$ means $\langle\Xcal[\bm A],\bm A\rangle\ge\nu\|\bm A\|_F^2$ for all $\bm A$. The gradients and Hessians of both losses are collected in Appendix~\ref{app_deriv}. Throughout, $\{\xi_t\}_{t\ge0}$ is strictly stationary on $\Xi$ with marginal law $\Pa$. Time-series dependence is the realistic regime for rolling panels, and i.i.d.\ sampling is the special case in which the dependence measures vanish. Recall the strong-mixing coefficients \citep[see][]{bra_05}
\begin{equation}\label{alphamixdef}
    \rmix(k)\coloneqq\sup\big\{|\Pa(A\cap B)-\Pa(A)\Pa(B)|:A\in\sigma(\xi_s,s\le t),\ B\in\sigma(\xi_s,s\ge t+k),\ t\ge0\big\},
\end{equation}
and that $\rmix(k)\to0$ implies ergodicity.

\subsection{The identity-link fit}\label{sec_theory_po}

We write $\Rcal(\bm U,\xi)\coloneqq\Lcal^{\rm Po}(\bm U,\xi)+\lambda\tr(\bm G\bm U)$ for the per-sample objective, so that $\Rcal_T$ in~\eqref{ermdef} is its average, and study any minimizer $\hat{\bm U}^{\rm Po}_T\in\argmin_{\bm U\in\Sy_+^m}\Rcal_T(\bm U)$, with fitted curve $\hat\ell^{\rm Po}(z)\coloneqq\langle\bm C(z),\hat{\bm U}^{\rm Po}_T\rangle$. The population risk is $\Ecal(\bm U)\coloneqq\E[\Rcal(\bm U,\xi_0)]$, written $\Ecal_\lambda$ where the penalty level matters, and $\bm U^{\rm Po}_{\star,\lambda}$ is a minimizer of $\Ecal_\lambda$ over $\Sy_+^m$, with curve $\ell^{\rm Po}_{\star,\lambda}$. At $\lambda=0$ it is the target $\bm U^{\rm Po}_\star$ of Section~\ref{sec_loss}. Theorem~\ref{thm_main} is stated at the penalty level of the fit, where $\bm U^{\rm Po}_\star$ stands for $\bm U^{\rm Po}_{\star,\lambda}$, the limit theory and the diagnostic of Section~\ref{sec_diag} take $\lambda=0$, and Corollary~\ref{cor_calib} compares the two. Wherever the margins are compared with the fitted curve, in Proposition~\ref{prop_target}\ref{prop_target_proj} and from Assumption~\ref{ass_clt} on, they are assumed positive, $p_{0,i}>0$ $\Pa$-a.s. If $\{\xi_t\}$ is ergodic and $\E[|\Rcal(\bm U,\xi_0)|]<\infty$ on an open convex set $\Ocal\subset\Sy^m$, then $\Ecal$ is finite and convex on $\Ocal$ and $\Rcal_T\to\Ecal$ uniformly on compacts of $\Ocal$, $\Pa$-a.s.\ (Lemma~\ref{lem_poprisk}). Minimization is over the cone and the local analysis lives on $\Ocal$, with no inclusion between the two. The loss is $+\infty$ at feasible points such as $\bm U=\bm 0$, while neighborhoods of a boundary target contain indefinite matrices at which it is finite. The only link is $\bm U^{\rm Po}_\star\in\Sy_+^m\cap\Ocal$.

The finite-sample analysis uses the following primitive conditions:
\begin{enumerate}
\item[(P1)] the dictionary is uniformly bounded, $B\coloneqq\sup_{z\in\Zcal}\|\bm\phi(z)\|_2<\infty$.
\item[(P2)] the cross-sectional and observed-subset sizes are bounded, $n_t\le\bar n$ and $|S_{t+1}|\le\bar s\le\bar n$, $\Pa$-a.s.
\item[(P3)] there is a constant $\underline c>0$ and a fixed open convex \emph{bounded} neighborhood $\Ocal$ of $\bm U^{\rm Po}_\star$ such that, $\Pa$-a.s.\ and for all $\bm U\in\Ocal$,
\begin{equation}\label{specbnd}
    \langle\bm C(z_{0,i}),\bm U\rangle\ \ge\ \underline c\qquad\text{for every }i\in S_{1}.
\end{equation}
\end{enumerate}
The bound~\eqref{specbnd} is required at the selected items only. With positive margins it bounds the curve away from zero at all items, and under exact representability it requires $p_{0,i}\ge\underline c$. On $\Ocal$ the matching upper bound $\ell_{\bm U}\le B^2\sup_{\bm U\in\Ocal}\|\bm U\|_F\eqqcolon\overline c$ is automatic. Interiority of the target is not needed. It suffices that the selected items keep their covariates away from the zero set of $\ell^{\rm Po}_\star$, so a singular target is compatible with (P1)--(P3). Define the constants
\begin{equation}\label{constants}
    \bar b\coloneqq B^2\bar n+\lambda\|\bm G\|_F,\qquad
    H\coloneqq\frac{\bar sB^4}{\underline c^{2}},\qquad
    M_\psi\coloneqq\frac{2B^2}{\underline c},
\end{equation}
which bound the linear coefficient of the loss and its Hessian on $\Ocal$ (Lemma~\ref{lem_vec}). Fix $r_0>0$ with $\bar B(\bm U^{\rm Po}_\star,r_0)\subseteq\Ocal$. We write $\bm\Psi^{\rm Po}_t\coloneqq\nabla\Rcal(\bm U^{\rm Po}_\star,\xi_t)-\nabla\Ecal(\bm U^{\rm Po}_\star)$ for the centered score, $\bar{\bm\Psi}^{\rm Po}_T\coloneqq\frac1T\sum_{t<T}\bm\Psi^{\rm Po}_t$, $\Psi_{t,j}\coloneqq\langle\bm E_j,\bm\Psi^{\rm Po}_t\rangle$ for its coordinates in a fixed orthonormal basis $\bm E_1,\dots,\bm E_M$ of $\Sy^m$, the only coordinates the paper uses, $\gamma(k)\coloneqq\E[\langle\bm\Psi^{\rm Po}_0,\bm\Psi^{\rm Po}_k\rangle]$, and
\begin{equation}\label{sigmaLRdef}
    \sigma^2\coloneqq\gamma(0),\qquad \sigma^2_{\rm LR}\coloneqq\gamma(0)+2\sum_{k\ge1}|\gamma(k)|.
\end{equation}
With $D^kf(\bm U)[\bm A_1,\dots,\bm A_k]$ the $k$-th directional derivative, under (P1)--(P3) the loss is \emph{pseudo-self-concordant} on $\Ocal$. $\Pa$-a.s., $\Rcal(\cdot,\xi)$ is three times continuously differentiable on $\Ocal$ and
\begin{equation}\label{pseudoSC}
    \big|D^3\Rcal(\bm U,\xi)[\bm A,\bm A,\bm\Delta]\big|\le M_\psi\,\|\bm\Delta\|_F\,D^2\Rcal(\bm U,\xi)[\bm A,\bm A]\qquad\text{for all }\bm U\in\Ocal,\ \bm A,\bm\Delta\in\Sy^m,
\end{equation}
and the same holds for $\Rcal_T$ and $\Ecal$ (Lemma~\ref{lem_SC}). The device is due to \citet{bach10} and was used by \citet{marteauferey19} and \citet{ost_bac_21} under independent sampling at a target where the gradient vanishes. Lemmas~\ref{lem_SC} and~\ref{lem_loc} adapt it to a linear penalty, a cone that the target may touch, and dependent samples.

\begin{theorem}[Consistency and finite-sample guarantees]\label{thm_main}
Let $\{\xi_t\}_{t\ge0}$ be strictly stationary and ergodic. Part~\ref{thm_main_cons} assumes the integrability hypothesis of Lemma~\ref{lem_poprisk} on an open convex set $\Ocal$ containing a point $\bm U^{\rm Po}_\star\in\Sy_+^m$. Parts~\ref{thm_main_det}--\ref{thm_main_mse} assume (P1)--(P3), $\bm U^{\rm Po}_\star\in\argmin_{\Sy_+^m}\Ecal$ at the penalty level of the fit, and
\begin{equation}\label{curv}
    \Hcal^{\rm Po}_\star\coloneqq\nabla^2\Ecal(\bm U^{\rm Po}_\star)\succeq\alpha\,\Ical\qquad\text{for some }\alpha>0.
\end{equation}
Let $\epsilon_0\coloneqq\frac\alpha6\min\{r_0,M_\psi^{-1}\}$, $\Acal_T^{\rm H}\coloneqq\{\opn{\nabla^2\Rcal_T(\bm U^{\rm Po}_\star)-\Hcal^{\rm Po}_\star}\le\alpha/2\}$ and $\Acal_T^{\rm s}\coloneqq\{\|\bar{\bm\Psi}^{\rm Po}_T\|_F<\epsilon_0\}$.
\begin{enumerate}
\item\label{thm_main_cons} \emph{Consistency.} If $\bm U^{\rm Po}_\star$ is the unique minimizer of $\Ecal$ over $\Sy_+^m\cap\Ocal$, then $\Pa$-a.s.\ $\Rcal_T$ has a minimizer over $\Sy_+^m$ for all $T$ large enough, and every sequence of minimizers converges to $\bm U^{\rm Po}_\star$. Minimality of $\bm U^{\rm Po}_\star$ over the whole cone is not assumed. Convexity of $\Rcal_T$ confines its minimizers over $\Sy_+^m$ to a ball on which $\Rcal_T$ is uniformly close to $\Ecal$.
\item\label{thm_main_det} \emph{Localization.} On $\Acal_T^{\rm H}\cap\Acal_T^{\rm s}$ the minimizer over $\Sy_+^m$ exists, is unique, and
\begin{equation}\label{detloc}
    \|\hat{\bm U}^{\rm Po}_T-\bm U^{\rm Po}_\star\|_F\le\frac6\alpha\|\bar{\bm\Psi}^{\rm Po}_T\|_F.
\end{equation}
\item\label{thm_main_conc} \emph{Concentration.} Let $F_\epsilon$ be the event that no minimizer exists, that it is not unique, or that it lies at distance at least $\epsilon$ from $\bm U^{\rm Po}_\star$, and $r_\epsilon\coloneqq\min(\alpha\epsilon/6,\epsilon_0)$. Then $F_\epsilon\subseteq\{\|\bar{\bm\Psi}^{\rm Po}_T\|_F\ge r_\epsilon\}\cup(\Acal_T^{\rm H})^{\rm c}$. Under i.i.d.\ sampling with $\E[\exp(\|\bm\Psi^{\rm Po}_0\|_F^2/\tau^2)]\le\e$,
\begin{equation}\label{conciid}
    \Pa[F_\epsilon]\le2\exp\Big(-\frac{Tr_\epsilon^2}{3\tau^2}\Big)+2M\exp\Big(-\frac{T\alpha^2}{32H^2+\tfrac83H\alpha}\Big).
\end{equation}
Under geometric mixing and bounded score coordinates a bound of the same form holds, with prefactors $M$ and $M^2$, factors $(\log T)^2$ in the exponents and an additional factor $M^2$ in the Hessian threshold from the entrywise argument (Lemma~\ref{lem_concmix}).
\item\label{thm_main_mse} \emph{Mean-squared error.} If $\sigma^2_{\rm LR}<\infty$, then $\E\big[\|\hat{\bm U}^{\rm Po}_T-\bm U^{\rm Po}_\star\|_F^2\,\mathbbm{1}_{\Acal_T^{\rm H}\cap\Acal_T^{\rm s}}\big]\le36\,\sigma^2_{\rm LR}/(\alpha^2T)$. With score coordinates bounded by $\bm b$, $|\gamma(k)|\le4\rmix(k)\|\bm b\|_2^2$, so summable mixing suffices.
\item\label{thm_main_curve} \emph{Fitted curve.} Under (P1), $\sup_{z\in\Zcal}|\hat\ell^{\rm Po}(z)-\ell^{\rm Po}_\star(z)|\le B^2\|\hat{\bm U}^{\rm Po}_T-\bm U^{\rm Po}_\star\|_F$, so parts~\ref{thm_main_cons}--\ref{thm_main_mse} hold for the curve, uniformly in $z$, with the factor $B^2$.
\item\label{thm_main_prim} \emph{Primitive conditions.} Under (P1)--(P3) the integrability hypothesis holds on $\Ocal$, $\Ecal$ is convex and lower semicontinuous on $\Sy^m$ with values in $(-\infty,+\infty]$, gradient and Hessian pass under the expectation on $\Ocal$, and
\begin{equation}\label{scorebnd}
    \|\bm\Psi^{\rm Po}_0\|_F\le2\Big(\frac{B^2\bar s}{\underline c}+\bar b\Big)\quad\Pa\text{-a.s.}
\end{equation}
Hence the sub-Gaussian condition holds, under i.i.d.\ sampling the first term of~\eqref{conciid} can be replaced by $2\exp(-Tr_\epsilon^2/(2b^2))$ with $b$ the right-hand side of~\eqref{scorebnd} \citep[Theorem~3.5]{pinelis94}, summable mixing gives $\sigma^2_{\rm LR}<\infty$, and the hypotheses of parts~\ref{thm_main_det}--\ref{thm_main_mse} reduce to~\eqref{curv} and $\bm U^{\rm Po}_\star\in\argmin_{\Sy_+^m}\Ecal$, a minimizer that~\eqref{curv} makes unique (Lemma~\ref{lem_loc}).
\item\label{thm_main_restr} All statements hold with $\Sy_+^m$ replaced by a closed convex subset of it and $\bm U^{\rm Po}_\star$ by the minimizer over that subset.
\end{enumerate}
\end{theorem}

Curvature enters only at the single point $\bm U^{\rm Po}_\star$, and pseudo-self-concordance converts it into localization, with no covering argument and no compactness of the feasible set. A single sample carries curvature only in the directions of its selected items, since its Hessian has rank at most $|S|$, so curvature accrues only by pooling across samples. Under exactly representable margins at an interior target, \eqref{curv} reads
\begin{equation}\label{curvexplicit}
    \E\Big[\sum_{i=1}^{n_0}\frac{\big(\bm\phi(z_{0,i})\bm A\bm\phi(z_{0,i})^\top\big)^2}{p_{0,i}}\Big]\ \ge\ \alpha\,\|\bm A\|_F^2\qquad\text{for all }\bm A\in\Sy^m,
\end{equation}
a checkable nondegeneracy of the quadratic forms generated by the dictionary. It fails when the products $\phi_k\phi_l$ are linearly dependent on the support of the covariates, and Appendix~\ref{app_ident} gives the limit theory of the curve in that case. Under i.i.d.\ sampling the guarantees become effective once $T\gtrsim\max\{\tau^2/\epsilon_0^2,H^2\log M/\alpha^2\}$. Under homogeneous margins $p_{0,i}\asymp p$, conditional independence and exchangeable items whose covariate law does not depend on $n_0$, $H/\alpha\asymp(\bar n/\E[n_0])B^4/(\alpha_1p)$ with $\alpha_1\coloneqq\lambda_{\min}(\E[\bm C_{0,1}\otimes\bm C_{0,1}])\le B^4/M$, so the cross-sectional size enters only through $\bar n/\E[n_0]$ and rare events cost a factor $p^{-2}$ in the threshold, under either link (Section~\ref{sec_theory_log}).

For the limit theory, define the score autocovariance operators and the long-run covariance operator on $\Sy^m$,
\begin{equation}\label{OmegaLRdef}
    \Gcal(k)[\bm A]\coloneqq\E\big[\langle\bm\Psi^{\rm Po}_k,\bm A\rangle\,\bm\Psi^{\rm Po}_0\big],\qquad
    \Qcal_{\rm LR}\coloneqq\Gcal(0)+\sum_{k\ge1}\big(\Gcal(k)+\Gcal(k)^*\big),
\end{equation}
with $\Gcal(k)^*$ the adjoint. Since $\gamma(k)=\tr\Gcal(k)$, the scalar $\sigma^2_{\rm LR}$ of~\eqref{sigmaLRdef} bounds $\tr\Qcal_{\rm LR}$.

\begin{assumption}\label{ass_clt}
Strict stationarity, (P1)--(P3) at $\lambda=0$, the curvature condition~\eqref{curv}, and
\begin{enumerate}
\item\label{ass_clt_i} \emph{score centering,} $\nabla\Ecal(\bm U^{\rm Po}_\star)=\bm 0$,
\item\label{ass_clt_ii} \emph{summable mixing,} $\sum_{k\ge1}\rmix(k)<\infty$.
\end{enumerate}
\end{assumption}

The first-order conditions over the cone read $\nabla\Ecal(\bm U^{\rm Po}_\star)=\bm\Lambda=\bm V\bm\Gamma\bm V^\top$ with $\bm\Gamma\in\Sy_+^{m-r}$, $r\coloneqq\rank\bm U^{\rm Po}_\star$, and the columns of $\bm V$ an orthonormal basis of $\ker\bm U^{\rm Po}_\star$, so part~\ref{ass_clt_i} is the statement $\bm\Gamma=\bm 0$. The radial direction lies in the range in any case, and $\langle\nabla\Ecal(\bm U^{\rm Po}_\star),\bm U^{\rm Po}_\star\rangle=0$ gives
\begin{equation}\label{aggcal}
    \E\Big[\sum_{i=1}^{n_0}\ell^{\rm Po}_\star(z_{0,i})\Big]=\E\Big[\sum_{i=1}^{n_0}p_{0,i}\Big]
\end{equation}
on the boundary of the cone as well as in its interior, so the fitted curve carries no aggregate bias. Part~\ref{ass_clt_i} holds at an interior target and, whatever the rank, under exact representability with positive margins, by~\eqref{orthocond}. Otherwise it generically holds only at an interior target, so interiority is a condition in its own right,
\begin{equation}\label{interior}
    \lambda_{\min}(\bm U^{\rm Po}_\star)>0,
\end{equation}
checked through $\lambda_{\min}(\hat{\bm U}^{\rm Po}_T)$ against the $T^{-1/2}$ sampling scale, the rule of Proposition~\ref{prop_plug}. A rank-deficient fit is evidence against~\eqref{interior}, while a full-rank fit does not establish it, since at a boundary target with zero multiplier the fit has full rank with positive limiting probability (Theorem~\ref{thm_clt}\ref{thm_clt_bdy}, Table~\ref{tab_e8thr}). At a boundary target the limit depends on the multiplier. With $\bm\Gamma=\bm 0$ it is the projection of a Gaussian onto the tangent cone
\begin{equation}\label{tangentcone}
    \Tcal(\bm U^{\rm Po}_\star)=\big\{\bm A\in\Sy^m:\ \bm V^\top\bm A\bm V\succeq\bm 0\big\},
\end{equation}
and under strict complementarity, $\bm\Gamma\succ\bm 0$, it is Gaussian on the tangent space of the rank-$r$ manifold, as for sample semidefinite programs \citep[Theorem~1]{shapiro19}. In the intermediate case the limit is the minimizer of the same quadratic, with the curvature term, over the critical cone \citep{shapiro89,bonnansshapiro00}, which we do not develop. Fix a compact convex $\Kcal_0\subseteq\Ocal$ with $\bm U^{\rm Po}_\star$ in its interior. The \emph{benchmark estimator} is any measurable $\hat{\bm U}^{\rm Po}_{\rm b}\in\argmin_{\bm U\in\Kcal_0}\Rcal_T(\bm U)$. It drops the cone constraint and coincides with the reported fit whenever the latter is positive definite and interior to $\Kcal_0$. In practice $\hat{\bm U}^{\rm Po}_{\rm b}$ is computed by Newton's method on $\Rcal_T$ started at the cone-constrained fit, with a backtracking line search that keeps the iterates in the domain of $\Rcal_T$. It coincides with the minimizer over $\Kcal_0$ with probability tending to one.

\begin{theorem}[Asymptotic distribution]\label{thm_clt}
Let Assumption~\ref{ass_clt} hold, with part~\ref{ass_clt_i} replaced by strict complementarity in part~\ref{thm_clt_sc}.
\begin{enumerate}
\item\label{thm_clt_par} \emph{Benchmark estimator.} The series~\eqref{OmegaLRdef} converges absolutely, $\hat{\bm U}^{\rm Po}_{\rm b}$ is unique and interior to $\Kcal_0$ with probability tending to one, $\nabla^2\Rcal_T(\hat{\bm U}^{\rm Po}_{\rm b})\to\Hcal^{\rm Po}_\star$ $\Pa$-a.s., and
\begin{equation}\label{cltdisplay}
    \sqrt T\,(\hat{\bm U}^{\rm Po}_{\rm b}-\bm U^{\rm Po}_\star)\ \longrightarrow_d\ \bm\eta,
\end{equation}
a centered, possibly degenerate Gaussian element of $\Sy^m$ with covariance operator $\Vcal\coloneqq{\Hcal^{\rm Po}_\star}^{-1}\Qcal_{\rm LR}{\Hcal^{\rm Po}_\star}^{-1}$.
\item\label{thm_clt_curve} \emph{Fitted curve.} With $\hat\ell^{\rm Po}_{\rm b}(z)\coloneqq\langle\bm C(z),\hat{\bm U}^{\rm Po}_{\rm b}\rangle$, $\sqrt T(\hat\ell^{\rm Po}_{\rm b}(z)-\ell^{\rm Po}_\star(z))\to_dN(0,\sigma^2(z))$, $\sigma^2(z)\coloneqq\langle\bm C(z),\Vcal[\bm C(z)]\rangle$, for every $z\in\Zcal$, and for continuous $\bm\phi$ in $C_b(\Zcal)$ to the finite-rank Gaussian process $\langle\bm C(\cdot),\bm\eta\rangle$.
\item\label{thm_clt_bdy} \emph{Boundary target, zero multiplier.} Without~\eqref{interior}, every measurable minimizer $\hat{\bm U}^{\rm Po}_T$ of $\Rcal_T$ over $\Sy_+^m$ at $\lambda=0$ satisfies
\begin{equation}\label{bdydisplay}
    \sqrt T\,(\hat{\bm U}^{\rm Po}_T-\bm U^{\rm Po}_\star)\ \longrightarrow_d\ \Pi_\star(\bm\eta)\coloneqq\argmin_{\bm A\in\Tcal(\bm U^{\rm Po}_\star)}\big\langle\bm A-\bm\eta,\Hcal^{\rm Po}_\star[\bm A-\bm\eta]\big\rangle,
\end{equation}
jointly with $\sqrt T(\hat\ell^{\rm Po}_T(z_j)-\ell^{\rm Po}_\star(z_j))\to_d\langle\bm C(z_j),\Pi_\star(\bm\eta)\rangle$ at finitely many points. If $r<m$ and $\Vcal$ is nondegenerate, the limit is not Gaussian and lies on the boundary of $\Tcal(\bm U^{\rm Po}_\star)$ with positive probability.
\item\label{thm_clt_sc} \emph{Boundary target, strict complementarity.} Let $r<m$ and $\nabla\Ecal(\bm U^{\rm Po}_\star)=\bm\Lambda=\bm V\bm\Gamma\bm V^\top$ with $\bm\Gamma\succ\bm 0$. Then $\rank\hat{\bm U}^{\rm Po}_T=r$ for all large $T$, $\Pa$-a.s., and
\begin{equation}\label{scdisplay}
    \sqrt T\,(\hat{\bm U}^{\rm Po}_T-\bm U^{\rm Po}_\star)\ \longrightarrow_d\ \bm\eta_\star,
\end{equation}
a centered Gaussian element of the tangent space $\Tcal_\star\coloneqq\{\bm A\in\Sy^m:\bm V^\top\bm A\bm V=\bm 0\}$ of the rank-$r$ manifold, with covariance operator
\[ \Vcal_\star\coloneqq\big(\Pcal(\Hcal^{\rm Po}_\star+\Scal_\star)\Pcal\big)^+\,\Pcal\,\Qcal_{\rm LR}\,\Pcal\,\big(\Pcal(\Hcal^{\rm Po}_\star+\Scal_\star)\Pcal\big)^+,\qquad
   \Scal_\star[\bm A]\coloneqq\bm\Lambda\bm A{\bm U^{\rm Po}_\star}^{+}+{\bm U^{\rm Po}_\star}^{+}\bm A\bm\Lambda, \]
$\Pcal$ the orthogonal projection onto $\Tcal_\star$ and $(\cdot)^+$ the Moore--Penrose inverse. The curve at $z$ is asymptotically $N(0,\sigma^2_\star(z))$, $\sigma^2_\star(z)\coloneqq\langle\bm C(z),\Vcal_\star[\bm C(z)]\rangle$.
\item\label{thm_clt_ci} \emph{Bands.} Let $\hat{\Hcal}_T\coloneqq\nabla^2\Rcal_T(\hat{\bm U}^{\rm Po}_{\rm b})$ and $\hat\Qcal_T$ be consistent for $\Qcal_{\rm LR}$, such as the empirical variance of the scores $\nabla\Rcal(\hat{\bm U}^{\rm Po}_{\rm b},\xi_t)$ under the condition of Proposition~\ref{prop_mds}, or in general a kernel estimator under the conditions of \citet[Theorem~1(a)]{andrews91} for estimated parameters, with bandwidth $h_T\to\infty$ and $h_T^2/T\to0$, which the Bartlett kernel meets under (P1)--(P3) and geometric mixing. Then
\[ \hat\ell^{\rm Po}_{\rm b}(z)\ \pm\ \zeta_{1-\beta/2}\,\sqrt{\big\langle\bm C(z),\hat{\Hcal}_T^{-1}\hat\Qcal_T\hat{\Hcal}_T^{-1}[\bm C(z)]\big\rangle/T} \]
is an asymptotic $(1-\beta)$ interval for $\ell^{\rm Po}_\star(z)$, with $\zeta_{1-\beta/2}$ the standard normal $(1-\beta/2)$-quantile, wherever $\sigma^2(z)>0$, which under~\eqref{csdef} is the margin of an item with covariate $z$. Under strict complementarity the interval centered at the cone-constrained fit $\hat\ell^{\rm Po}_T(z)$, with $\sigma^2_\star(z)$ in place of $\sigma^2(z)$ and the Hessian and score covariance compressed at that fit as in Remark~\ref{rem_faceband}, is an asymptotic $(1-\beta)$ interval for $\ell^{\rm Po}_\star(z)$. Under a zero multiplier it does not apply, and Proposition~\ref{prop_plug} gives a consistent band.
\end{enumerate}
\end{theorem}

Because the curve is linear in the parameter, part~\ref{thm_clt_curve} needs no delta-method step. The limit of part~\ref{thm_clt_bdy} can be simulated once the kernel of the target is estimated, and the kernel of the fit does not serve. With a zero multiplier the kernel block of $\sqrt T(\hat{\bm U}^{\rm Po}_T-\bm U^{\rm Po}_\star)$ converges to the positive semidefinite block $\bm V^\top\Pi_\star(\bm\eta)\bm V$, whose rank is random when $\Vcal$ is nondegenerate, which the assumptions allow to fail, so the rank of the fit then does not converge to $r$, and the tangent cone at the fit is not a consistent estimate of the tangent cone at the target. Proposition~\ref{prop_plug} in Appendix~\ref{app_projband} estimates the rank by a threshold that vanishes more slowly than the $T^{-1/2}$ scale of the kernel eigenvalues, the device of \citet[Section~6.4]{andrews1999b} (see also \citealp{fangsantos19} for estimated directional derivatives, and \citealp{donaldfortunapipiras07,chenfang19} for rank estimation), with the sampling scale of the smallest eigen-coordinate as its unit, and projects Gaussian draws onto the cone of the estimated kernel. Subsampling is the alternative \citep{politisromanowolf99}. Misspecified fits such as those of the application call for the bands of part~\ref{thm_clt_sc}, under which the rank of the fit identifies the kernel of the target, as does the rule of Proposition~\ref{prop_plug} in either regime. In the application strict complementarity is rarely supported, and we report no bands (Section~\ref{sec_app_data}).
Remark~\ref{rem_penclt} in Appendix~\ref{app_faceband} gives the limits under a penalty of order $T^{-1/2}$.

The long-run correction can be dropped under a condition weaker than correct specification of the joint law. Let $\{\Fcal_t\}_{t\ge0}$ be a filtration with $\sigma(\xi_0,\dots,\xi_{t-1},n_t,\bm z_t)\subseteq\Fcal_t$ and $S_{t+1}$ measurable with respect to $\Fcal_{t+1}$, and let
\begin{equation}\label{csdef}
    \text{(CS)}\qquad\text{there is }\bm U^\circ\in\Sy_+^m\text{ with }\Pa\big(i\in S_{t+1}\mid\Fcal_t\big)=\ell_{\bm U^\circ}(z_{t,i})>0\quad\Pa\text{-a.s.\ for all }i\le n_t.
\end{equation}
(CS) is the dynamic form of $\epsilon_{\rm class}^{\rm Po}=0$, dynamic completeness in the terminology of \citet{Wooldridgebook2001}, with the joint law of $S_{t+1}$ unrestricted beyond its margins. It fails whenever the margins depend on cross-sectional composition, by Proposition~\ref{prop_target}\ref{prop_target_comp}, which is why Theorems~\ref{thm_main} and~\ref{thm_clt} do without it. Under (CS), $\bm U^\circ$ minimizes the Poisson risk by Proposition~\ref{prop_target}\ref{prop_target_proj}, so Assumption~\ref{ass_clt}\ref{ass_clt_i} holds, and under~\eqref{curv} $\bm U^\circ=\bm U^{\rm Po}_\star$ and $\ell^{\rm Po}_\star(z_{t,i})=p_{t,i}$.

\begin{proposition}[Martingale scores under correct margins]\label{prop_mds}
Let Assumption~\ref{ass_clt} and~(CS) hold.
\begin{enumerate}
\item\label{prop_mds_i} The score $\bm\Psi^{\rm Po}_t=\sum_{i=1}^{n_t}\big(1-\mathbbm{1}_{\{i\in S_{t+1}\}}/p_{t,i}\big)\bm C_{t,i}$ is a martingale difference with respect to $\{\Fcal_{t+1}\}$, so $\Gcal(k)=\bm 0$ for $k\ge1$ and $\Qcal_{\rm LR}=\Gcal(0)$, whatever the temporal dependence of covariates and cross-sections.
\item\label{prop_mds_ii} With $\pi_{t,ij}\coloneqq\Pa(i,j\in S_{t+1}\mid\Fcal_t)$, $\pi_{t,ii}=p_{t,i}$,
\begin{equation}\label{Gamma0}
    \Gcal(0)[\bm A]=\E\Big[\sum_{i,j=1}^{n_0}\frac{\pi_{0,ij}-p_{0,i}p_{0,j}}{p_{0,i}\,p_{0,j}}\,\langle\bm C_{0,j},\bm A\rangle\,\bm C_{0,i}\Big]
    =\E\Big[\sum_{i=1}^{n_0}\frac{1-p_{0,i}}{p_{0,i}}\langle\bm C_{0,i},\bm A\rangle\bm C_{0,i}\Big]+\bm\Upsilon[\bm A],
\end{equation}
where the cross terms $\bm\Upsilon$ vanish under conditional independence.
\end{enumerate}
\end{proposition}

The margins alone govern the first-order behavior of the fit, and its variance~\eqref{Gamma0} reflects pairwise dependence only. The cross-sectional half is the working-independence approach of \citet{liangzeger86}, and part~\ref{prop_mds_i} adds that under (CS) the time-series half collapses too. The equality $\Qcal_{\rm LR}=\Gcal(0)$ is therefore the null of a specification check, and a long-run estimate materially different from the empirical score variance is evidence against dynamic completeness. The converse fails, since with samples independent across $t$ the autocovariances vanish whether or not (CS) holds. Section~\ref{sec_sim_cover} reports what ignoring the gap costs in simulation. Under (CS), ergodicity replaces summable mixing in Theorem~\ref{thm_clt} through the martingale central limit theorem of \citet{billingsley61}, and the score branch of Theorem~\ref{thm_main}\ref{thm_main_conc} holds without independence or mixing. At $\lambda=0$ the score is then a martingale difference bounded by $b_\Psi\coloneqq B^2\bar s/\underline c+B^2\bar n$, the bound of~\eqref{scorebnd} without the factor two for centering, so \citet[Theorem~3.5]{pinelis94} gives $\Pa[\|\bar{\bm\Psi}^{\rm Po}_T\|_F\ge r_\epsilon]\le2\exp(-Tr_\epsilon^2/(2b_\Psi^2))$. The Hessian branch still requires independence or mixing, since its summands are not martingale differences.

The preceding results combine into a decomposition of the distance between fitted and true margins, the population part of the calibration that requirement (R1) asks for.

\begin{corollary}[Calibration decomposition]\label{cor_calib}
Let the hypotheses of Theorem~\ref{thm_main} hold at the penalty level $\lambda\ge0$ of the fit, with target $\bm U^{\rm Po}_{\star,\lambda}$, unique by Lemma~\ref{lem_loc}.
\begin{enumerate}
\item\label{cor_calib_dict} \emph{Class gap.} If $\bm U^{\rm Po}_{\star,0}$ exists, as under part~\ref{cor_calib_lam}, and $\ell^{\rm Po}_{\star,0}\le\overline c$, then
\[ \E\Big[\sum_{i=1}^{n_0}\big(p_{0,i}-\ell^{\rm Po}_{\star,0}(z_{0,i})\big)^2\Big]\le2\max(1,\overline c)\,\epsilon_{\rm class}^{\rm Po}. \]
\item\label{cor_calib_lam} \emph{Regularization bias.} If $3\lambda\|\bm G\|_F/\alpha<\min\{r_0,M_\psi^{-1}\}$, then $\Ecal_0$ has a unique minimizer $\bm U^{\rm Po}_{\star,0}$ over $\Sy_+^m$, $\|\bm U^{\rm Po}_{\star,\lambda}-\bm U^{\rm Po}_{\star,0}\|_F\le3\lambda\|\bm G\|_F/\alpha$, and $\lambda\mapsto\tr(\bm G\bm U^{\rm Po}_{\star,\lambda})$ is nonincreasing.
\item\label{cor_calib_total} \emph{Total.} Let $\xi=(n,S,\bm z)$ be a fresh sample with law $\Pa$, independent of $\Dcal_T$, with margins $p_i\coloneqq\Pa(i\in S\mid n,\bm z)$, and $\|f\|\coloneqq(\E_\xi[\sum_{i\le n}f_i^2]/\E[n])^{1/2}$, the root mean square per item. On every sample on which $\hat{\bm U}^{\rm Po}_T$ exists,
\[ \big\|p_{\cdot}-\hat\ell^{\rm Po}(z_{\cdot})\big\|\le\underbrace{\big(2\max(1,\overline c)\,\epsilon_{\rm class}^{\rm Po}/\E[n]\big)^{1/2}}_{\text{class (composition and dictionary)}}+\underbrace{B^2\big\|\bm U^{\rm Po}_{\star,\lambda}-\bm U^{\rm Po}_{\star,0}\big\|_F}_{\text{regularization, }O(\lambda)}+\underbrace{B^2\big\|\hat{\bm U}^{\rm Po}_T-\bm U^{\rm Po}_{\star,\lambda}\big\|_F}_{\text{estimation, }O_\Pa(T^{-1/2})}. \]
The last term obeys Theorems~\ref{thm_main} and~\ref{thm_clt}, and the first does not vanish as the dictionary grows when the margins depend on cross-sectional composition (Proposition~\ref{prop_target}\ref{prop_target_comp}).
\end{enumerate}
\end{corollary}

A reliability curve \citep{murphywinkler77} plots, bin by bin of the fitted probability, the realized frequency of the event against the mean fitted probability.
% With bins fixed in advance as cells of the covariate space, the binned calibration errors in the Gaussian regimes of Theorem~\ref{thm_clt} would support Wald-type tests in the spirit of \citet{hosmerlemeshow80}, which we do not develop, while bins formed from the fitted probability have random edges and select on the estimate \citep{dimitriadisgneitingjordan21}.

\subsection{The log-link fit}\label{sec_theory_log}

The loss enters the results of Section~\ref{sec_theory_po} only through condition (P3), the constants $M_\psi$ and $H$ of~\eqref{constants}, and the score bound~\eqref{scorebnd}. The cone, the localization argument, the limit theory at a boundary target and the martingale structure do not depend on it. The log-link fit inherits them under the mirror image of (P3):
\begin{enumerate}
\item[(P3$^{\log}$)] there is a constant $\underline c>0$ and a fixed open convex bounded neighborhood $\Ocal$ of the log-link target $\bm U^{\log}_\star$, a minimizer over $\Sy_+^m$ of $\bm U\mapsto\E[\Lcal^{\log}(\bm U,\xi_0)]+\lambda\langle\bm G,\bm U\rangle$, such that, $\Pa$-a.s.\ and for all $\bm U\in\Ocal$, $\langle\bm C(z_{0,i}),\bm U\rangle\ge\underline c$ for every $i\le n_0$ with $i\notin S_1$.
\end{enumerate}
The bound is required at the non-events and keeps the fitted probabilities away from one, $p\le\e^{-\underline c}$ on $\Ocal$, the benign condition for a rare event. In the application no fitted probability of the log-link fit, in any dictionary, exceeds $0.57$. We write $\Rcal^{\log}$, $\Rcal^{\log}_T$, $\Ecal^{\log}$ and $\bm\Psi^{\log}_t$ for the per-sample objective, its average, the population risk and the centered score of Section~\ref{sec_theory_po} with $\Lcal^{\log}$ in place of $\Lcal^{\rm Po}$. Under (P1), (P2) and (P3$^{\log}$) the per-sample log-link objective $\Rcal^{\log}$ is pseudo-self-concordant with $M^{\log}_\psi\coloneqq B^2\coth(\underline c/2)$, its Hessian is bounded by $H^{\log}\coloneqq\bar nB^4/(4\sinh^2(\underline c/2))$, its linear coefficient by $\bar b^{\log}\coloneqq B^2\bar s+\lambda\|\bm G\|_F$, and its centered score by $2(B^2\bar n/\underline c+\bar b^{\log})$ (Lemma~\ref{lem_SC_log}). The roles of the events and the non-events are exchanged. The curvature of the log-link loss comes from the non-events, of which a cross-section has many, and that of the Poisson loss from the events, of which it may have few.

\begin{theorem}[The log-link fit]\label{thm_log}
Under (P1), (P2), (P3$^{\log}$) and the remaining hypotheses of each result, with $\bm U^{\log}_\star$ in place of $\bm U^{\rm Po}_\star$ and the constants of Lemma~\ref{lem_SC_log} in place of~\eqref{constants} and~\eqref{scorebnd}, Theorems~\ref{thm_main} and~\ref{thm_clt} and Proposition~\ref{prop_mds} hold for the log-link fit. In Proposition~\ref{prop_mds}, (CS) requires a $\bm U^\circ\in\Sy_+^m$ with $\Pa(i\in S_{t+1}\mid\Fcal_t)=\e^{-\ell_{\bm U^\circ}(z_{t,i})}\in(0,1)$, the score is that of Proposition~\ref{prop_loglink}\ref{prop_loglink_score}, and the weights of~\eqref{Gamma0} become $(\pi_{0,ij}-p_{0,i}p_{0,j})/((1-p_{0,i})(1-p_{0,j}))$. Corollary~\ref{cor_calib} holds with part~\ref{cor_calib_dict} replaced by
\[ \E\Big[\sum_{i=1}^{n_0}\big(p_{0,i}-\e^{-\ell^{\log}_{\star,0}(z_{0,i})}\big)^2\Big]\le\tfrac12\,\epsilon^{\log}_{\rm class}. \]
Since $x\mapsto\e^{-x}$ is $1$-Lipschitz on $[0,\infty)$ and the curve of the cone-constrained fit is nonnegative, the bounds on the curve hold for its fitted probability, whose limit law at $z$ is that of the curve times $-\e^{-\ell^{\log}_\star(z)}$. A minimizer exists whenever the objective is finite somewhere on the cone, as with a constant in the dictionary, and either $\lambda>0$ and $\bm G\succ\bm 0$ or $\lambda=0$ and the event design $\sum_{t<T}\sum_{i\in S_{t+1}}\bm C_{t,i}$ is positive definite, which holds eventually, $\Pa$-a.s., if $\E\big[\sum_{i\le n_0}p_{0,i}\bm C(z_{0,i})\big]\succ\bm 0$.
\end{theorem}

The condition for $\lambda=0$ is the analogue of the absence of separation in logistic regression, since a direction of the cone that no event sees lowers the loss of the non-events without bound. Three statements specific to the identity link change their form. The radial first-order condition gives $\E\big[\sum_{i\le n_0}\ell^{\log}_\star(z_{0,i})\,(p_{0,i}-p^{\log}_{\star,i})/(1-p^{\log}_{\star,i})\big]=0$ with $p^{\log}_{\star,i}\coloneqq\e^{-\ell^{\log}_\star(z_{0,i})}$ in place of~\eqref{aggcal}, so the log-link fit reproduces a weighted count of the events, not the count itself. Under homogeneous margins the Hessian branch of the threshold of Theorem~\ref{thm_main} carries no rare-event factor, since its curvature comes from the non-events, but the score branch still does. The log-link loss is a likelihood where the Poisson loss is a quasi-likelihood. If the items of a cross-section are independent given $\Fcal_t$ and (CS) holds in its log-link form, so that the score is a martingale difference by Proposition~\ref{prop_mds}, the expected Hessian equals the long-run variance of the score and the sandwich of Theorem~\ref{thm_clt} is the inverse information.

For the log link the penalty is itself a sample. With $\bm g_k$ the rows of $\bm G^{1/2}$, $\tr(\bm G\bm U)=\sum_k\bm g_k\bm U\bm g_k^\top$, and since the event term of~\eqref{lossloglink} is linear in $\bm U$, $T\,\Rcal^{\log}_T$ is the unpenalized loss with $\lambda T$ \emph{phantom events} at each $\bm g_k$. Bands and prediction call for different counts. For bands, Corollary~\ref{cor_calib}\ref{cor_calib_lam} and Theorem~\ref{thm_main}\ref{thm_main_mse} bound the regularization bias and the estimation error, and the ratio of the two bounds is $\lambda\|\bm G\|_F\sqrt T/(2\sigma_{\rm LR})$. Remark~\ref{rem_penclt} gives the exact shift of the limit law. When $\sigma^2_{\rm LR}$ grows in proportion to the events per period, this is of order $\lambda T/\sqrt N$, with $N$ the events in the window, and a fixed ratio calls for $\lambda\propto\sqrt N/T$. The ratio of bounds is a heuristic, and Table~\ref{tab_e14pen} confirms the scaling in simulation. For prediction, with only a constant and $N$ events among $n$ item-periods, the fit is $(N+\lambda T)/(n+\lambda T)$, and one phantom event minimizes the formal second-order expansion in $1/n$ of the expected excess log score.\footnote{With $N\sim{\rm Bin}(n,p)$ and $\hat p=(N+\lambda T)/(n+\lambda T)$, the exact expectation $\E[d_{\rm Be}(p\,\|\,\hat p)]$ is infinite, since $N=n$ has positive probability. The same holds for every dictionary, since a window without non-events has positive probability and then $\bm U=\bm 0$ minimizes the log-link objective, phantom events included. The excess risks of Table~\ref{tab_lamstudy} are Monte Carlo averages, which estimate the risk conditional on a window with at least one non-event. In its formal expansion in powers of $1/n$, the part of the coefficient of $n^{-2}$ that depends on $\lambda T$ is $-\lambda T+(\lambda T)^2(1-p)/(2p)-\lambda T(1-2p)/p$, which is minimized at $\lambda T=1$ for every $p$.} At a target of rank below $m$ the penalty also prunes directions that few events see, so the risk-optimal count rises with the dictionary in the simulations of Appendix~\ref{sec_sim_lampred}, though not in every dictionary of the application (Table~\ref{tab_lamstudy_app}).

The two links suit different ranges of probabilities. The constants of the identity-link fit grow as the fitted curve approaches zero at the events, and those of the log-link fit as the fitted probabilities approach one. For events of moderate or high probability the identity link is the natural choice, since its bands need no delta-method step and it reproduces the expected count exactly by~\eqref{aggcal}. For rare events the log link is, since its conditions, local constants and calibration weights stay benign where those of the identity link degrade, it is a likelihood rather than a quasi-likelihood, and to first order it is a logistic regression with index $-\ell_{\bm U}$ (Section~\ref{sec_loglink}). The application of Section~\ref{sec_app} is of the second kind. The failure rate is $0.249$ percent and the fitted probabilities span three orders of magnitude, while no log-link probability in any dictionary exceeds $0.57$, so the fitted values satisfy $\hat\ell\ge-\log0.57\approx0.56$ and the log-link constants are moderate. The identity-link fit is in the opposite regime, with calibration weights that vary by two orders of magnitude across the deciles of the fitted probability (Section~\ref{sec_app_calib}).

%%%%%%%%%%%%%%%%%%%%%%%%%%%%%%%%%%%%%%%%%%%%%%%%%%%%%%%%%%%%%%%%%%%%%%%%
\section{Diagnosing the model class}\label{sec_diag}
% Question: how do I know the class is wrong without knowing p? The two links of Section sec_estimator, compared.

The class gap $\epsilon^{\rm Po}_{\rm class}$ of Proposition~\ref{prop_target}\ref{prop_target_proj} is a population quantity defined through the margins, which are unknown. A richer dictionary tests it only against the composition floor (Section~\ref{sec_loss}), and a reliability curve sees it only through sampling noise (Corollary~\ref{cor_calib}). This section builds a diagnostic that needs no knowledge of $p$ from the two fits of Section~\ref{sec_estimator}, which represent different curves of the same margins, $p$ and $-\log p$. The two fits agree $\mu_1$-a.e.\ when both classes contain $\bar p$, and when they disagree, half the disagreement bounds from below the larger of their two errors on the probability scale, and through it the larger of the two dictionary gaps.

Each of the two fits estimates the margins up to the gap of its own class, so their disagreement is specification error observable without knowing $p$. Write $\ell^{\rm Po}_\star$ for the target of the Poisson loss, $\ell^{\rm log}_\star$ for that of the log link, $p^{\rm log}_\star\coloneqq\e^{-\ell^{\rm log}_\star}$ for its probability reading, and
\begin{equation}\label{Ddef}
    D(z)\coloneqq\ell^{\rm Po}_\star(z)-\e^{-\ell^{\rm log}_\star(z)},\qquad
    \hat D(z)\coloneqq\hat\ell^{\rm Po}(z)-\e^{-\hat\ell^{\rm log}(z)},
\end{equation}
for the population discrepancy and its empirical counterpart, both on the probability scale.

\begin{proposition}[The link discrepancy diagnostic]\label{prop_scalediag}
Let $\|f\|\coloneqq\big(\int f^2\,d\mu_1/\E[n_0]\big)^{1/2}$ for $f$ on $\Zcal$, the root mean square per item.
\begin{enumerate}
\item\label{prop_scalediag_bd} \emph{Lower bound on the dictionary gaps.} Under the hypotheses of Proposition~\ref{prop_loglink}\ref{prop_loglink_proj}, let both unregularized targets exist and $\ell^{\rm Po}_\star$ be bounded, as under (P1). Then
\begin{equation}
 \begin{split}
 \tfrac12\|D\|\ &\le\ \max\big(\|\ell^{\rm Po}_\star-\bar p\|,\ \|p^{\rm log}_\star-\bar p\|\big) \\
 &\le\ \max\Big(\sqrt{2\max\big(1,\|\ell^{\rm Po}_\star\|_\infty\big)\,\epsilon^{\rm Po}_{\rm dict}/\E[n_0]},\ \sqrt{\epsilon^{\rm log}_{\rm dict}/(2\E[n_0])}\Big).
 \end{split}
 \end{equation}
\item\label{prop_scalediag_null} \emph{The degenerate null.} If both classes contain $\bar p$, in particular if the margins are functions of the own covariate that both classes contain, then $D=0$ $\mu_1$-a.e. At $p\equiv c\in(0,1)$ with a constant in the dictionary the two scores are collinear and the variance $\sigma^2_{D,{\rm LR}}$ of Corollary~\ref{cor_scalediag} vanishes identically. Under the hypotheses of Corollary~\ref{cor_scalediag}\ref{prop_scalediag_alt} for both blocks, $\hat D=O_\Pa(T^{-1})$ for the benchmark fits.
\end{enumerate}
\end{proposition}

\begin{corollary}[Limit of the discrepancy]\label{cor_scalediag}
\begin{enumerate}
\item\label{prop_scalediag_alt} \emph{Limit distribution.} Stack the two fits as $(\bm U^{\rm Po},\bm U^{\rm log})$, and call each component, with its score and Hessian, a block. Let Assumption~\ref{ass_clt} hold for the identity-link block and its analogue for the log-link block at a minimizer $\bm U^{\rm log}_\star$ of its population risk, with (P3$^{\log}$) in place of (P3) and score centering at $\bm U^{\rm log}_\star$, with Hessian $\Hcal^{\rm log}_\star\succeq\alpha^{\rm log}\Ical$ and $\ell_{\bm U}\ge\underline c$ at every unselected item near $\bm U^{\rm log}_\star$. With the block-diagonal operator $\Hcal_{\rm j}\coloneqq\Hcal^{\rm Po}_\star\oplus\Hcal^{\rm log}_\star$ on $\Sy^m\times\Sy^m$, $\bm d(z)\coloneqq\big(\bm C(z),\,p^{\rm log}_\star(z)\,\bm C(z)\big)$ and $\Qcal_{{\rm j},{\rm LR}}$ the long-run covariance of the stacked centered score, the benchmark fits, and the cone-constrained fits at interior targets, satisfy
\[ \sqrt T\big(\hat D(z)-D(z)\big)\ \longrightarrow_d\ N\big(0,\sigma^2_{D,{\rm LR}}(z)\big),\qquad \sigma^2_{D,{\rm LR}}(z)=\big\langle\bm d(z),\Hcal_{\rm j}^{-1}\Qcal_{{\rm j},{\rm LR}}\Hcal_{\rm j}^{-1}[\bm d(z)]\big\rangle, \]
and $\hat D(z)\pm1.96\,\hat\sigma_{D,{\rm LR}}(z)/\sqrt T$, with the plug-in of Theorem~\ref{thm_clt}\ref{thm_clt_ci} for the stacked score, is an asymptotic $95\%$ interval for $D(z)$ wherever $\sigma_{D,{\rm LR}}(z)>0$.
\item\label{prop_scalediag_face} \emph{Boundary targets.} If the target of a block is singular and strict complementarity holds there as in Theorem~\ref{thm_clt}\ref{thm_clt_sc}, part~\ref{prop_scalediag_alt} holds for measurable minimizers of the two unregularized problems with $\Hcal_{\rm j}^{-1}$ replaced by $\mathcal W_{\rm j}\Pcal_{\rm j}$, where $\Pcal_{\rm j}$ and $\mathcal W_{\rm j}\coloneqq\big(\Pcal_{\rm j}(\Hcal_{\rm j}+\Scal_{\rm j})\Pcal_{\rm j}\big)^+$ are formed from the blockwise operators of that theorem, $\Ical$ and $0$ for a positive definite block, and with the plug-in of Remark~\ref{rem_faceband} for the interval. We write $\sigma^2_{D,\star}(z)$ for the resulting variance. If the multiplier of a singular block vanishes, neither interval is justified.
\end{enumerate}
\end{corollary}

By part~\ref{prop_scalediag_bd}, $\|D\|>0$ implies that at least one class fails to contain $\bar p$, and half the discrepancy bounds the larger probability-scale error, and through it the larger divergence, without knowledge of $p$. The bound concerns only the dictionary parts of the two gaps, since both targets, and with them $D$, are functionals of $(\mu_1,\bar p)$ by Propositions~\ref{prop_target}\ref{prop_target_comp} and~\ref{prop_loglink}\ref{prop_loglink_proj}, so $D$ carries no information about $\epsilon^{\rm Po}_{\rm comp}$ or $\epsilon^{\rm log}_{\rm comp}$, and it holds a fortiori for the class gaps $\epsilon^{\rm Po}_{\rm class}$ and $\epsilon^{\rm log}_{\rm class}$. The stacked score is not a martingale difference in general. Its autocovariances vanish when the samples are independent across $t$, which excludes, for instance, panels whose names exit after an event. Otherwise the long-run covariance of the stacked score requires a kernel estimator.

Two qualifications apply. First, the null $D=0$ $\mu_1$-a.e.\ is agreement of the two projections of $\bar p$. Joint representability, $p$ and $-\log p$ both in the sum-of-squares cone of the dictionary, implies it, but so does any misspecification orthogonal to the span of the products $\bm C$ under both weightings, to which the diagnostic is blind as it is to composition. A zero discrepancy is therefore no evidence of correct specification. At the one tractable null, constant margins with a constant in the dictionary, the limit of the benchmark fits is degenerate by part~\ref{prop_scalediag_null}, and the representing matrices $\Diag(c,0,\dots,0)$ and $\Diag(-\log c,0,\dots,0)$ have rank one, on the boundary of the cone. For the cone-constrained fits of~\eqref{Ddef} the multipliers vanish there, so Theorem~\ref{thm_clt}\ref{thm_clt_bdy} applies to both fits, and $\sqrt T\hat D(z)$ has the non-Gaussian limit $\langle\bm C(z),\Pi_\star(\bm\eta)+\Pi_\star(-\bm\eta)\rangle$, with $\bm\eta$ the limit of the identity-link benchmark. Its mean is orthogonal to the dictionary in $L^2(\mu_1)$ and so vanishes at some covariate values. The boundary effect of Section~\ref{sec_sim_diag} therefore shows in the norm of $\hat D$ (Table~\ref{tab_e5}), not pointwise. We therefore report no test of $D=0$ and read $\hat D$ as an estimate of the link discrepancy $D$ with the standard error of Corollary~\ref{cor_scalediag}\ref{prop_scalediag_alt}, and part~\ref{prop_scalediag_bd} turns it into a lower bound on the dictionary gaps. Second, the discrepancy diagnoses the model class, not the dependence structure, and it does not say which class fails. Neither loss needs independence for its own-scale calibration, so a nonzero discrepancy is no evidence against independence.

%%%%%%%%%%%%%%%%%%%%%%%%%%%%%%%%%%%%%%%%%%%%%%%%%%%%%%%%%%%%%%%%%%%%%%%%
\section{Simulation study}\label{sec_sim}

We check the rate and the composition floor, the bands and the diagnostic. Appendix~\ref{app_simapp} repeats the application with a known truth.

\subsection{Design}\label{sec_sim_design}

The dictionary has $m=4$ functions $\phi_k(z)=\cos(\omega_kz+\varsigma_k)$ with $\bm\omega=(1.5,4.0,8.5,14.0)$ and $\bm\varsigma=(0.3,1.1,2.0,0.7)$. Condition~\eqref{curv} holds under the independent data-generating process DGP-I below, since by~\eqref{curvexplicit} the Hessian dominates $\E[n_0]/\max p^\circ$ times the design operator $\E[\bm C(z)\otimes\bm C(z)]$ of an item under the baseline covariate law, whose smallest eigenvalue is $0.012$. A variant with $m=8$ adds the frequencies $2.7,6.2,11.0,17.5$ with phases $1.9,0.5,1.4,2.6$. Its products are linearly independent, but its design operator is singular to working precision, the smallest eigenvalue being below $10^{-16}$ of the largest, so in practice only its curve is identified. The true curve is $p^\circ=\bm\phi\bm U^\circ\bm\phi^\top$ with $\bm U^\circ=s(\bm W_0\bm W_0^\top+0.6\,\bm I_4)$, $\bm W_0$ a fixed standard normal draw and $s$ such that $p^\circ$ ranges over $[0.021,0.35]$. Panels have $n_t$ uniform on $\{15,\dots,45\}$ and $z_{t,i}=\tanh(a_t+\varepsilon_{t,i})$, $\varepsilon_{t,i}\sim N(0,0.7^2)$, with $a_t=0$ in the baseline design and $a_t=0.9\,a_{t-1}+N(0,0.4^2)$ in the dependent design. Under DGP-I inclusions are independent with probability $p^\circ(z_{t,i})$. Under DGP-D the subset is a determinantal point process \citep{macchi1975coincidence,kulesza2012determinantal} with $L$-ensemble kernel $\bm L_t=\Diag(\sqrt{\bm\varpi_t})\,\bm K^{\rm rbf}_t\,\Diag(\sqrt{\bm\varpi_t})$, where $\varpi_{t,i}=p^\circ(z_{t,i})/(1-p^\circ(z_{t,i}))$ and $[\bm K^{\rm rbf}_t]_{ij}=\exp(-(z_{t,i}-z_{t,j})^2/(2\cdot0.15^2))$, a locally repulsive law whose margins depend on the composition of the cross-section. In the null design inclusions are independent with $p\equiv c\in\{0.05,0.20\}$ on the dictionary $\{1,\phi_1,\dots,\phi_4\}$, so both the identity and the log class contain the margins, through rank-one matrices on the boundary of the cone, and we also fit the unconstrained benchmark estimators of Theorem~\ref{thm_clt}. The diagonal specification $\bm U=\Diag(\bm u)$, the additive model $\ell=\sum_ku_k\phi_k^2$, is a convex subcone to which Section~\ref{sec_theory} applies, and it is misspecified in every experiment in which it is fitted, with $47$ percent of the entrywise mass of $\bm U^\circ$ off the diagonal.

Fits are unregularized unless stated and use the solvers of the application (Section~\ref{sec_app_data}).\ifnum\joesubmission=1{} The code of the simulations and of the application was written and tested with the assistance of Claude (Anthropic) and reviewed by the author.\fi{} Accuracy is the root mean squared error (RMSE) of the fitted curve against the true margins on a fixed test panel of $200$ cross-sections. Population targets are computed on $20{,}000$ cross-sections with the exact margins in place of the inclusions, and $\bar p$ as the mean exact margin in $200$ covariate bins.

\subsection{Rate and the composition floor}\label{sec_sim_rate}

Table~\ref{tab_e1} reports the RMSE over $200$ replications of the baseline design. Under DGP-I the psd model attains the parametric rate, $\sqrt T\times$RMSE being flat near $0.18$ (Theorem~\ref{thm_main}). Under DGP-D its error settles near $0.021$, and doubling the dictionary lowers it only to $0.019$. Proposition~\ref{prop_target}\ref{prop_target_comp} splits this floor in population. The composition part, the root mean square of $p_{t,i}-\bar p(z_{t,i})$, is $0.0186$ whatever the dictionary. The dictionary part, that of $\bar p-\ell^{\rm Po}_\star$, falls from $0.0081$ at $m=4$ to $0.0015$ at $m=8$. The two add in squares to $0.0203$ and $0.0187$, the population totals up to the binning of $\bar p$, so composition, which no dictionary removes, makes up $84$ to $99$ percent of the squared floor. The diagonal specification shows the other kind of gap. Its RMSE is $0.049$ under DGP-I at every $T$, and under DGP-D it falls from $0.038$ to $0.033$ with $m=8$.

\begin{table}[t]
\centering
\caption{RMSE of the fitted curve against the true margins, $200$ replications, baseline design. DGP-I: independent, margins representable by the psd model. DGP-D: repulsive determinantal law. Monte Carlo standard errors are at most $0.0002$. The DGP-I columns use $m=4$.}
\label{tab_e1}
\begin{small}
\begin{tabular}{rccccccc}
\toprule
 & \multicolumn{4}{c}{psd model} & \multicolumn{3}{c}{diagonal specification} \\
\cmidrule(lr){2-5}\cmidrule(lr){6-8}
 & \multicolumn{2}{c}{DGP-I} & \multicolumn{2}{c}{DGP-D} & DGP-I & \multicolumn{2}{c}{DGP-D} \\
\cmidrule(lr){2-3}\cmidrule(lr){4-5}\cmidrule(lr){7-8}
$T$ & RMSE & $\sqrt T\times$RMSE & $m=4$ & $m=8$ & RMSE & $m=4$ & $m=8$ \\
\midrule
$250$  & $0.0110$ & $0.175$ & $0.0225$ & $0.0229$ & $0.0492$ & $0.0379$ & $0.0340$ \\
$1000$ & $0.0056$ & $0.176$ & $0.0210$ & $0.0198$ & $0.0487$ & $0.0375$ & $0.0333$ \\
$4000$ & $0.0029$ & $0.181$ & $0.0207$ & $0.0189$ & $0.0486$ & $0.0375$ & $0.0332$ \\
\bottomrule
\end{tabular}
\end{small}
\end{table}

\subsection{Coverage}\label{sec_sim_cover}

Table~\ref{tab_e2} reports the coverage of nominal $95$ percent pointwise intervals for $\ell^{\rm Po}_\star(z)$ at each of the points $z\in\{-0.5,0,0.5\}$ over $500$ replications of the dependent design. Under DGP-I the martingale variance of Proposition~\ref{prop_mds} covers at $0.93$ to $0.96$. Under DGP-D, where (CS) fails and the score is serially correlated, the fit is unbiased for its target to within $0.001$ at $T=1000$, yet at $z=0.5$ the martingale interval covers in only $0.81$ to $0.84$ of the replications. A Bartlett long-run variance recovers part of the loss, and the equal-weighted cosine estimator of \citet{lazarusetal18} most, $0.90$ to $0.93$. The intervals of the diagonal specification and those at $m=8$, with the Moore--Penrose inverse of the Hessian, cover at $0.93$ to $0.97$ at every $T$.

\begin{table}[t]
\centering
\caption{Coverage of nominal $95$ percent intervals for the target curve, $500$ replications, dependent design. Variances: martingale (Proposition~\ref{prop_mds}), Bartlett at the Newey--West lag $\lfloor4(T/100)^{2/9}\rfloor$ or at $\lceil1.3\,T^{1/2}\rceil$, which violates $h_T^2/T\to0$ of Theorem~\ref{thm_clt}\ref{thm_clt_ci} and serves as a robustness check, and equal-weighted cosine (EWC) with $\nu=\lfloor0.4\,T^{2/3}\rfloor$ and $t_\nu$ critical values \citep{lazarusetal18}.}
\label{tab_e2}
\begin{small}
\begin{tabular}{lrccc}
\toprule
 & $T$ & $z=-0.5$ & $z=0$ & $z=0.5$ \\
\midrule
DGP-I, martingale & $250$  & $0.932$ & $0.958$ & $0.950$ \\
 & $1000$ & $0.958$ & $0.954$ & $0.942$ \\
 & $4000$ & $0.948$ & $0.950$ & $0.938$ \\
\addlinespace
DGP-D, martingale & $250$  & $0.950$ & $0.926$ & $0.838$ \\
 & $1000$ & $0.966$ & $0.924$ & $0.806$ \\
 & $4000$ & $0.958$ & $0.930$ & $0.826$ \\
\addlinespace
DGP-D, Bartlett, Newey--West lag & $250$  & $0.936$ & $0.930$ & $0.862$ \\
 & $1000$ & $0.966$ & $0.928$ & $0.856$ \\
 & $4000$ & $0.954$ & $0.930$ & $0.872$ \\
\addlinespace
DGP-D, Bartlett, lag $1.3\,T^{1/2}$ & $250$  & $0.910$ & $0.910$ & $0.864$ \\
 & $1000$ & $0.952$ & $0.920$ & $0.878$ \\
 & $4000$ & $0.958$ & $0.930$ & $0.910$ \\
\addlinespace
DGP-D, EWC & $250$  & $0.948$ & $0.938$ & $0.912$ \\
 & $1000$ & $0.960$ & $0.936$ & $0.904$ \\
 & $4000$ & $0.962$ & $0.934$ & $0.928$ \\
\bottomrule
\end{tabular}
\end{small}
\end{table}

Appendix~\ref{app_faceband} checks the bands at targets of rank $4$ on $\{1,\phi_1,\dots,\phi_4\}$. Under strict complementarity the face band of Remark~\ref{rem_faceband} covers at $0.93$ to $0.96$ from $T=1000$ on, while the full-space band is $14$ to $15$ percent too wide, and $29$ percent at $T=250$. With a zero multiplier the rank of the fit does not settle, and the projection band of Proposition~\ref{prop_plug} covers at $0.89$ to $0.91$ at $T=250$ and $0.93$ to $0.95$ at $T=4000$ (Table~\ref{tab_e8thr}).

Fitting DGP-I at $\lambda=\kappa\,T^{-1/2}$, $\kappa\in\{0,1/4,1,4\}$, keeps the bias at a fixed fraction of the standard deviation as $T$ grows, at most $0.12$ for $\kappa\le1$ and $0.23$ for $\kappa=4$ (Corollary~\ref{cor_calib}\ref{cor_calib_lam}), and coverage and the fitted rank do not change.

\subsection{The diagnostic}\label{sec_sim_diag}

On the dictionary $\{1,\phi_1,\dots,\phi_4\}$ let $p^\circ=\e^{-\ell^\circ}$ with $\ell^\circ=c_0+s\,\bm\phi\bm U_4\bm\phi^\top$, $\bm\phi=(\phi_1,\dots,\phi_4)$ and $\bm U_4\coloneqq\bm W_1\bm W_1^\top+0.6\,\bm I_4$ for a second fixed standard normal $4\times4$ draw $\bm W_1$, a sum of squares and $p^\circ\in[0.02,0.30]$, so that model~\eqref{loglinkmodel} is correctly specified. In the dependent design with independent inclusions the log-link fit attains the parametric rate, $\sqrt T\times$RMSE being $0.25$ to $0.26$, and its martingale bands cover at $0.93$ to $0.96$ (Theorem~\ref{thm_log}). The identity-link fit, misspecified here, stalls at an RMSE of $0.033$.

Table~\ref{tab_e5} reports the discrepancy $\hat D$ of~\eqref{Ddef} at $T=1000$ in four designs, three with independent inclusions given the covariates and the repulsive DGP-D. In this log-correct design the identity class is wrong, in DGP-I the log class is wrong, in DGP-D both are, and in the null design both are correct. The last three columns evaluate Proposition~\ref{prop_scalediag}\ref{prop_scalediag_bd} in $L^2(\mu_1)$ on the $20{,}000$-period panel. Where one class is correct, $D$ is the error of the other class, so $\tfrac12\|D\|$ equals half the larger gap by construction. Under DGP-D, where both classes miss, it is $0.493$ of it. The divergence bound is $1.5$ to $2.4$ times the larger gap, so the first inequality carries the information. The discrepancy is estimated with a ratio of $\|\hat D\|$ to its standard error of $12$ to $14$ under DGP-I and DGP-D, and its stacked standard error matches the Monte Carlo dispersion there within ten percent in root mean square over the evaluation points, although both targets are singular in these designs and outside Corollary~\ref{cor_scalediag}\ref{prop_scalediag_alt}. In the log-correct design the identity-link target has rank three of five and satisfies strict complementarity, so the face standard error of Corollary~\ref{cor_scalediag}\ref{prop_scalediag_face} applies. It matches the Monte Carlo dispersion within five percent at each of the four points, and the ratio is $5.4$. The full-space standard error understates the dispersion by half there. At the null the cone-constrained fits settle on rank three around rank-one targets, and their discrepancy, $0.0035$ and $0.0063$, is a boundary effect. The benchmark fits, to which part~\ref{prop_scalediag_null} refers, give $\|\hat D\|$ of $0.0003$, with a standard error that overstates the dispersion by a third and a Wald test that rejects in $0.2$ and $0$ percent of the replications.

\begin{table}[t]
\centering
\caption{The link discrepancy $\hat D=\hat\ell^{\rm Po}-\e^{-\hat\ell^{\rm log}}$ of~\eqref{Ddef}, $T=1000$, baseline design, $200$ replications ($100$ under DGP-D, $500$ at the null), dictionary $\{1,\phi_1,\dots,\phi_4\}$. Root mean squares over the four points $z\in\{\pm0.25,\pm0.75\}$ or, for DGP-D, over a $60$-period test panel. ${\rm se}$: stacked-score standard error in the full space (Corollary~\ref{cor_scalediag}\ref{prop_scalediag_alt}). Ratio: $\|\hat D\|$ over the root mean square of the pointwise standard errors, a signal-to-noise summary and not a test statistic. The last three columns are the terms of Proposition~\ref{prop_scalediag}\ref{prop_scalediag_bd} in $L^2(\mu_1)$, per item, from the population targets on $20{,}000$ periods: half the discrepancy, the larger gap $\max(\|\ell^{\rm Po}_\star-\bar p\|,\|p^{\rm log}_\star-\bar p\|)$, and the divergence bound. At the null the rows are the benchmark fits of Theorem~\ref{thm_clt}\ref{thm_clt_par}. $^{\ast}$In the log-correct row the identity-link target is singular, and ${\rm se}$ is the face standard error of Corollary~\ref{cor_scalediag}\ref{prop_scalediag_face}. The full-space standard error there is $0.0027$. In the DGP-I and DGP-D rows both targets are singular, of rank four for the identity link and three for the log link, outside Corollary~\ref{cor_scalediag}\ref{prop_scalediag_alt}, and under DGP-I the identity-link multiplier vanishes, outside part~\ref{prop_scalediag_face} as well. In these rows ${\rm se}$ is the full-space standard error.}
\label{tab_e5}
\begin{small}
\setlength{\tabcolsep}{4pt}
\begin{tabular}{lccccccc}
\toprule
 & & \multicolumn{2}{c}{dispersion of $\hat D$} & & \multicolumn{3}{c}{part~\ref{prop_scalediag_bd} in population} \\
\cmidrule(lr){3-4}\cmidrule(lr){6-8}
design & $\|\hat D\|$ & ${\rm se}$ & MC s.d. & ratio & $\tfrac12\|D\|$ & larger gap & div.\ bound \\
\midrule
log-correct & $0.0297$ & $0.0055^{\ast}$ & $0.0054$ & $5.4$ & $0.0151$ & $0.0303$ & $0.0714$ \\
DGP-I (identity-correct) & $0.0456$ & $0.0037$ & $0.0040$ & $12.4$ & $0.0307$ & $0.0613$ & $0.0899$ \\
DGP-D & $0.0403$ & $0.0028$ & $0.0031$ & $14.3$ & $0.0208$ & $0.0422$ & $0.0730$ \\
null, $p\equiv0.05$ & $0.00028$ & $0.00041$ & $0.00031$ & $0.7$ & $0$ & $0$ & $0$ \\
null, $p\equiv0.20$ & $0.00025$ & $0.00039$ & $0.00028$ & $0.6$ & $0$ & $0$ & $0$ \\
\bottomrule
\end{tabular}
\end{small}
\end{table}

%%%%%%%%%%%%%%%%%%%%%%%%%%%%%%%%%%%%%%%%%%%%%%%%%%%%%%%%%%%%%%%%%%%%%%%%
\section{Application: failures of US common stocks}\label{sec_app}

We apply the estimator to monthly failures of US common stocks, an unbalanced panel whose cross-section changes every month, whose events are rare, and for which established competitors exist \citep{shumway01,chavajarrow04,DuffieSaitaWang2007,campbellhilscherszilagyi08}. All parameters are estimated out of sample.

\subsection{Data, event, and design}\label{sec_app_data}

The panel is the US common-stock universe of \citet{jensenkellypedersen23}, CRSP common stock on NYSE, AMEX and NASDAQ from January $1963$ to November $2024$. The covariate vector $z_{t,i}$ collects the $121$ of their $153$ characteristics that are missing in less than $30$ percent of the name-months, rank-transformed to $(-1,1)$ within month, so $d=121$ and the dictionary is bounded, condition (P1). This screen uses the full sample. A second screen drops the name-months that miss more than a third of the characteristics, $15$ percent of the CRSP name-months and $24$ percent of the performance delistings, whose median dropped name misses four fifths of the characteristics. The estimand therefore refers to names with at least two thirds of the characteristics observed. Characteristics at the end of month $t$ are paired with the event of month $t{+}1$. Missing coordinates enter through the masked dictionary~\eqref{maskeddict} (Appendix~\ref{app_masked}). The ordinary logit keeps the standard imputation at the median rank, and the CHS logit below, fitted on the complete name-months of each window, sets a missing predictor to its winsorized median in the window. Trained on every name-month with an indicator of each missing predictor, it gains $0.016$ ($t=-5.3$), and Tables~\ref{tab_app_main}, \ref{tab_app_chs} and~\ref{tab_app_robust} report this variant as a robustness row. On imputed covariates, with every missing characteristic at the median rank in all fits, the log score of the log-link fit of arc, $1.137$, is the same, and the results of every table of this section are available upon request. Five summaries enter beside the dictionary. Two describe the name, the within-month ranks, in $(-1,1)$, of the norm of $z_{t,i}$ with missing coordinates at the median rank and of its share of missing characteristics. Three describe the state of the market, which the rank transform removes, namely the share of names that lose more than $20$ percent or delist for cause in month $t$, its twelve-month mean, and the cross-sectional interquartile range of log returns, each mapped to its percentile over the expanding history of the panel.

Name $i$ enters $S_{t+1}$ if it is delisted for performance reasons in month $t{+}1$, CRSP delisting codes $500$ and $520$ to $599$, the delisting part of the failure indicator of \citet{campbellhilscherszilagyi08}. These are the performance codes of \citet{shumway97}, $500$ and $520$ to $584$, together with codes $585$, $587$ and $591$, which carry $2.7$ percent of the failures. Over the $646$ test months the failure rate is $0.249$ percent of $2{,}965{,}609$ name-months, drifting from $0.15$ percent over $1971$--$1987$ to $0.35$ over $1988$--$2006$ and $0.19$ after. The median month has seven failures.

The grid dictionaries, named for the grid of sizes $m$ below, use the kernel $k(z,z')=1+\langle z,z'\rangle/d$, so the fitted curve is a sum of squares of affine functions of $z$, on the information set of the classical distress logit. The dictionary is the pivoted Cholesky selection of $m\in\{5,10,20,40\}$ kernel sections at training covariates \citep{filipovicschneider24mom}, followed by the five summaries, so that the dictionary has $m+5$ functions. In this section $m$ counts the functions built from the covariates, without the summaries. For each of the $646$ test months from March $1971$ to December $2024$ the model is fitted on the preceding $97$ months, with $940$ failures at the median window, and evaluated on the unseen cross-section of about $4{,}600$ names. Both fits use the nonmonotone spectral projected gradient method of \citet{birginmartinezraydan00} and pass an ex post check of the first-order conditions. The empirical design operator is positive definite at $m\le20$, so these sizes are identified (Appendix~\ref{app_ident}). The identity-link fits have rank $2$ to $3$ on the grid and $4$ in the representation the solver returns for the richer dictionaries, and the log-link fits have numerical rank $6$ to $22$ on the grid, below $m+5$, which suggests that the targets lie on the boundary of the cone. Standard errors would then be face standard errors (Remark~\ref{rem_faceband}), which require strict complementarity at the target.

The competitors are the unconditional rate of the window, the ordinary logit on the $121$ characteristics and the five summaries, kernel logistic regression (KLR) \citep{zhuhastie05} in each dictionary, and the failure logit of \citet{campbellhilscherszilagyi08} (CHS) on their eight predictors, built from CRSP and Compustat as in their paper and winsorized at the $5$th and $95$th percentiles of the window. Each window uses the predictors observed in at least half of its name-months, all eight from the early $1980$s on.

The theory treats the dictionary as fixed. Here the pivots, the principal directions, the bandwidths, the completions of the masked dictionary and the percentiles of the state summaries are recomputed in each window from its covariates, which the theory does not cover. The event rate drifts across the test period. The three state summaries let the fitted curve move with the month, and Table~\ref{tab_e14} shows what a drifting rate costs the bands. Out of sample the log-link fit of arc predicts $71$, $91$ and $113$ percent of the realized failures over $1971$--$1987$, $1988$--$2006$ and after. The pattern follows the drift of the rate, which a $97$-month window tracks with a lag, but in sample the log-link fit reproduces a weighted count of the events rather than the count (Section~\ref{sec_theory_log}), and its penalty shifts the level, so the ratios do not separate drift from class error and regularization. Strict complementarity is generic for linear semidefinite programs \citep{alizadehhaeberlyoverton97}, and our programs are not linear. The data support it only weakly. On the grid the smallest eigenvalue of the estimated multiplier $\hat{\bm Q}_\perp^\top\nabla\Rcal_T(\hat{\bm U})\hat{\bm Q}_\perp$ of Remark~\ref{rem_faceband} exceeds two Bartlett standard errors in $4$ to $20$ percent of the windows for the identity-link fit and in at most $4$ percent for the log-link fit. Where the multiplier vanishes, Corollary~\ref{cor_scalediag}\ref{prop_scalediag_face} justifies no interval for the discrepancy, so we report no standard errors in the application. Condition (P3) is strained as well, since at $0.5$ to $0.6$ percent of the test events the identity-link grid fits lie below $10^{-4}$.

\subsection{Calibration and prediction}\label{sec_app_calib}

\begin{figure}[!tbp]
\centering
\figorbox[width=0.5\textwidth]{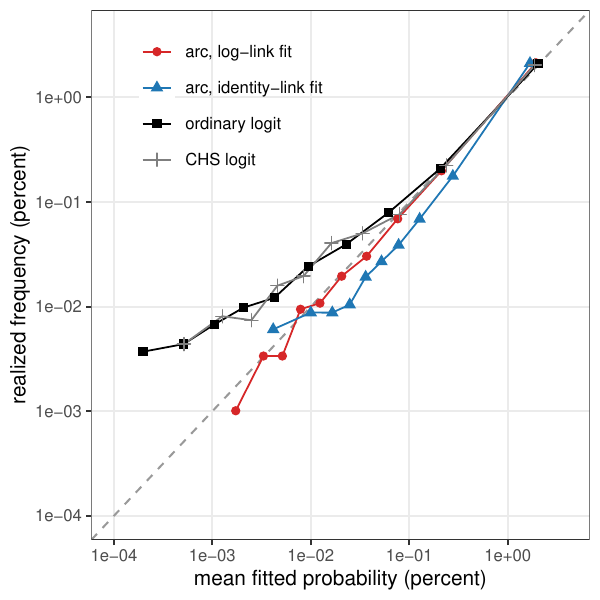}
\caption{Out-of-sample reliability of the fitted failure probabilities, $646$ test months pooled by bin, logarithmic axes. Realized frequency against mean fitted probability in ten equal-count bins per month, each model binned by its own fitted probability, for the log-link and the identity-link fit of the dictionary arc at $\lambda=T^{-1/2}$, the ordinary logit on all covariates and the CHS logit. Dashed: the diagonal.}
\label{fig_reliability}
\end{figure}

\begin{figure}[!tbp]
\centering
\figorbox[width=0.5\textwidth]{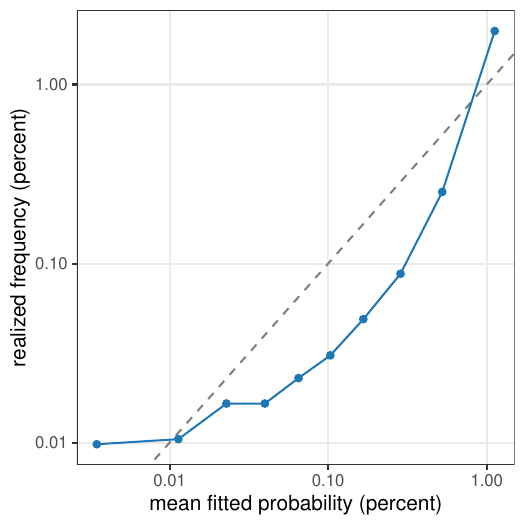}
\caption{Out-of-sample reliability of the identity-link grid fit at $m=10$, unregularized, $646$ test months pooled by bin, logarithmic axes. Realized frequency against mean fitted probability in ten equal-count bins per month by the fitted probability. Dashed: the diagonal.}
\label{fig_reliability_grid}
\end{figure}

Figure~\ref{fig_reliability} reports out-of-sample reliability curves for the dictionary arc of Section~\ref{sec_app_dicts}, with the test names split each month into deciles of the fitted probability and pooled.  The log-link fit stays within ten percent of the diagonal in the top two deciles. Its top decile, $1.94$ percent fitted against $2.14$ realized, dominates its calibration RMSE of $0.65\times10^{-3}$ in Table~\ref{tab_app_main}, which exceeds those of the ordinary logit and of KLR in arc, while the logistic competitors are far off in relative terms in their lowest decile, where failures occur about $19$ times as often as the ordinary logit predicts. The identity-link fit overstates the third to the ninth decile by a factor of $1.6$ to $2.4$ and understates the top decile, $1.69$ against $2.12$ percent, and the grid fit at $m=10$ does so more strongly (Figure~\ref{fig_reliability_grid}), by factors of $2.1$ to $3.4$ and with $1.12$ against $1.99$ percent in the top decile. Both projections tolerate larger absolute errors where the fitted probability is high, so these misses reflect the gap of the identity class (Section~\ref{sec_app_diag}). In the simulation of Appendix~\ref{app_simapp}, selection, the null set and a drifting rate keep the identity-link fit within eight percent of the realized frequency in every bin of its fit, apart from $22$ percent in the lowest bin at $m=40$. Bins of the fitted probability select on the estimate (Table~\ref{tab_e14}), so we attach no bands to them.

Table~\ref{tab_app_main} reports the average log score, its difference to the ordinary logit, the AUC, the area under the receiver operating characteristic curve, and the calibration RMSE. The log-link fit of arc has the best log score, $1.137$ against $1.185$ for the ordinary logit ($t=-11.5$, Newey--West with six lags on the monthly name-weighted differences, a test of equal predictive ability for forecasts from rolling windows of fixed length \citep{dieboldmariano95,giacominiwhite06}, applied to log scores as in \citealp{amisanogiacomini07}), about a tenth of the gain of the ordinary logit over the unconditional rate. KLR in arc, unpenalized, scores $1.163$. How the $0.048$ gain splits between the rank-angle dictionary and the sum-of-squares index depends on the clip, the dictionary's share being $63$ percent at a clip of $10^{-4}$ and $29$ percent at $10^{-10}$, and on how both sides are regularized, so on the $121$ characteristics we do not attribute it. On the eight CHS predictors the dictionary contributes almost all of it (Section~\ref{sec_app_chs}). Its mean monthly AUC is $0.940$ against $0.932$ ($t=3.8$, Newey--West with six lags on the monthly differences) and $0.914$ for the CHS logit ($t=8.5$), and its top percentile holds $39$ percent of all failures against $36$. With the equal-weighted cosine estimator of Table~\ref{tab_e2} the lead over the ordinary logit keeps $t=-7.8$, and the log-link fit of arc has the best score at every clip from $10^{-4}$ to $10^{-10}$ (Appendix~\ref{app_extra}). The identity-link fit of arc ties the ordinary logit in log score despite its calibration error, and the CHS logit is $0.015$ behind ($t=2.4$). With missingness indicators it ties the ordinary logit ($t=-0.3$), and the lead of the log-link fit of arc over it falls to $0.047$ ($t=-9.3$). On the grid both fits trail the logit on the same features, by $0.10$ to $0.11$ points for the identity-link fit and $0.02$ to $0.09$ for the log-link fit, since a sum of squares of affine functions cannot represent a monotone logistic or log-linear index (Section~\ref{sec_loglink}) and its unpenalized parameter count grows quadratically, $(m+5)(m+6)/2$, $1{,}035$ at $m=40$, against $m+6=46$ for the logit on the same features. The log-link fit worsens from $m=20$ to $m=40$ while the logit improves. Table~\ref{tab_app_robust} in Appendix~\ref{app_extra} shows that the lead of the log-link fit holds in both halves of the test period, at $0.051$ and $0.046$, and on the names with all CHS predictors observed, at $0.041$, where the CHS logit ties the ordinary logit. Sorted by the log-link fit of arc, the realized failure rate never falls across deciles and rises from $0.001$ to $2.14$ percent.

\begin{table}[!tbp]
\centering
\caption{Out-of-sample predictive performance on the failure event, $0.249$ percent of name-months, $646$ test months. Log score: average negative logarithmic score in units of $10^{-2}$, with every fitted probability clipped to $[10^{-6},1-10^{-6}]$. $\Delta$: difference of the log score to the ordinary logit, with the Newey--West $t$-statistic (six lags) of the monthly name-weighted differences in brackets. AUC: mean monthly AUC. RMSE: calibration RMSE across the ten reliability bins, pooled over months, each fit binned by its own fitted probability, in units of $10^{-3}$. Grid dictionaries: pivoted Cholesky selections of $m$ sections of the linear kernel, unregularized, against the logit with a linear index in the same $m+5$ features. Richer dictionaries: as described in Section~\ref{sec_app_dicts}, with the sum-of-squares fits at $\lambda=T^{-1/2}$ on their own norm against kernel logistic regression (KLR) in the same kernel, directions and summaries, unpenalized. The ordinary logit uses the $121$ characteristics and the five summaries, unpenalized, the CHS logit the eight predictors of \citet{campbellhilscherszilagyi08}, and its variant with missingness indicators is trained on every name-month with a missing predictor at its winsorized median and an indicator of each missing predictor. The unconditional rate of the training window scores $1.723$.}
\label{tab_app_main}
\begin{footnotesize}
\setlength{\tabcolsep}{2.5pt}
\resizebox{\textwidth}{!}{%
\begin{tabular}{lccccr@{\,}lccccr@{\,}lcc}
\toprule
 & \multicolumn{3}{c}{logit or KLR} & \multicolumn{5}{c}{identity link} & \multicolumn{5}{c}{log link} \\
\cmidrule(lr){2-4}\cmidrule(lr){5-9}\cmidrule(lr){10-14}
 & log score & AUC & RMSE & log score & \multicolumn{2}{c}{$\Delta$} & AUC & RMSE & log score & \multicolumn{2}{c}{$\Delta$} & AUC & RMSE \\
\midrule
ordinary logit & $1.185$ & $0.932$ & $0.23$ & & & & & & & & & & \\
CHS logit, eight predictors & $1.200$ & $0.914$ & $0.51$ & & & & & & & & & & \\
\quad with missingness indicators & $1.184$ & $0.919$ & $0.29$ & & & & & & & & & & \\
\addlinespace
\multicolumn{14}{l}{\emph{grid dictionaries: logit on the same features, identity-link and log-link fits}} \\
\quad $m=5$ & $1.275$ & & $0.51$ & $1.373$ & & & & $2.90$ & $1.295$ & & & & $0.64$ \\
\quad $m=10$ & $1.242$ & & $0.51$ & $1.353$ & & & & $2.98$ & $1.264$ & & & & $0.54$ \\
\quad $m=20$ & $1.226$ & & $0.42$ & $1.340$ & & & & $2.81$ & $1.266$ & & & & $0.42$ \\
\quad $m=40$ & $1.214$ & & $0.36$ & $1.326$ & & & & $2.49$ & $1.302$ & & & & $0.45$ \\
\addlinespace
\multicolumn{14}{l}{\emph{richer dictionaries: KLR, identity-link and log-link fits at $\lambda=T^{-1/2}$}} \\
\quad lin, $m=122$ & $1.185$ & $0.931$ & $0.21$ & $1.270$ & $+0.085$ & $(8.8)$ & $0.910$ & $1.97$ & $1.202$ & $+0.016$ & $(3.6)$ & $0.930$ & $0.87$ \\
\quad lg, $m=200$ & $1.180$ & $0.931$ & $0.19$ & $1.229$ & $+0.043$ & $(6.2)$ & $0.919$ & $1.69$ & $1.192$ & $+0.007$ & $(1.6)$ & $0.933$ & $0.82$ \\
\quad lgz, $m=200$ & $1.179$ & $0.933$ & $0.23$ & $1.212$ & $+0.026$ & $(4.5)$ & $0.918$ & $1.40$ & $1.174$ & $-0.011$ & $(-3.2)$ & $0.936$ & $0.84$ \\
\quad arc, $m=243$ & $1.163$ & $0.934$ & $0.15$ & $1.180$ & $-0.005$ & $(-1.0)$ & $0.926$ & $1.44$ & $1.137$ & $-0.048$ & $(-11.5)$ & $0.940$ & $0.65$ \\
\quad half, $m=243$ & $1.155$ & $0.934$ & $0.13$ & $1.184$ & $-0.002$ & $(-0.3)$ & $0.930$ & $1.82$ & $1.172$ & $-0.013$ & $(-3.5)$ & $0.933$ & $0.86$ \\
\bottomrule
\end{tabular}}
\end{footnotesize}
\end{table}

\subsection{Richer dictionaries}\label{sec_app_dicts}

Five richer dictionaries compare bases at a common penalty rule. \emph{lin} keeps the linear kernel $k(z,z')=1+\langle z,z'\rangle/d$ of the grid with all $121$ directions, $m=122$. Its dictionary spans the affine functions of $z$, so KLR in lin is the ordinary logit, except that missing characteristics enter at their conditional mean (Proposition~\ref{prop_masked}\ref{prop_masked_lin}) rather than at the median rank. The two score $1.185$ in Table~\ref{tab_app_main}. \emph{lg} adds a Gaussian term on the direction $u=z/\|z\|$, $k(z,z')=1+\langle z,z'\rangle/d+\exp(-\|u-u'\|^2/(2l^2))$ with $l$ set by the median heuristic, reduced to $m=200$ principal directions of its Nystr\"om features, and \emph{lgz} puts the Gaussian term on $z$. \emph{arc} and \emph{half} are the rank-angle dictionaries. Every rank dictionary is bounded, and these two are built from the rank angle. With $\theta_j\coloneqq\arccos z_j$, arc is the kernel
\[ k(z,z')=1+\frac1d\sum_{j=1}^d\cos(\theta_j-\theta'_j), \]
with the features $1$, $z_j/\sqrt d$ and $\sqrt{1-z_j^2}/\sqrt d$, so $m=2d+1=243$. Each term $b_jz_j+c_j\sqrt{1-z_j^2}$ of its span is a phase-shifted cosine of the rank angle, monotone, hump-shaped or U-shaped, and its sums of squares add pairwise interactions. Half uses $\cos((\theta_j-\theta'_j)/2)$, with the monotone features $\sqrt{(1\pm z_j)/2}$. In both, every affine function that is nonnegative on $[-1,1]^d$ is a sum of squares, so the identity link nests the linear probability model and the log link nests relative-risk regression on the whole box, which lin cannot. The summaries enter arc and half lifted like covariates, as two features each, so these dictionaries have $m+10$ functions.

The dictionaries are orthonormal in the RKHS of their kernel, so with $\bm G=\bm I$ the trace penalty is a sum of squared RKHS norms. We set $\lambda=T^{-1/2}$. This rule, the dictionaries and the summaries were chosen before the failure event was studied, on test-month scores of a related distress event, a monthly loss of more than $20$ percent or a delisting for cause, which contains every failure and is about $25$ times as frequent. No bandwidth was tuned. For the log-link fit this $\lambda$ places $\sqrt T=9.85$ phantom events in each direction (Section~\ref{sec_theory_log}), more than prediction or the bands call for. KLR and the ordinary logit penalized in the same way lose $0.04$ to $0.12$ points, so Table~\ref{tab_app_main} reports them unpenalized. Neither the sum-of-squares fits nor the logistic competitors are tuned per window. Table~\ref{tab_lamstudy_app} in Appendix~\ref{app_extra} reports, ex post, the log-link fits at other penalties. Two phantom events per direction improve every log-link fit, by $0.005$ to $0.031$ points, and leave the ranking of the dictionaries unchanged. With $800$ to $1{,}000$ face coordinates against $97$ months, sandwich bands are not reported for these dictionaries.

The rank-angle dictionaries are the best on this event. The log-link fit leads the ordinary logit in arc, half and lgz, ties it in lg and trails it in lin. The identity-link fit ties the ordinary logit in arc and half and trails it elsewhere. Unpenalized KLR is the strongest logistic model at every clip from $10^{-4}$ to $10^{-10}$, in half $0.008$ ahead of arc ($t=-1.6$). The log-link fit of arc leads KLR in half by $0.018$ ($t=-2.3$) and KLR in arc by $0.025$ ($t=-2.2$), with neither side tuned, though neither lead is significant with the EWC estimator ($t=-1.4$ and $-1.5$).

\subsection{The diagnostic on the panel}\label{sec_app_diag}

Table~\ref{tab_app_main} shows that the log-link fit beats the identity-link fit in log score in every dictionary, by $0.01$ to $0.09$ points. One property of the event favors it. Its fitted probabilities stay far from one, the largest being $0.57$, in arc, so the constants of Section~\ref{sec_theory_log} stay moderate. A logistic regression whose index is minus a sum of squares trails the log-link fit by only $0.003$ points (Table~\ref{tab_app_qlogit} in Appendix~\ref{app_extra}), so at this penalty the gain comes from the index class more than from the link.

Table~\ref{tab_app_diag} reports the link discrepancy~\eqref{Ddef} as point estimates. We attach no standard errors. Those of Corollary~\ref{cor_scalediag}\ref{prop_scalediag_face} assume strict complementarity in both fits, which the log-link fits rarely support (Section~\ref{sec_app_data}), and the bins are those of the identity-link fit, so they select on its estimate, which Table~\ref{tab_e14} shows can move coverage far from nominal. The bound of Proposition~\ref{prop_scalediag}\ref{prop_scalediag_bd} holds in any seminorm used on both sides. In the seminorm of the ten-bin RMSE, half the discrepancy is $0.33$ to $0.38$ of the calibration RMSE of the identity-link fit on the grid. In the same bins the log-link fit misses by $0.75$ to $1.09\times10^{-3}$, against $1.40$ to $2.98\times10^{-3}$ for the identity-link fit, so the identity-link error is the larger. In the sample $\tfrac12\|\hat D\|_{\rm bin}$ is at most the larger of the two binned RMSEs, by the triangle inequality. Read as a statement about the dictionary gaps, the ratio also needs small estimation, selection and drift parts, so it is indicative even on the grid. Sorted by the identity-link fit, $\hat D$ is positive in the middle deciles and negative in the top decile, the pattern of Figures~\ref{fig_reliability} and~\ref{fig_reliability_grid}. With the log-link fit close to calibrated in the upper deciles, the discrepancy is consistent with a class gap of the identity-link fit. If the log-link class were correct, $\|\hat D\|_{\rm bin}$ itself, not half of it, would estimate that gap, two thirds to three quarters of the calibration RMSE in the binned seminorm.

\begin{table}[!tbp]
\centering
\caption{The link discrepancy $\hat D=\hat\ell^{\rm Po}-\e^{-\hat\ell^{\rm log}}$ of~\eqref{Ddef} on the panel, failure event, psd model, $646$ test months. All entries are point estimates, and no standard errors are reported (Section~\ref{sec_app_diag}). The rows of the richer dictionaries use $\lambda=T^{-1/2}$, which biases the log-link fit (Appendix~\ref{sec_sim_app} and Table~\ref{tab_e14pen}), and are not estimates of the $D$ of Proposition~\ref{prop_scalediag}. $\|\hat D\|$: root mean square over all test names. $\|\hat D\|_{\rm bin}$: root mean square over the ten pooled bins of the identity-link fit of the bin means of $\hat D$, the seminorm of the calibration RMSE. RMSE: calibration RMSE of the identity-link fit (Table~\ref{tab_app_main}). RMSE log: calibration RMSE of the log-link fit in the same bins. Ratio: $\tfrac12\|\hat D\|_{\rm bin}$, the bound of Proposition~\ref{prop_scalediag}\ref{prop_scalediag_bd} in that seminorm, divided by the RMSE. On the richer dictionaries it is descriptive. Top: mean $\hat D$ in the top bin of the identity-link fit, pooled over months, in percentage points.}
\label{tab_app_diag}
\begin{small}
\setlength{\tabcolsep}{4pt}
\begin{tabular}{lcccccc}
\toprule
 & $\|\hat D\|$ & $\|\hat D\|_{\rm bin}$ & RMSE & RMSE log & ratio & top \\
\midrule
\quad $m=5$ & $0.0048$ & $0.0021$ & $0.0029$ & $0.0008$ & $0.37$ & $-0.60$ \\
\quad $m=10$ & $0.0055$ & $0.0023$ & $0.0030$ & $0.0008$ & $0.38$ & $-0.64$ \\
\quad $m=20$ & $0.0063$ & $0.0021$ & $0.0028$ & $0.0008$ & $0.38$ & $-0.59$ \\
\quad $m=40$ & $0.0071$ & $0.0016$ & $0.0025$ & $0.0010$ & $0.33$ & $-0.43$ \\
\addlinespace
\quad lin, $m=122$ & $0.0044$ & $0.0011$ & $0.0020$ & $0.0009$ & $0.28$ & $-0.31$ \\
\quad lg, $m=200$ & $0.0037$ & $0.0008$ & $0.0017$ & $0.0009$ & $0.24$ & $-0.23$ \\
\quad lgz, $m=200$ & $0.0039$ & $0.0005$ & $0.0014$ & $0.0009$ & $0.19$ & $-0.15$ \\
\quad arc, $m=243$ & $0.0055$ & $0.0007$ & $0.0014$ & $0.0008$ & $0.24$ & $-0.20$ \\
\quad half, $m=243$ & $0.0042$ & $0.0007$ & $0.0018$ & $0.0011$ & $0.20$ & $-0.22$ \\
\bottomrule
\end{tabular}
\end{small}
\end{table}

\subsection{The classical predictors}\label{sec_app_chs}

To compare the model classes at identical information, every model is refitted on the eight CHS predictors alone, rank-transformed as the characteristics are, with the five summaries, the two name-level ones computed from the eight predictors, and the same names, months and events, so lin has $m=9$ and arc and half have $m=17$. Table~\ref{tab_app_chs} reports the results. The CHS logit uses the levels of the predictors without the summaries. The logit on the ranked predictors with the summaries is $0.008$ ahead of it ($t=-1.6$). Without the summaries it is $0.022$ behind ($t=4.3$), so the rank transform loses information, and the summaries more than compensate for it. The log-link fit of arc on the eight predictors is the best model of the paper at the pre-specified penalty, at $1.128$, $0.072$ ahead of the CHS logit ($t=-11.0$) and $0.056$ ahead of its variant with missingness indicators ($t=-9.7$), $0.064$ ahead of the logit on the same ranked predictors and summaries and $0.009$ ahead of the same fit on the $121$ characteristics ($t=-3.0$). Most of the gain belongs to the dictionary. In arc every logistic model is within $0.004$ of the log-link fit, since the features $(z_j,\sqrt{1-z_j^2})$ let each predictor act nonmonotonically, and the sum of squares adds at most $0.004$. In lin the log-link fit trails the logit on the ranks by $0.017$, since an increasing log-linear index is not minus a sum of squares of affine functions. The identity-link fit trails the log-link fit in every dictionary, by $0.03$ to $0.10$ points.

\begin{table}[!t]
\centering
\caption{The failure event on the eight predictors of \citet{campbellhilscherszilagyi08}, $646$ test months. Log score as in Table~\ref{tab_app_main}, its difference to the CHS logit with the Newey--West $t$-statistic in brackets, mean monthly AUC, calibration RMSE in units of $10^{-3}$, and share of failures in the top percentile. The predictors enter as within-month ranks with the five summaries, and in the CHS logit as levels without them. Missing predictors enter the sum-of-squares fits and KLR through the masked dictionary and the other logits as in Table~\ref{tab_app_main}. Logit, $-$SOS index: ${\rm expit}(-\bm\phi\bm U\bm\phi^\top)$. KLR in $k^2$: ${\rm expit}(\bm\phi\bm A\bm\phi^\top)$ with a ridge on $\|\bm A\|_F$. All fits but KLR use $\lambda=T^{-1/2}$. The last two rows use the $121$ characteristics.}
\label{tab_app_chs}
\begin{footnotesize}
\setlength{\tabcolsep}{3.5pt}
\begin{tabular}{llcr@{\,}lccc}
\toprule
 & & log score & \multicolumn{2}{c}{vs CHS logit} & AUC & RMSE & top $1\%$ \\
\midrule
\multicolumn{2}{l}{CHS logit, levels} & $1.200$ & & & $0.914$ & $0.51$ & $0.358$ \\
\multicolumn{2}{l}{\quad with missingness indicators} & $1.184$ & $-0.016$ & $(-5.3)$ & $0.919$ & $0.29$ & $0.364$ \\
\multicolumn{2}{l}{logit on the ranks} & $1.192$ & $-0.008$ & $(-1.6)$ & $0.926$ & $0.34$ & $0.363$ \\
\addlinespace
lin, $m=9$ & identity link & $1.313$ & $+0.113$ & $(9.5)$ & $0.915$ & $2.72$ & $0.306$ \\
 & log link & $1.209$ & $+0.009$ & $(1.7)$ & $0.927$ & $0.78$ & $0.346$ \\
 & logit, $-$SOS index & $1.211$ & $+0.011$ & $(2.1)$ & $0.927$ & $0.86$ & $0.344$ \\
 & KLR in $k^2$ & $1.196$ & $-0.005$ & $(-0.7)$ & $0.925$ & $0.59$ & $0.373$ \\
 & KLR & $1.192$ & $-0.008$ & $(-1.6)$ & $0.927$ & $0.38$ & $0.366$ \\
\addlinespace
arc, $m=17$ & identity link & $1.176$ & $-0.024$ & $(-5.0)$ & $0.926$ & $1.08$ & $0.399$ \\
 & \textbf{log link} & $\bm{1.128}$ & $\bm{-0.072}$ & $(-11.0)$ & $0.934$ & $0.39$ & $0.414$ \\
 & logit, $-$SOS index & $1.129$ & $-0.071$ & $(-11.0)$ & $0.934$ & $0.46$ & $0.413$ \\
 & KLR in $k^2$ & $1.129$ & $-0.071$ & $(-10.7)$ & $0.934$ & $0.35$ & $0.416$ \\
 & KLR & $1.132$ & $-0.068$ & $(-10.6)$ & $0.934$ & $0.35$ & $0.413$ \\
\addlinespace
half, $m=17$ & identity link & $1.175$ & $-0.025$ & $(-5.3)$ & $0.928$ & $1.43$ & $0.410$ \\
 & log link & $1.147$ & $-0.053$ & $(-9.0)$ & $0.932$ & $0.57$ & $0.401$ \\
 & logit, $-$SOS index & $1.150$ & $-0.050$ & $(-8.9)$ & $0.932$ & $0.66$ & $0.400$ \\
 & KLR in $k^2$ & $1.142$ & $-0.058$ & $(-9.4)$ & $0.931$ & $0.50$ & $0.405$ \\
 & KLR & $1.130$ & $-0.070$ & $(-10.8)$ & $0.935$ & $0.37$ & $0.417$ \\
\addlinespace
\multicolumn{2}{l}{\emph{121 characteristics:} arc, log link} & $1.137$ & $-0.063$ & $(-10.9)$ & $0.940$ & $0.65$ & $0.392$ \\
\multicolumn{2}{l}{\emph{121 characteristics:} ordinary logit} & $1.185$ & $-0.015$ & $(-2.4)$ & $0.932$ & $0.23$ & $0.357$ \\
\bottomrule
\end{tabular}
\end{footnotesize}
\end{table}

%%%%
\section{Conclusion}\label{sec_conclusion}

Two limits are structural. Once the margins depend on the composition of the cross-section, every marginal fit carries a composition error that no dictionary removes, about two percentage points under a repulsive law in simulation. The floor is relative to the information set, so covariates that summarize the cross-section lower it, and a model of the subset is the alternative where no such summary is available. The model class may also miss the margins. Half the disagreement between the two links bounds the larger of their errors on the probability scale from below, and through it the larger of the two dictionary gaps, with a standard error and without knowledge of the margins, and it is blind to the composition part.

On monthly failures of US common stocks the log-link fit in a rank-angle dictionary predicts better in log score and in ranking than the ordinary logit on $121$ characteristics and the failure logit of \citet{campbellhilscherszilagyi08}, in log score by about a tenth of the logit's gain over the unconditional rate and significantly in every subsample, and it stays close to the diagonal of the reliability curve in relative terms, although the logistic competitors have the smaller calibration RMSE. On the eight predictors of that failure logit it leads it by $0.072$ points, most of it owed to the rank-angle dictionary. The identity-link fit ranks the names less well and is calibrated loosely from the third decile up, where a class gap of the identity link shows, and the link discrepancy is consistent with such a gap of the identity class. Next steps include bankruptcy filings and default ratings, longer horizons, richer summaries of the cross-section, a joint model of the observed subsets, and applications to occupancy data and to inclusion-probability weighting in survey and missing-data problems.

%%%%%%%%%%%%%%%%%%%%%%%%%%%%%%%%%%%%%%%%%%%%%%%%%%%%%%%%%%%%%%%%%%%%%%%%
\clearpage
\begin{appendices}

\section{Proofs}\label{app_proofs}

\subsection{Gradients and Hessians}\label{app_deriv}

Two exact identities are used throughout,
\begin{equation}\label{cbound}
    \|\bm C(z)\|_F=\|\bm\phi(z)\|_2^2,\qquad
    \sum_{i=1}^{n_t}\bm C_{t,i}=\bm\Phi_t^\top\bm\Phi_t .
\end{equation}
The first is why we work in $\|\cdot\|_F$, since a representation of $\Sy^m$ by its distinct entries is not an isometry and inflates the off-diagonal entries. The second identifies the linear term of the Poisson loss with the Gram matrix of the cross-section.

\begin{lemma}[The regularized risk and its derivatives]\label{lem_vec}
The per-sample regularized loss $\Rcal(\bm U,\xi)=\Lcal^{\rm Po}(\bm U,\xi)+\lambda\tr(\bm G\bm U)$ of Section~\ref{sec_theory_po} satisfies
\begin{equation}\label{Rvec}
    \Rcal(\bm U,\xi)=-\sum_{i\in S}\log\langle\bm C_i,\bm U\rangle+\langle\bm B(\xi),\bm U\rangle,\qquad
    \bm B(\xi)=\bm\Phi^\top\bm\Phi+\lambda\bm G\in\Sy_+^m,
\end{equation}
with $\bm C_i\coloneqq\bm C(z_i)$ and $\bm\Phi\coloneqq\bm\phi(\bm z)\in\R^{n\times m}$ the design and feature matrices of the generic sample $\xi=(n,S,\bm z)$, under the convention $-\log x=+\infty$ for $x\le0$. The function $\Rcal(\cdot,\xi)$ is convex and lower semicontinuous on $\Sy^m$, and twice continuously differentiable on the open convex set $\{\bm U\in\Sy^m:\langle\bm C_i,\bm U\rangle>0\ \forall i\in S\}$, where, with $\ell_i\coloneqq\langle\bm C_i,\bm U\rangle$ and $\bm W\coloneqq\Diag\big(\mathbbm{1}_{\{i\in S\}}/\ell_i\big)\in\R^{n\times n}$,
\begin{equation}\label{Rvecgrad}
    \nabla\Rcal(\bm U,\xi)=\bm\Phi^\top(\bm I_n-\bm W)\bm\Phi+\lambda\bm G,\qquad
    \nabla^2\Rcal(\bm U,\xi)=\sum_{i\in S}\frac{\bm C_i\otimes\bm C_i}{\ell_i^{2}}\ \succeq\ \bm 0,
\end{equation}
the gradient an element of $\Sy^m$ and the Hessian a self-adjoint positive semidefinite operator on $\Sy^m$, acting as $\nabla^2\Rcal[\bm A]=\sum_{i\in S}\langle\bm C_i,\bm A\rangle\bm C_i/\ell_i^2$. Moreover, $\xi\mapsto\Rcal(\bm U,\xi)$ is measurable for every $\bm U$.
\end{lemma}

\begin{proof}
The representation follows from~\eqref{cdef} and~\eqref{cbound} applied to each term of~\eqref{lossdef}. Each summand $-\log\langle\bm C_i,\bm U\rangle$, extended by $+\infty$, is convex and lower semicontinuous as the composition of $-\log$ with a linear functional, and the linear term is smooth. Differentiating $\bm U\mapsto\log\langle\bm C_i,\bm U\rangle$ gives $\bm C_i/\ell_i$, and $\sum_{i\in S}\bm C_i/\ell_i=\bm\Phi^\top\bm W\bm\Phi$, which with the second identity of~\eqref{cbound} yields the gradient. Differentiating once more gives the Hessian, self-adjoint and positive semidefinite because $\langle\nabla^2\Rcal[\bm A],\bm A\rangle=\sum_{i\in S}\langle\bm C_i,\bm A\rangle^2/\ell_i^2\ge0$. Measurability holds because the dictionary functions are Borel and $\Xi$ carries the $\sigma$-algebra of the disjoint union of the product Borel $\sigma$-algebras on $\{n\}\times2^{\{1,\dots,n\}}\times\Zcal^n$.
\end{proof}

By~\eqref{cbound}, $\|\bm B(\xi)\|_F\le\bar b$ and, under~\eqref{specbnd}, $\opn{\nabla^2\Rcal(\bm U,\xi)}\le H$ uniformly on $\Ocal$, $\Pa$-a.s., with the constants~\eqref{constants}. The Hessian has rank at most $|S|$ as an operator on the $M$-dimensional space $\Sy^m$, so per-sample strong convexity is unrealistic whenever the observed subsets are small relative to $M$. Since $\langle\bm C_i,\bm A\rangle=[\bm\Phi\bm A\bm\Phi^\top]_{ii}$, its action on a direction is
\begin{equation}\label{hessop}
    \nabla^2\Rcal(\bm U,\xi)[\bm A]=\bm\Phi^\top\Diag\Big(\tfrac{\mathbbm{1}_{\{i\in S\}}}{\ell_i^{2}}\,[\bm\Phi\bm A\bm\Phi^\top]_{ii}\Big)\bm\Phi ,
\end{equation}
computable in $O(nm^2)$ operations and $O(m^2)$ memory. The operator is therefore never assembled as an $M\times M$ array, which would cost $O(m^4)$ memory and $O(m^6)$ to factorize, and all second-order quantities of Section~\ref{sec_theory}, the sandwich variance of Theorem~\ref{thm_clt} included, can be evaluated by conjugate gradients through~\eqref{hessop}.

\begin{lemma}[The constants of the log-link loss]\label{lem_SC_log}
With $g(x)\coloneqq-\log(1-\e^{-x})$ for $x>0$ and $g(x)\coloneqq+\infty$ for $x\le0$, the per-sample regularized log-link loss is
\begin{equation}\label{Rlogvec}
\begin{split}
    \Rcal^{\log}(\bm U,\xi)&\coloneqq\Lcal^{\log}(\bm U,\xi)+\lambda\langle\bm G,\bm U\rangle=\langle\bm B^{\log}(\xi),\bm U\rangle+\sum_{i\notin S}g(\langle\bm C_i,\bm U\rangle),\\
    \bm B^{\log}(\xi)&\coloneqq\sum_{i\in S}\bm C_i+\lambda\bm G\in\Sy_+^m,
\end{split}
\end{equation}
convex and lower semicontinuous on $\Sy^m$, nonnegative on $\Sy_+^m$, measurable in $\xi$, and infinitely differentiable on $\{\bm U:\langle\bm C_i,\bm U\rangle>0\ \forall i\notin S\}$, with $\nabla\Rcal^{\log}=\bm B^{\log}(\xi)+\sum_{i\notin S}g'(\ell_i)\bm C_i$ and $\nabla^2\Rcal^{\log}=\sum_{i\notin S}g''(\ell_i)\,\bm C_i\otimes\bm C_i\succeq\bm 0$, the counterpart of Lemma~\ref{lem_vec}. Under (P1), (P2) and (P3$^{\log}$), $\Pa$-a.s., for all $\bm U\in\Ocal$ and $\bm A,\bm\Delta\in\Sy^m$:
\begin{enumerate}
\item\label{lem_SC_log_sc} $\big|D^3\Rcal^{\log}(\bm U,\xi)[\bm A,\bm A,\bm\Delta]\big|\le M^{\log}_\psi\,\|\bm\Delta\|_F\,D^2\Rcal^{\log}(\bm U,\xi)[\bm A,\bm A]$ with $M^{\log}_\psi\coloneqq B^2\coth(\underline c/2)\le B^2(1+2/\underline c)$;
\item\label{lem_SC_log_h} $\opn{\nabla^2\Rcal^{\log}(\bm U,\xi)}\le H^{\log}\coloneqq\bar nB^4\big/\big(4\sinh^2(\underline c/2)\big)\le\bar nB^4/\underline c^{2}$, and $\|\bm B^{\log}(\xi)\|_F\le\bar b^{\log}\coloneqq B^2\bar s+\lambda\|\bm G\|_F$;
\item\label{lem_SC_log_s} the centered score at the target satisfies $\|\bm\Psi^{\log}_0\|_F\le2\big(B^2\bar n/\underline c+\bar b^{\log}\big)$, and $|\Rcal^{\log}(\bm U,\xi)|\le\bar n\,g(\underline c)+\bar b^{\log}\sup_{\Ocal}\|\bm U\|_F$.
\end{enumerate}
\end{lemma}

\begin{proof}
Representation~\eqref{Rlogvec}: the selected summands of~\eqref{lossloglink} are $\langle\bm C_i,\bm U\rangle$ by~\eqref{cdef}. The extended $g$ is convex and lower semicontinuous on $\R$, since $g''>0$ on $(0,\infty)$ and $g(x)\to+\infty$ as $x\downarrow0$, so each non-event summand is convex and lower semicontinuous as the composition with a linear functional, and the linear term is smooth. On $\Sy_+^m$ every $\ell_i\ge0$, $g\ge0$ and $\langle\bm B^{\log}(\xi),\bm U\rangle\ge0$, so $\Rcal^{\log}\ge0$. The derivatives follow by the chain rule, and measurability as in Lemma~\ref{lem_vec}. With $\ell_i\coloneqq\langle\bm C_i,\bm U\rangle$, the selected summands and the penalty are linear in $\bm U$, so $D^2\Rcal^{\log}[\bm A,\bm A]=\sum_{i\notin S}g''(\ell_i)\langle\bm C_i,\bm A\rangle^2$ and $D^3\Rcal^{\log}[\bm A,\bm A,\bm\Delta]=\sum_{i\notin S}g'''(\ell_i)\langle\bm C_i,\bm A\rangle^2\langle\bm C_i,\bm\Delta\rangle$, with
\[ g'(x)=-\frac1{\e^x-1},\qquad g''(x)=\frac{\e^x}{(\e^x-1)^2}=\frac1{4\sinh^2(x/2)},\qquad g'''(x)=-\coth(x/2)\,g''(x). \]
On $\Ocal$, $\ell_i\ge\underline c$ at every non-event by (P3$^{\log}$), $\coth$ is decreasing on $(0,\infty)$, and $|\langle\bm C_i,\bm\Delta\rangle|\le B^2\|\bm\Delta\|_F$ by~\eqref{cbound}, which gives part~\ref{lem_SC_log_sc}. The inequality $\coth y=1+2/(\e^{2y}-1)\le1+1/y$ gives the bound on $M^{\log}_\psi$. Part~\ref{lem_SC_log_h}: $g''$ is decreasing, so $g''(\ell_i)\le g''(\underline c)=1/(4\sinh^2(\underline c/2))\le1/\underline c^{2}$ by $\sinh y\ge y$, and the rank-one operator $\bm C_i\otimes\bm C_i$ has norm $\|\bm C_i\|_F^2\le B^4$, for at most $\bar n$ non-events. The bound on $\bm B^{\log}$ is~\eqref{cbound} for at most $\bar s$ events. Part~\ref{lem_SC_log_s}: $\nabla\Rcal^{\log}=\bm B^{\log}(\xi)+\sum_{i\notin S}g'(\ell_i)\bm C_i$ with $|g'(\ell_i)|=1/(\e^{\ell_i}-1)\le1/\ell_i\le1/\underline c$, and centering at most doubles the bound. The bound on $|\Rcal^{\log}|$ follows from $0\le g(\ell_i)\le g(\underline c)$ at the non-events and Cauchy--Schwarz for the linear term.
\end{proof}

\subsection{Proofs for Section~\ref{sec_estimator}}\label{app_proofs_est}

\begin{proof}[Proof of Proposition~\ref{prop_target}]
\ref{prop_target_convex}: Each $\ell_{t,i}$ is linear in $\bm U$ by~\eqref{cdef}, and $x\mapsto x-a\log x$, extended by $+\infty$ on $x\le0$, is convex and lower semicontinuous on $\R$ for $a=1$ and linear for $a=0$.
\ref{prop_target_unc}: The sums are finite, so expectation and sum interchange, and $\E[\mathbbm{1}_{\{i\in S\}}]=p_i$ gives the closed form. The risk separates across $i$, and each summand $\ell\mapsto\ell-p\log\ell$ is strictly convex on $(0,\infty)$ for $p>0$, with derivative $1-p/\ell$ vanishing uniquely at $\ell=p$.
\ref{prop_target_proj}: By~\eqref{lossdef}, $\Lcal^{\rm Po}(\bm U,\xi_t)=\sum_{i=1}^{n_t}\big(\ell_{t,i}-\mathbbm{1}_{\{i\in S_{t+1}\}}\log\ell_{t,i}\big)$, and $\ell_{t,i}$ is measurable with respect to $(n_t,\bm z_t)$. On $\Sy_+^m$ each term is nonnegative, $\ell\ge0$ if $i\notin S_{t+1}$ and $\ell-\log\ell\ge1$ otherwise, so conditioning on $(n_t,\bm z_t)$ is legitimate in $[0,\infty]$ and gives $\E[\Lcal^{\rm Po}(\bm U,\xi_0)]=\E\big[\sum_i(\ell_{0,i}-p_{0,i}\log\ell_{0,i})\big]$, with $\ell-p\log\ell=d_{\rm KL}(p\|\ell)+(p-p\log p)$ termwise, including $\ell=0$ under the convention $d_{\rm KL}(p\|0)=+\infty$. Both summands are nonnegative, so taking expectations produces~\eqref{klident} in $[0,\infty]$. The constant lies in $[0,\bar n]$ because $p-p\log p\in(0,1]$ on $(0,1]$ and $n_t\le\bar n$. The projection statement follows from $d_{\rm KL}\ge0$, with equality only at $\ell=p$. For~\eqref{orthocond}, differentiate $\bm U\mapsto\E[\Lcal^{\rm Po}(\bm U,\xi_0)]$ under the integral sign at an interior point, where $\partial_{\bm U}\ell_{\bm U}(z)=\bm\phi(z)^\top\bm\phi(z)$. The interchange requires the gradient of $\bm U\mapsto\Lcal^{\rm Po}(\bm U,\xi_0)$ to be dominated by an integrable variable on a neighborhood of $\bm U^{\rm Po}_\star$. This is assumed here, and it holds under the primitive conditions of Section~\ref{sec_theory_po} by Theorem~\ref{thm_main}\ref{thm_main_prim}.
\ref{prop_target_comp}: This is the Bregman identity for conditional means, and we verify it in $[0,\infty]$. First, $\bar p>0$ $\mu_1$-a.e.\ by positivity of the margins, \eqref{pbarident} with $f=\mathbbm{1}_{\{\bar p=0\}}$, and the null-set statement that follows it, so $\bar p(z_{0,i})>0$ for all $i\le n_0$, $\Pa$-a.s. Second, $\epsilon^{\rm Po}_{\rm comp}\le\bar n/\e$. Indeed $d_{\rm KL}(p\|\bar p)=p\log p+p(-\log\bar p)+(\bar p-p)$, where $\E[\sum_ip_{0,i}\log p_{0,i}]\le0$, $\E[\sum_ip_{0,i}(-\log\bar p(z_{0,i}))]=\int(-\bar p\log\bar p)\,d\mu_1\le\bar n/\e$ by~\eqref{pbarident} with the nonnegative $f=-\log\bar p$, and $\E[\sum_i(\bar p(z_{0,i})-p_{0,i})]=0$ by~\eqref{pbarident} and the definition of $\mu_1$. If $\mu_1(\{\ell_{\bm U}=0\})>0$, both sides of~\eqref{compsplit} are infinite, since $d_{\rm KL}(p\|0)=+\infty$ for $p>0$ and $\bar p>0$ $\mu_1$-a.e. Otherwise write $X_i\coloneqq d_{\rm KL}(p_{0,i}\|\ell_{\bm U}(z_{0,i}))$ and $Y_i\coloneqq d_{\rm KL}(p_{0,i}\|\bar p(z_{0,i}))$, both nonnegative. Then $X_i-Y_i=p_{0,i}u(z_{0,i})+v(z_{0,i})$ with $u\coloneqq\log(\bar p/\ell_{\bm U})$ and $v\coloneqq\ell_{\bm U}-\bar p$, affine in the margin, and $d_{\rm KL}(\bar p\|\ell_{\bm U})=\bar pu+v$. On $A_k\coloneqq\{|u|\le k,\ |v|\le k\}$ the functions are bounded, so~\eqref{pbarident}, which extends to bounded signed $f$ because $\mu_1$ is finite, and the definition of $\mu_1$ give
\[ \E\Big[\sum_i(X_i-Y_i)\mathbbm{1}_{A_k}(z_{0,i})\Big]=\int_{A_k}d_{\rm KL}(\bar p\|\ell_{\bm U})\,d\mu_1, \]
both expectations on the left being finite since $X_i\mathbbm{1}_{A_k}\le Y_i+2k$. As $k\to\infty$, $A_k$ increases to a set of full $\mu_1$-measure, so $z_{0,i}\in\bigcup_kA_k$ for all $i$, $\Pa$-a.s., by the null-set statement. The right-hand side and $\E[\sum_iX_i\mathbbm{1}_{A_k}(z_{0,i})]$ increase to their unrestricted counterparts by monotone convergence, and $\E[\sum_iY_i\mathbbm{1}_{A_k}(z_{0,i})]\to\epsilon^{\rm Po}_{\rm comp}$ by dominated convergence, which gives~\eqref{compsplit} in $[0,\infty]$.
For the gap, take the infimum over $\bm U$ in~\eqref{compsplit}, where the composition constant does not involve $\bm U$. It is at most $\bar n/\e$ by the bound above, and a sum of nonnegative terms each of which vanishes only at $p_{0,i}=\bar p(z_{0,i})$. By~\eqref{klident} the population risk differs from the left-hand side of~\eqref{compsplit} by a constant, so its minimizers over the cone are those of $\int d_{\rm KL}(\bar p\|\ell_{\bm U})\,d\mu_1$, which involves the law only through $(\mu_1,\bar p)$. Since $d_{\rm KL}(\bar p\|\cdot)$ is strictly convex wherever $\bar p>0$, that is, $\mu_1$-a.e., the target curve is unique $\mu_1$-a.e.\ when $\epsilon^{\rm Po}_{\rm class}<\infty$. For existence let $F(\bm U)\coloneqq\int d_{\rm KL}(\bar p\|\ell_{\bm U})\,d\mu_1$, which is finite for some $\bm U\in\Sy_+^m$ when $\epsilon^{\rm Po}_{\rm class}<\infty$, by~\eqref{compsplit}. On the cone, $\int\ell_{\bm U}\,d\mu_1=\langle\int\bm C\,d\mu_1,\bm U\rangle\ge c_1\tr\bm U\ge c_1\|\bm U\|_F$ with $c_1\coloneqq\lambda_{\min}(\int\bm C\,d\mu_1)>0$, and $\ell_{\bm U}\le B^2\|\bm U\|_F$ under (P1). Since $d_{\rm KL}(p\|\ell)=(p\log p-p)+\ell-p\log\ell$ and $p\log p-p\ge-1$ on $[0,1]$, $F(\bm U)\ge c_1\|\bm U\|_F-\log(B^2\|\bm U\|_F)\int\bar p\,d\mu_1-\E[n_0]$, so the sublevel sets of $F$ on the cone are bounded. If $\bm U_n\to\bm U$, then $\ell_{\bm U_n}\to\ell_{\bm U}$ pointwise, and the lower semicontinuity of $d_{\rm KL}(\bar p\|\cdot)$ on $[0,\infty)$ and Fatou's lemma give $F(\bm U)\le\liminf_nF(\bm U_n)$. Hence $F$ attains its infimum on the closed cone, and by~\eqref{klident} and~\eqref{compsplit} so does the population risk.
\ref{prop_target_score}: $\frac{d}{d\ell}(\ell-y\log\ell)=(\ell-y)/\ell$, and $\partial\ell_{t,i}/\partial\bm U=\bm C_{t,i}$ by~\eqref{cdef}. For the calibration equations, a minimizer of a differentiable convex function over the closed convex cone satisfies $\langle\bm Z,\bm U-\hat{\bm U}\rangle\ge0$ for all $\bm U\in\Sy_+^m$, with $\bm Z\coloneqq\nabla\Rcal_T(\hat{\bm U})$, which on a cone is equivalent to $\bm Z\succeq\bm 0$ and $\langle\bm Z,\hat{\bm U}\rangle=0$, hence $\bm Z\hat{\bm U}=\bm 0$. For $\bm v$ in the range of $\hat{\bm U}$ this gives $\bm Z\bm v=\bm 0$, so $\bm v'^\top\bm Z\bm v=0$ for every $\bm v'\in\R^m$, and $\bm v'^\top\bm C_{t,i}\bm v=(\bm\phi(z_{t,i})\bm v)(\bm\phi(z_{t,i})\bm v')$ gives~\eqref{calibeq}. Summing over an orthonormal eigenbasis of the range of $\hat{\bm U}$ weighted by the eigenvalues gives $\langle\bm Z,\hat{\bm U}\rangle=0$ with $\langle\bm C_{t,i},\hat{\bm U}\rangle=\hat\ell_{t,i}$, and $(y/\hat\ell-1)\hat\ell=y-\hat\ell$ gives~\eqref{aggcalemp}.
\end{proof}

\subsection{Auxiliary results for Section~\ref{sec_theory}}\label{app_aux}

\begin{lemma}[Well-defined population risk]\label{lem_poprisk}
Assume $\{\xi_t\}_{t\ge0}$ is strictly stationary and ergodic and there exists an open convex set $\Ocal\subset\Sy^m$ with $\E[|\Rcal(\bm U,\xi_0)|]<\infty$ for every $\bm U\in\Ocal$. Then the population risk $\Ecal(\bm U)\coloneqq\E[\Rcal(\bm U,\xi_0)]$ is finite and convex on $\Ocal$, and $\Rcal_T\to\Ecal$ as $T\to\infty$, uniformly on compacts of $\Ocal$, $\Pa$-almost surely.
\end{lemma}

\begin{proof}
For fixed $\bm U\in\Ocal$ the sequence $\{\Rcal(\bm U,\xi_t)\}$ is stationary, ergodic, and integrable, so Birkhoff's ergodic theorem \citep[Chapter~1]{kre_85} gives $\Rcal_T(\bm U)\to\Ecal(\bm U)$ $\Pa$-a.s., and linearity of expectation gives convexity of $\Ecal$. Let $\Ocal_{\rm d}\subset\Ocal$ be countable and dense. On the intersection of the corresponding full-measure events, $\Rcal_T\to\Ecal$ pointwise on $\Ocal_{\rm d}$, and each $\Rcal_T$ is finite on $\Ocal$ by convexity, since every point of $\Ocal$ lies in the interior of a simplex with vertices in $\Ocal_{\rm d}$. \citet[Theorem~10.8]{roc_97}, applied in the $M$-dimensional space $\Sy^m$, upgrades pointwise convergence of finite convex functions on a dense subset to convergence at every point of $\Ocal$, uniform on compacts, along the same sample path. The limit is finite, convex, and continuous on $\Ocal$, and coincides with $\Ecal$ on $\Ocal_{\rm d}$, hence everywhere on $\Ocal$.
\end{proof}

\begin{lemma}[Pseudo-self-concordance]\label{lem_SC}
$\Pa$-a.s., $\bm U\mapsto\Rcal(\bm U,\xi)$ is three times continuously differentiable on $\Ocal$ and satisfies~\eqref{pseudoSC} for all $\bm U\in\Ocal$ and $\bm A,\bm\Delta\in\Sy^m$.
The same inequality, with the same constant, holds for $\Rcal_T$ and, by dominated convergence, for $\Ecal$ on $\Ocal$.
\end{lemma}

\begin{proof}
With $\ell_i\coloneqq\langle\bm C_i,\bm U\rangle$, Lemma~\ref{lem_vec} gives $D^2\Rcal[\bm A,\bm A]=\sum_{i\in S}\langle\bm C_i,\bm A\rangle^2/\ell_i^2$ and $D^3\Rcal[\bm A,\bm A,\bm\Delta]=-2\sum_{i\in S}\langle\bm C_i,\bm A\rangle^2\langle\bm C_i,\bm\Delta\rangle/\ell_i^3$, since the linear term contributes nothing. On $\Ocal$, Cauchy--Schwarz and~\eqref{cbound} give $|\langle\bm C_i,\bm\Delta\rangle|/\ell_i\le B^2\|\bm\Delta\|_F/\underline c$ for $i\in S$ by~\eqref{specbnd}, and bounding each summand yields~\eqref{pseudoSC}. Unlike standard self-concordance, which degrades under scaling by $1/T$ \citep[Chapter~5]{nesterov18}, the class defined by~\eqref{pseudoSC} is closed under averaging with the same constant. The statement for $\Ecal$ follows from the uniform Hessian bound $H$ and dominated convergence.
\end{proof}

\begin{lemma}[Deterministic localization]\label{lem_loc}
Let $f:\Sy^m\to(-\infty,+\infty]$ be convex and lower semicontinuous, three times continuously differentiable on an open convex set containing $\bar B(\bm U^{\rm Po}_\star,r_0)$ where it satisfies~\eqref{pseudoSC} with constant $M_\psi$, and let $\Ccal\subseteq\Sy^m$ be closed convex with $\bm U^{\rm Po}_\star\in\Ccal$. Assume $\nabla^2f(\bm U^{\rm Po}_\star)\succeq\nu\Ical$ with $\nu>0$, and $\nabla f(\bm U^{\rm Po}_\star)=\bm P+\bm Q$ with $\langle\bm P,\bm U-\bm U^{\rm Po}_\star\rangle\ge0$ for all $\bm U\in\Ccal$. If $3\|\bm Q\|_F/\nu<r_1\coloneqq\min\{r_0,M_\psi^{-1}\}$, then $\argmin_\Ccal f$ is a singleton $\{\hat{\bm U}\}$ and $\|\hat{\bm U}-\bm U^{\rm Po}_\star\|_F\le3\|\bm Q\|_F/\nu$.
\end{lemma}

\begin{proof}
Fix $\bm U\in\Ccal$ with $\bm\Delta\coloneqq\bm U-\bm U^{\rm Po}_\star$, $\delta\coloneqq\|\bm\Delta\|_F\le r_1$, and set $\psi(s)\coloneqq D^2f(\bm U^{\rm Po}_\star+s\bm\Delta)[\bm\Delta,\bm\Delta]$. By~\eqref{pseudoSC}, $\psi'(s)\ge-M_\psi\delta\psi(s)$, so Gr\"onwall's inequality gives $\psi(s)\ge\psi(0)\e^{-M_\psi\delta s}$, and integrating twice,
\[ f(\bm U)-f(\bm U^{\rm Po}_\star)-\langle\nabla f(\bm U^{\rm Po}_\star),\bm\Delta\rangle\ge D^2f(\bm U^{\rm Po}_\star)[\bm\Delta,\bm\Delta]\,\chi(M_\psi\delta),\qquad\chi(t)\coloneqq\frac{\e^{-t}+t-1}{t^2}, \]
with $\chi$ decreasing and $\chi(t)\ge\e^{-1}\ge1/3$ for $t\le1$. Since $\langle\nabla f(\bm U^{\rm Po}_\star),\bm\Delta\rangle\ge\langle\bm Q,\bm\Delta\rangle\ge-\|\bm Q\|_F\,\delta$,
\[ f(\bm U)-f(\bm U^{\rm Po}_\star)\ge\delta\Big(\frac{\nu\delta}{3}-\|\bm Q\|_F\Big)>0\qquad\text{whenever }\frac{3\|\bm Q\|_F}{\nu}<\delta\le r_1. \]
For $\bm U\in\Ccal$ with $\|\bm\Delta\|_F>r_1$, the point $\bm W$ where the segment $[\bm U^{\rm Po}_\star,\bm U]$ crosses the sphere of radius $r_1$ lies in $\Ccal$ and satisfies $f(\bm W)>f(\bm U^{\rm Po}_\star)$, so convexity gives $f(\bm U)>f(\bm U^{\rm Po}_\star)$. Every minimizer thus lies in the compact set $\Ccal\cap\bar B(\bm U^{\rm Po}_\star,3\|\bm Q\|_F/\nu)$, on which lower semicontinuity yields existence. Finally, the two-sided Hessian bound of pseudo-self-concordant functions \citep[Proposition~1]{bach10} gives $\nabla^2f(\bm U)\succeq\e^{-M_\psi\|\bm\Delta\|_F}\nabla^2f(\bm U^{\rm Po}_\star)\succeq\e^{-1}\nu\Ical\succ\bm 0$ on $\bar B(\bm U^{\rm Po}_\star,r_1)$, so $f$ is strictly convex there and the minimizer is unique.
\end{proof}

\begin{lemma}[Concentration under geometric mixing]\label{lem_concmix}
Under the hypotheses of Theorem~\ref{thm_main}\ref{thm_main_conc}, let $\rmix(k)\le\e^{-2c_{\rm mix}k}$ for some $c_{\rm mix}>0$ and $|\Psi_{0,j}|\le b_j$, $\Pa$-a.s. With $\bar b_g\coloneqq\max_jb_j$, $\bar v^2\coloneqq\max_j\sum_{k\in\mathbb Z}|\E[\Psi_{0,j}\Psi_{k,j}]|$ and $v_H^2\coloneqq4H^2(1+8\sum_{k\ge1}\rmix(k))$, there is $C_{\rm mix}>0$ depending only on $c_{\rm mix}$ such that, for $T\ge4$,
\begin{equation}\label{concmix}
\begin{split}
    \Pa[F_\epsilon]\le\ &2M\exp\Big(-\frac{C_{\rm mix}Tr_\epsilon^2}{M(\bar v^2+\bar b_g^2/T)+\sqrt M\,\bar b_g(\log T)^2r_\epsilon}\Big)\\
    &+2M^2\exp\Big(-\frac{C_{\rm mix}T\alpha^2}{4M^2(v_H^2+4H^2/T)+4M\alpha H(\log T)^2}\Big).
\end{split}
\end{equation}
It is proved with Theorem~\ref{thm_main}\ref{thm_main_conc} in Appendix~\ref{app_proofs_theory}.
\end{lemma}

\subsection{Proofs for Section~\ref{sec_theory}}\label{app_proofs_theory}

\begin{proof}[Proof of Theorem~\ref{thm_main}]
\ref{thm_main_cons}: $\Rcal(\cdot,\xi)$ is convex and lower semicontinuous on $\Sy^m$ (Lemma~\ref{lem_vec}) and measurable in $\xi$, hence random lower semicontinuous, and so is $\Rcal_T$. Fix $r_0>0$ with $\bar B(\bm U^{\rm Po}_\star,r_0)$ inside the integrability set $\Ocal$ of this part, and work on the $\Pa$-full event of Lemma~\ref{lem_poprisk} on which $\Rcal_T\to\Ecal$ uniformly on $\bar B(\bm U^{\rm Po}_\star,r_0)$. Fix $\delta\in(0,r_0]$ and the sphere section $\Ucal_\delta\coloneqq\Sy_+^m\cap\{\|\bm U-\bm U^{\rm Po}_\star\|_F=\delta\}$, which is compact and, for the cone, nonempty, since $\bm U^{\rm Po}_\star+s\bm I_m\in\Sy_+^m$ for every $s\ge0$. $\Ecal$ is continuous on it, and uniqueness of $\bm U^{\rm Po}_\star$ gives $\Delta_\delta\coloneqq\min_{\Ucal_\delta}\Ecal-\Ecal(\bm U^{\rm Po}_\star)>0$. For $T$ large, $\sup_{\bar B(\bm U^{\rm Po}_\star,r_0)}|\Rcal_T-\Ecal|<\Delta_\delta/3$, whence $\Rcal_T(\bm U^{\rm Po}_\star)<\Ecal(\bm U^{\rm Po}_\star)+\Delta_\delta/3$ and $\min_{\Ucal_\delta}\Rcal_T>\Ecal(\bm U^{\rm Po}_\star)+2\Delta_\delta/3$. For $\bm U\in\Sy_+^m$ outside the ball, the crossing point $\bm W\in\Ucal_\delta$ of the segment $[\bm U^{\rm Po}_\star,\bm U]$ lies in $\Sy_+^m$, and convexity of $\Rcal_T$ gives $\Rcal_T(\bm U)\ge\Rcal_T(\bm U^{\rm Po}_\star)+\Delta_\delta/3$. Hence all minimizers lie in the compact set $\Sy_+^m\cap\bar B(\bm U^{\rm Po}_\star,\delta)$, where lower semicontinuity yields existence. Applying this along $\delta_k\downarrow0$ gives $\hat{\bm U}^{\rm Po}_T\to\bm U^{\rm Po}_\star$ $\Pa$-a.s. No compactness, level-boundedness, or coercivity is required, because convexity confines the minimizers once the empirical risk dominates on the sphere \citep[cf.][Theorems~5.3--5.4]{sha_den_rus_14}.

\ref{thm_main_det}: On $\Acal_T^{\rm H}$, $\nabla^2\Rcal_T(\bm U^{\rm Po}_\star)\succeq\Hcal^{\rm Po}_\star-\frac\alpha2\Ical\succeq\frac\alpha2\Ical$. Decompose $\nabla\Rcal_T(\bm U^{\rm Po}_\star)=\nabla\Ecal(\bm U^{\rm Po}_\star)+\bar{\bm\Psi}^{\rm Po}_T$, and note that minimality of $\bm U^{\rm Po}_\star$ over $\Sy_+^m$ gives $\langle\nabla\Ecal(\bm U^{\rm Po}_\star),\bm U-\bm U^{\rm Po}_\star\rangle\ge0$ on $\Sy_+^m$. Lemma~\ref{lem_loc} with $f=\Rcal_T$ (which satisfies~\eqref{pseudoSC} by Lemma~\ref{lem_SC}), $\Ccal=\Sy_+^m$, $\nu=\alpha/2$, $\bm P=\nabla\Ecal(\bm U^{\rm Po}_\star)$, $\bm Q=\bar{\bm\Psi}^{\rm Po}_T$ applies precisely on $\Acal_T^{\rm s}$ and yields existence, uniqueness, and~\eqref{detloc}.

\ref{thm_main_conc}: The inclusion follows from part~\ref{thm_main_det} as in the display. Under i.i.d.\ sampling, the deviation bound for sums of independent norm-sub-Gaussian vectors of \citet[Theorem~1]{mil_23}, applied in the $M$-dimensional Euclidean space $\Sy^m$, gives $\Pa[\|\bar{\bm\Psi}^{\rm Po}_T\|_F\ge r_\epsilon]\le2\exp(-Tr_\epsilon^2/(3\tau^2))$, and $\sigma^2\le\tau^2$ follows from Jensen's inequality. The matrix Bernstein inequality \citep[Theorem~1.6.2]{tropp15}, applied to the independent centered self-adjoint summands $\frac1T(\nabla^2\Rcal(\bm U^{\rm Po}_\star,\xi_t)-\Hcal^{\rm Po}_\star)$ acting on the $M$-dimensional space $\Sy^m$, with norm bound $2H/T$ and variance parameter $4H^2/T$, bounds $\Pa[(\Acal_T^{\rm H})^{\rm c}]$, yielding~\eqref{conciid}. Under geometric mixing, apply the Bernstein inequality for strongly mixing sequences \citep[Theorem~2]{mer_pel_rio_09} to the coordinate sequences $\{\Psi_{t,j}\}$ with $x=Tr_\epsilon/\sqrt M$, and to the entries $\langle\bm E_j,(\nabla^2\Rcal(\bm U^{\rm Po}_\star,\xi_t)-\Hcal^{\rm Po}_\star)[\bm E_l]\rangle$, bounded by $2H$, with $x=T\alpha/(2M)$, using the covariance inequality for bounded mixing sequences \citep[Section~1.2.2]{dou_94} to dominate the variance proxies by $\bar v^2$ and $v_H^2$. The strong-mixing coefficient of \citet{mer_pel_rio_09} is~\eqref{alphamixdef}, and their bound is two-sided. Union bounds over the $M$ coordinates and the $M^2$ entries give~\eqref{concmix}. The factor two in the prefactors, and the count $M^2$ rather than the $M(M+1)/2$ distinct entries of a self-adjoint operator, are conservative. Orthonormality of $\{\bm E_j\}$ is used twice, to pass from $\|\bar{\bm\Psi}^{\rm Po}_T\|_F$ to its coordinates and from the operator norm to the entries, through $\opn{\Xcal}\le M\max_{j,l}|\langle\bm E_j,\Xcal[\bm E_l]\rangle|$, and it keeps the constants free of any basis-dependent distortion.

\ref{thm_main_mse}: By stationarity, $\E[\|\bar{\bm\Psi}^{\rm Po}_T\|_F^2]=\frac1T\sum_{|k|\le T-1}(1-|k|/T)\gamma(k)\le\sigma^2_{\rm LR}/T$, and the bound follows from~\eqref{detloc}. The covariance bound $|\gamma(k)|\le4\rmix(k)\|\bm b\|_2^2$ is the covariance inequality for bounded strongly mixing variables \citep[Section~1.2.2]{dou_94} applied coordinatewise, since the mixing coefficients of $\{\bm\Psi^{\rm Po}_t\}$ are dominated by those of $\{\xi_t\}$.

\ref{thm_main_curve}: Immediate from the linearity $\hat\ell^{\rm Po}(z)-\ell^{\rm Po}_\star(z)=\langle\bm C(z),\hat{\bm U}^{\rm Po}_T-\bm U^{\rm Po}_\star\rangle$, the Cauchy--Schwarz inequality, and the identity $\|\bm C(z)\|_F=\|\bm\phi(z)\|_2^2\le B^2$ of~\eqref{cbound}, which (P1) supplies. The transfer is exact, with no remainder term.

\ref{thm_main_prim}: On $\Ocal$, $|\log\langle\bm C_i,\bm U\rangle|\le\max\{|\log\underline c|,|\log\overline c|\}$ for $i\in S$, so $|\Rcal(\bm U,\xi)|\le\bar s\max\{|\log\underline c|,|\log\overline c|\}+\bar b\sup_{\Ocal}\|\bm U\|_F$ is bounded uniformly on the bounded set $\Ocal$, giving integrability. For $\bm U\in\Sy^m$, $\langle\bm B(\xi),\bm U\rangle\ge-\bar b\|\bm U\|_F$ and $\langle\bm C_i,\bm U\rangle\le B^2\|\bm U\|_F$, so $\Rcal(\bm U,\xi)\ge-\bar s\max\{0,\log(B^2\|\bm U\|_F)\}-\bar b\|\bm U\|_F$ $\Pa$-a.s. On $\Sy_+^m$ the loss is even nonnegative, $\Lcal^{\rm Po}\ge|S|$, as in the proof of Proposition~\ref{prop_target}. The bound is uniform on bounded sets, so $\Ecal$ is well defined in $(-\infty,+\infty]$ and convex on $\Sy^m$, and lower semicontinuous by Fatou's lemma and the lower semicontinuity of $\Rcal(\cdot,\xi)$, as Lemma~\ref{lem_loc} requires when it is applied to $\Ecal$, $\Ecal_0$ and $\Ecal_\lambda$. The gradient bound $\|\nabla\Rcal(\bm U,\xi)\|_F\le B^2\bar s/\underline c+\bar b$ on $\Ocal$ gives~\eqref{scorebnd} after centering, and with the Hessian bound $H$ justifies both interchanges by dominated convergence. A bounded random element of a Euclidean space is norm-sub-Gaussian for some $\tau$, and summable mixing gives $\sigma^2_{\rm LR}<\infty$ by part~\ref{thm_main_mse}.

\ref{thm_main_restr}: Lemma~\ref{lem_loc} is stated for an arbitrary closed convex set, and the only use of the feasible set in parts~\ref{thm_main_det}--\ref{thm_main_mse} is the variational inequality at $\bm U^{\rm Po}_\star$, which holds on any closed convex $\Ccal\ni\bm U^{\rm Po}_\star$ over which $\bm U^{\rm Po}_\star$ minimizes $\Ecal$. Part~\ref{thm_main_cons} uses one further property of the cone, that $\Ucal_\delta$ is nonempty for every $\delta\le r_0$. For a general $\Ccal$ the argument is unchanged when $\Ccal\cap\{\|\bm U-\bm U^{\rm Po}_\star\|_F=\delta\}\ne\emptyset$, and when it is empty, $\Ccal\subseteq\bar B(\bm U^{\rm Po}_\star,\delta)$ and the confinement of the minimizers is immediate.
\end{proof}

\begin{proof}[Proof of Theorem~\ref{thm_clt}]
\ref{thm_clt_par}: \emph{Absolute convergence.} The centered score is bounded (Theorem~\ref{thm_main}\ref{thm_main_prim}), so the covariance inequality for bounded strongly mixing variables \citep[Section~1.2.2]{dou_94}, applied to the coordinates of Section~\ref{sec_theory_po}, gives $|\langle\bm E_j,\Gcal(k)[\bm E_l]\rangle|\le4\rmix(k)b_jb_l$, whence $\opn{\Gcal(k)}\le4\rmix(k)\|\bm b\|_2^2$, summable by Assumption~\ref{ass_clt}\ref{ass_clt_ii}.

By Assumption~\ref{ass_clt}\ref{ass_clt_i} and convexity, $\bm U^{\rm Po}_\star$ minimizes $\Ecal$ over $\Kcal_0$, uniquely. The Gr\"onwall argument of Lemma~\ref{lem_loc} applied to $\Ecal$ gives $\nabla^2\Ecal\succeq\e^{-1}\alpha\Ical$ on $\bar B(\bm U^{\rm Po}_\star,r_1)$, and convexity extends strict minimality to $\Kcal_0$. Birkhoff's theorem applied to the bounded stationary sequences $\{\bm\Psi^{\rm Po}_t\}$ and $\{\nabla^2\Rcal(\bm U^{\rm Po}_\star,\xi_t)\}$, coordinatewise in the orthonormal basis, gives $\bar{\bm\Psi}^{\rm Po}_T\to\bm 0$ and $\nabla^2\Rcal_T(\bm U^{\rm Po}_\star)\to\Hcal^{\rm Po}_\star$ $\Pa$-a.s., so the events $\Acal_T^{\rm H},\Acal_T^{\rm s}$ hold eventually. Lemma~\ref{lem_loc} with $\Ccal=\Kcal_0$, $\bm P=\bm 0$, $\bm Q=\bar{\bm\Psi}^{\rm Po}_T$ then yields a unique minimizer with $\|\hat{\bm U}^{\rm Po}_{\rm b}-\bm U^{\rm Po}_\star\|_F\le(6/\alpha)\|\bar{\bm\Psi}^{\rm Po}_T\|_F\to0$, eventually interior, satisfying $\nabla\Rcal_T(\hat{\bm U}^{\rm Po}_{\rm b})=\bm 0$.

With $\bm\Delta\coloneqq\hat{\bm U}^{\rm Po}_{\rm b}-\bm U^{\rm Po}_\star$ and $\widetilde{\Hcal}_T\coloneqq\int_0^1\nabla^2\Rcal_T(\bm U^{\rm Po}_\star+s\bm\Delta)\,ds$, the first-order condition gives $\bm 0=\bar{\bm\Psi}^{\rm Po}_T+\widetilde{\Hcal}_T[\bm\Delta]$. Pseudo-self-concordance sandwiches $\e^{-M_\psi\|\bm\Delta\|_F}\nabla^2\Rcal_T(\bm U^{\rm Po}_\star)\preceq\nabla^2\Rcal_T(\bm U^{\rm Po}_\star+s\bm\Delta)\preceq\e^{M_\psi\|\bm\Delta\|_F}\nabla^2\Rcal_T(\bm U^{\rm Po}_\star)$, whence $\widetilde{\Hcal}_T\to\Hcal^{\rm Po}_\star$ and $\nabla^2\Rcal_T(\hat{\bm U}^{\rm Po}_{\rm b})\to\Hcal^{\rm Po}_\star$ $\Pa$-a.s. For fixed $\bm A\in\Sy^m$, $\{\langle\bm A,\bm\Psi^{\rm Po}_t\rangle\}$ is stationary, centered, and bounded, with mixing coefficients dominated by $\rmix$. The central limit theorem for bounded strongly mixing sequences with summable mixing coefficients \citep[Theorem~18.5.4]{ibr_lin_71} and the Cram\'er--Wold device applied in $\Sy^m$ give $\sqrt T\bar{\bm\Psi}^{\rm Po}_T\to_d$ a centered Gaussian element with covariance operator $\Qcal_{\rm LR}$, and in a degenerate direction $\var(\sqrt T\langle\bm A,\bar{\bm\Psi}^{\rm Po}_T\rangle)\to0$. Slutsky's theorem applied to $\sqrt T\bm\Delta=-\widetilde{\Hcal}_T^{-1}[\sqrt T\bar{\bm\Psi}^{\rm Po}_T]$ yields~\eqref{cltdisplay}.

\ref{thm_clt_curve}: It follows from part~\ref{thm_clt_par} and the continuous mapping theorem applied to the linear map $\bm W\mapsto\langle\bm C(\cdot),\bm W\rangle$, which is a bounded linear map from $\Sy^m$ to $C_b(\Zcal)$, of norm at most $B^2$ by~\eqref{cbound}, when $\bm\phi$ is continuous and bounded.

\ref{thm_clt_bdy}: Let $\bm\Delta_T\coloneqq\sqrt T(\hat{\bm U}^{\rm Po}_T-\bm U^{\rm Po}_\star)$, the minimizer over $\sqrt T(\Sy_+^m-\bm U^{\rm Po}_\star)$ of the convex process $Z_T(\bm A)\coloneqq T\big[\Rcal_T(\bm U^{\rm Po}_\star+T^{-1/2}\bm A)-\Rcal_T(\bm U^{\rm Po}_\star)\big]$. The argument is that of \citet[Lemma~4.1 and Theorem~4.4]{gey_94} for i.i.d.\ data and of \citet[Theorem~3]{andrews1999b} under high-level conditions that admit dependence. Both assume that the gradient of the population objective vanishes at the target, Geyer's Assumption~A and Andrews's Assumption~3 with norming $T^{1/2}$, which is Assumption~\ref{ass_clt}\ref{ass_clt_i}. The sequence $\bm\Delta_T$ is tight by the localization of Theorem~\ref{thm_main}\ref{thm_main_det}, which holds eventually as in the proof of part~\ref{thm_clt_par}, and the central limit theorem of that proof. By (P3), Assumption~\ref{ass_clt}\ref{ass_clt_i} and the pseudo-self-concordant sandwich of the same proof, $Z_T$ is finite on compact sets for $T$ large and equals $\langle\sqrt T\bar{\bm\Psi}^{\rm Po}_T,\bm A\rangle+\frac12\langle\bm A,\Hcal^{\rm Po}_\star[\bm A]\rangle+o_\Pa(1)$ uniformly there, with $\sqrt T\bar{\bm\Psi}^{\rm Po}_T\to_d\bm\Psi_\infty\sim N(\bm 0,\Qcal_{\rm LR})$, so it converges in distribution uniformly on compact sets \citep{hjortpollard93}. The constraint enters as the indicator of the feasible set $\sqrt T(\Sy_+^m-\bm U^{\rm Po}_\star)$ of $\bm\Delta_T$. These sets increase to the cone of feasible directions, whose closure is~\eqref{tangentcone} \citep[Section~2.2]{bonnansshapiro00}, so their indicators epi-converge to that of $\Tcal(\bm U^{\rm Po}_\star)$, which is Chernoff regularity \citep{gey_94}, and the limit is strictly convex by~\eqref{curv}. The argmin theorem for constrained convex processes \citep[Theorem~4.4]{gey_94} then gives $\bm\Delta_T\to_d\argmin_{\Tcal(\bm U^{\rm Po}_\star)}\{\langle\bm\Psi_\infty,\bm A\rangle+\frac12\langle\bm A,\Hcal^{\rm Po}_\star[\bm A]\rangle\}$, which is~\eqref{bdydisplay} after completing the square with $\bm\eta\coloneqq-{\Hcal^{\rm Po}_\star}^{-1}[\bm\Psi_\infty]$, and the joint statement follows by linearity of the fitted values. If $r<m$ and $\Vcal$ is nondegenerate, the limit law has a nonzero absolutely continuous part, on $\{\bm\eta\in{\rm int}\,\Tcal(\bm U^{\rm Po}_\star)\}$, and a nonzero part on the Lebesgue-null boundary of $\Tcal(\bm U^{\rm Po}_\star)$, the image of $\{\bm\eta\notin\Tcal(\bm U^{\rm Po}_\star)\}$, so it is neither a nondegenerate nor a degenerate Gaussian law.

\ref{thm_clt_sc}: Write $\bm U^{\rm Po}_\star=\bm Q_r\bm\Theta_r\bm Q_r^\top$ with $\bm\Theta_r\succ\bm 0$ diagonal and $\bm Q\coloneqq[\bm Q_r\ \bm V]$ orthogonal. \emph{Localization.} The centered score, the Hessian and the central limit theorem of part~\ref{thm_clt_par} do not use Assumption~\ref{ass_clt}\ref{ass_clt_i}. As there, the events of Theorem~\ref{thm_main} hold eventually, and Theorem~\ref{thm_main}\ref{thm_main_det}, whose proof takes $\bm P=\nabla\Ecal(\bm U^{\rm Po}_\star)=\bm\Lambda$, gives $\|\hat{\bm U}^{\rm Po}_T-\bm U^{\rm Po}_\star\|_F\le(6/\alpha)\|\bar{\bm\Psi}^{\rm Po}_T\|_F$, which tends to zero $\Pa$-a.s.\ and is $O_\Pa(T^{-1/2})$. \emph{Rank.} By the first-order conditions in the proof of Proposition~\ref{prop_target}\ref{prop_target_score}, $\hat{\bm\Lambda}_T\coloneqq\nabla\Rcal_T(\hat{\bm U}^{\rm Po}_T)\succeq\bm 0$ and $\hat{\bm\Lambda}_T\hat{\bm U}^{\rm Po}_T=\bm 0$. Since $\nabla\Rcal_T(\bm U^{\rm Po}_\star)=\bm\Lambda+\bar{\bm\Psi}^{\rm Po}_T$ and the Hessian is bounded by $H$ on $\Ocal$, $\|\hat{\bm\Lambda}_T-\bm\Lambda\|_F\le H\|\hat{\bm U}^{\rm Po}_T-\bm U^{\rm Po}_\star\|_F+\|\bar{\bm\Psi}^{\rm Po}_T\|_F\to0$. The rank is lower semicontinuous, so eventually $\rank\hat{\bm\Lambda}_T\ge\rank\bm\Gamma=m-r$, whence $\rank\hat{\bm U}^{\rm Po}_T\le r$, and $\hat{\bm U}^{\rm Po}_T\to\bm U^{\rm Po}_\star$ gives $\rank\hat{\bm U}^{\rm Po}_T\ge r$. \emph{Chart.} For $\bm\theta=(\bm a,\bm b)$ with $\bm a\in\Sy^r$ and $\bm b\in\R^{r\times(m-r)}$ near zero let
\[ \iota(\bm\theta)\coloneqq\bm Q\begin{pmatrix}\bm\Theta_r+\bm a&\bm b\\ \bm b^\top&\bm b^\top(\bm\Theta_r+\bm a)^{-1}\bm b\end{pmatrix}\bm Q^\top. \]
A symmetric matrix whose upper left block in the basis $\bm Q$ is positive definite is positive semidefinite of rank $r$ if and only if its Schur complement vanishes, so $\iota$ is a smooth bijection of a neighborhood of zero onto a neighborhood of $\bm U^{\rm Po}_\star$ in the positive semidefinite matrices of rank $r$, with $\iota(\bm 0)=\bm U^{\rm Po}_\star$. Its differential $J\coloneqq D\iota(\bm 0)$, $J[\bm\theta]=\bm Q\big(\begin{smallmatrix}\bm a&\bm b\\ \bm b^\top&\bm 0\end{smallmatrix}\big)\bm Q^\top$, is an isometry onto $\Tcal_\star$ when the parameters carry the inner product $\langle\bm\theta,\bm\theta'\rangle\coloneqq\tr(\bm a\bm a')+2\tr(\bm b^\top\bm b')$, so that $J^*J$ is the identity and $JJ^*=\Pcal$, and $D^2\iota(\bm 0)[\bm\theta,\bm\theta]=\bm Q\big(\begin{smallmatrix}\bm 0&\bm 0\\ \bm 0&2\bm b^\top\bm\Theta_r^{-1}\bm b\end{smallmatrix}\big)\bm Q^\top$. \emph{Expansion.} Eventually $\hat{\bm U}^{\rm Po}_T=\iota(\hat{\bm\theta}_T)$ with $\hat{\bm\theta}_T\to\bm 0$, and since $\iota$ maps into the cone, $\hat{\bm\theta}_T$ is an unconstrained local minimizer of $\Rcal_T\circ\iota$, which is three times continuously differentiable near zero by (P3), so $\nabla (\Rcal_T\circ\iota)(\hat{\bm\theta}_T)=\bm 0$. The multiplier is orthogonal to $\Tcal_\star$, $\langle\bm\Lambda,\bm A\rangle=\tr(\bm\Gamma\bm V^\top\bm A\bm V)=0$, so $\nabla (\Rcal_T\circ\iota)(\bm 0)=J^*[\bm\Lambda+\bar{\bm\Psi}^{\rm Po}_T]=J^*[\bar{\bm\Psi}^{\rm Po}_T]$. The Hessian of $(\Rcal_T\circ\iota)$ acts as $D^2(\Rcal_T\circ\iota)(\bm\theta)[\bm\vartheta,\bm\vartheta]=D^2\Rcal_T(\iota(\bm\theta))[D\iota(\bm\theta)\bm\vartheta,D\iota(\bm\theta)\bm\vartheta]+\langle\nabla\Rcal_T(\iota(\bm\theta)),D^2\iota(\bm\theta)[\bm\vartheta,\bm\vartheta]\rangle$. On balls shrinking to zero it converges, by the bound $H$, the pseudo-self-concordant sandwich and Birkhoff's theorem as in part~\ref{thm_clt_par}, to $\Kcal_\theta\coloneqq J^*(\Hcal^{\rm Po}_\star+\Scal_\star)J$. Indeed, with $\bm A=J[\bm\theta]$, the curvature term is $\langle\bm\Lambda,D^2\iota(\bm 0)[\bm\theta,\bm\theta]\rangle=2\tr(\bm\Gamma\bm b^\top\bm\Theta_r^{-1}\bm b)=2\tr(\bm\Lambda\bm A{\bm U^{\rm Po}_\star}^{+}\bm A)=\langle\bm A,\Scal_\star[\bm A]\rangle$. The limit $\Kcal_\theta$ is positive definite by~\eqref{curv}, $\Scal_\star\succeq\bm 0$ and injectivity of $J$. The mean value theorem then gives $\sqrt T\hat{\bm\theta}_T=-\Kcal_\theta^{-1}J^*[\sqrt T\bar{\bm\Psi}^{\rm Po}_T]+o_\Pa(1)$, and $\sqrt T(\hat{\bm U}^{\rm Po}_T-\bm U^{\rm Po}_\star)=J[\sqrt T\hat{\bm\theta}_T]+O_\Pa(T^{-1/2})$. With $\sqrt T\bar{\bm\Psi}^{\rm Po}_T\to_d\bm\Psi_\infty\sim N(\bm 0,\Qcal_{\rm LR})$ and the identity $J\Kcal_\theta^{-1}J^*=\big(\Pcal(\Hcal^{\rm Po}_\star+\Scal_\star)\Pcal\big)^+$, which holds because $J$ is an isometry onto $\Tcal_\star$, the limit is $-\big(\Pcal(\Hcal^{\rm Po}_\star+\Scal_\star)\Pcal\big)^+\Pcal\bm\Psi_\infty$, which is~\eqref{scdisplay} with covariance $\Vcal_\star$. The curve statement follows by linearity of the fitted values.

\ref{thm_clt_ci}: It combines part~\ref{thm_clt_curve} with $\hat{\Hcal}_T\to\Hcal^{\rm Po}_\star$ (part~\ref{thm_clt_par}), consistency of $\hat\Qcal_T$, and Slutsky's theorem. The empirical score variance is consistent under the martingale condition. On $\Ocal$, $\|\nabla\Rcal(\hat{\bm U}^{\rm Po}_{\rm b},\xi_t)-\nabla\Rcal(\bm U^{\rm Po}_\star,\xi_t)\|_F\le H\|\hat{\bm U}^{\rm Po}_{\rm b}-\bm U^{\rm Po}_\star\|_F\to0$, $\nabla\Rcal(\bm U^{\rm Po}_\star,\xi_t)=\bm\Psi^{\rm Po}_t$ is bounded by~\eqref{scorebnd}, and Birkhoff's theorem applies to $\langle\bm\Psi^{\rm Po}_t,\bm A\rangle\bm\Psi^{\rm Po}_t$. The final identification of $\ell^{\rm Po}_\star$ with $p$ is condition~(CS). Under strict complementarity the plug-ins of Remark~\ref{rem_faceband} are consistent because the rank of the fit identifies $r$ eventually, by part~\ref{thm_clt_sc}, so that the projections onto the range and the kernel of the fit converge to those of $\bm U^{\rm Po}_\star$ and $\hat{\bm\Lambda}_T\to\bm\Lambda$.
\end{proof}

\begin{proof}[Proof of Proposition~\ref{prop_mds}]
\ref{prop_mds_i}: $\bm C_{t,i}$ and $p_{t,i}$ are $\Fcal_t$-measurable and $\E[\mathbbm{1}_{\{i\in S_{t+1}\}}\mid\Fcal_t]=p_{t,i}$, so $\E[\bm\Psi^{\rm Po}_t\mid\Fcal_t]=\bm 0$, the matrix form following from the second identity of~\eqref{cbound}. The score is centered by Assumption~\ref{ass_clt}\ref{ass_clt_i} and Proposition~\ref{prop_target}\ref{prop_target_proj}. For $k\ge1$ and any $\bm A$, $\Gcal(k)[\bm A]=\E\big[\langle\E[\bm\Psi^{\rm Po}_k\mid\Fcal_k],\bm A\rangle\,\bm\Psi^{\rm Po}_0\big]=\bm 0$ by the tower property, $\bm\Psi^{\rm Po}_0$ being $\Fcal_k$-measurable for $k\ge1$.
\ref{prop_mds_ii}: Expanding and conditioning on $\Fcal_0$, $\E[(\mathbbm{1}_{\{i\in S_1\}}/p_{0,i}-1)(\mathbbm{1}_{\{j\in S_1\}}/p_{0,j}-1)\mid\Fcal_0]=(\pi_{0,ij}-p_{0,i}p_{0,j})/(p_{0,i}p_{0,j})$, which is~\eqref{Gamma0}, the diagonal terms using $\pi_{0,ii}=p_{0,i}$.
\end{proof}

\begin{proof}[Proof of Corollary~\ref{cor_calib}]
\ref{cor_calib_dict}: For $p,\ell>0$, $d_{\rm KL}(p\|\ell)=\int_p^\ell(s-p)/s\,ds$, and $1/s\ge1/\max(p,\ell)$ on the segment between $p$ and $\ell$ gives $d_{\rm KL}(p\|\ell)\ge(p-\ell)^2/(2\max(p,\ell))$, where $\max(p,\ell)\le\max(1,\overline c)$. Evaluating~\eqref{epsdictdef} at its minimizer, summing, and taking expectations gives the claim.
\ref{cor_calib_lam}: The regularized population risk is $\Ecal_\lambda=\Ecal_0+\lambda\langle\bm G,\cdot\rangle$, a linear perturbation. Minimality of $\bm U^{\rm Po}_{\star,\lambda}$ gives the variational inequality $\langle\nabla\Ecal_0(\bm U^{\rm Po}_{\star,\lambda})+\lambda\bm G,\bm U-\bm U^{\rm Po}_{\star,\lambda}\rangle\ge0$ on the cone, so $\bm U^{\rm Po}_{\star,\lambda}$ is an approximate stationary point of $\Ecal_0$ with residual $\lambda\|\bm G\|_F$, and Lemma~\ref{lem_loc} with $f=\Ecal_0$, centered at $\bm U^{\rm Po}_{\star,\lambda}$, where $\nabla^2\Ecal_0=\nabla^2\Ecal_\lambda\succeq\alpha\Ical$ by~\eqref{curv} since the penalty is linear, and $\nu=\alpha$ yields the bound. For the sign, compare the optimality of $\bm U^{\rm Po}_{\star,\lambda_1}$ and $\bm U^{\rm Po}_{\star,\lambda_2}$ in each other's objectives and add the two inequalities to obtain $(\lambda_2-\lambda_1)\big(\tr(\bm G\bm U^{\rm Po}_{\star,\lambda_2})-\tr(\bm G\bm U^{\rm Po}_{\star,\lambda_1})\big)\le0$.
\ref{cor_calib_total}: The triangle inequality in $\|\cdot\|$ applied to $p-\hat\ell^{\rm Po}=(p-\ell^{\rm Po}_{\star,0})+(\ell^{\rm Po}_{\star,0}-\ell^{\rm Po}_{\star,\lambda})+(\ell^{\rm Po}_{\star,\lambda}-\hat\ell^{\rm Po})$. The first term is part~\ref{cor_calib_dict}, since $\xi$ has the law of $\xi_0$. The last two are bounded uniformly in $z$ by $B^2$ times the respective parameter distances (Theorem~\ref{thm_main}\ref{thm_main_curve}), the parameters being fixed for the expectation over $\xi$, and hence contribute at most $B^2$ times those distances in the per-item norm.
\end{proof}

\subsection{Proofs for Section~\ref{sec_theory_log}}\label{app_proofs_log}

\begin{proof}[Proof of Theorem~\ref{thm_log}]
We go through the proofs of Appendix~\ref{app_proofs_theory} and Lemmas~\ref{lem_poprisk}--\ref{lem_loc} and list every place where the loss enters. All other steps use only convexity, the cone, the curvature condition~\eqref{curv} at the target, and stationarity and mixing.

Lemma~\ref{lem_poprisk} needs integrability on $\Ocal$, which the bound of Lemma~\ref{lem_SC_log}\ref{lem_SC_log_s} on $|\Rcal^{\log}|$ gives. Lemma~\ref{lem_loc} needs a convex lower semicontinuous function that is pseudo-self-concordant on a neighborhood of the target. $\Rcal^{\log}_T$ is one by~\eqref{Rlogvec} and part~\ref{lem_SC_log_sc}, and the inequality passes to averages with the same constant and to $\Ecal^{\log}$ by dominated convergence under the uniform bounds of parts~\ref{lem_SC_log_sc} and~\ref{lem_SC_log_h}, as in Lemma~\ref{lem_SC}.

\emph{Theorem~\ref{thm_main}.} Part~\ref{thm_main_cons} uses convexity, lower semicontinuity and measurability of the per-sample loss, from~\eqref{Rlogvec}. Part~\ref{thm_main_det} uses Lemma~\ref{lem_loc} with $f=\Rcal^{\log}_T$. Part~\ref{thm_main_conc} uses the score bound, which makes the score norm-sub-Gaussian and its coordinates bounded, and the Hessian bound $H$ in the matrix Bernstein inequality and in the entrywise bound $2H$ under mixing. Parts~\ref{lem_SC_log_h} and~\ref{lem_SC_log_s} supply $H^{\log}$ and the score bound. Part~\ref{thm_main_mse} is generic. Part~\ref{thm_main_curve} is the linearity of the curve in the parameter, unchanged. The transfer to the probability uses $|\e^{-a}-\e^{-b}|\le|a-b|$ for $a,b\ge0$. Part~\ref{thm_main_prim} uses integrability on $\Ocal$, a lower bound on all of $\Sy^m$, here $\Rcal^{\log}(\bm U,\xi)\ge-\bar b^{\log}\|\bm U\|_F$, since the non-event terms are nonnegative where finite and the event and penalty terms are linear with coefficient bounded by $\bar b^{\log}$, and the gradient and Hessian bounds for the interchanges. Part~\ref{thm_main_restr} is generic.

\emph{Theorem~\ref{thm_clt}.} Absolute convergence uses boundedness of the score. Localization uses Lemma~\ref{lem_loc} and the ergodic theorem for the bounded score and Hessian. The expansion uses the pseudo-self-concordant sandwich, and the boundary part uses (P3$^{\log}$) for finiteness of the local process on compacts and~\eqref{curv} for strict convexity of its limit. None of these uses the form of the loss beyond the facts listed. The curve statement for the probability is the delta method for $x\mapsto\e^{-x}$, and the identification in part~\ref{thm_clt_ci} is $\e^{-\ell^{\log}_\star}=p$ under the log-link reading of~(CS).

\emph{Proposition~\ref{prop_mds}.} Part~\ref{prop_mds_i} uses only that the score is a sum of $\Fcal_t$-measurable weights times the residuals $\mathbbm{1}_{\{i\in S_{t+1}\}}-p_{t,i}$. With the weights $-w_{\rm log}(\ell^{\log}_{\star,t,i})=1/(1-p_{t,i})$ the tower property argument is unchanged, and centering follows from the log-link reading of~(CS), under which the score $\sum_i(y_{t,i}-p_{t,i})\bm C_{t,i}/(1-p_{t,i})$ has conditional mean zero. Part~\ref{prop_mds_ii} follows by the same expansion with these weights, $\E[(\mathbbm{1}_{\{i\in S_1\}}-p_{0,i})(\mathbbm{1}_{\{j\in S_1\}}-p_{0,j})\mid\Fcal_0]=\pi_{0,ij}-p_{0,i}p_{0,j}$.

\emph{Corollary~\ref{cor_calib}.} Part~\ref{cor_calib_dict} is Pinsker's inequality for two-point laws, $d_{\rm Be}(p\|q)\ge2(p-q)^2$, at the minimizer of $\epsilon^{\log}_{\rm class}$. Part~\ref{cor_calib_lam} is Lemma~\ref{lem_loc} for $\Ecal^{\log}_0$, the penalty being linear, and part~\ref{cor_calib_total} is the triangle inequality with the $1$-Lipschitz map $x\mapsto\e^{-x}$ on $[0,\infty)$.

\emph{Theorem~\ref{thm_clt}\ref{thm_clt_sc} and Remark~\ref{rem_faceband}} use the Hessian, the score covariance and the multiplier of the fit, the localization and the pseudo-self-concordant sandwich, and hold as stated.

For $\lambda>0$ the objective is at least $\lambda\,\lambda_{\min}(\bm G)\tr\bm U$ on $\Sy_+^m$, since $g\ge0$ and the selected summands are nonnegative there, so it is coercive and lower semicontinuous on the closed cone. For $\lambda=0$ it is at least $\frac1T\lambda_{\min}\big(\sum_{t<T}\sum_{i\in S_{t+1}}\bm C_{t,i}\big)\tr\bm U$, and by the ergodic theorem and the tower property $\frac1T\sum_{t<T}\sum_{i\in S_{t+1}}\bm C_{t,i}\to\E\big[\sum_{i\le n_0}p_{0,i}\bm C(z_{0,i})\big]$ $\Pa$-a.s., positive definite under the stated condition.
\end{proof}

\subsection{Proofs for Sections~\ref{sec_loglink} and~\ref{sec_diag}}\label{app_proofs_diag}

\begin{proof}[Proof of Proposition~\ref{prop_loglink}]
\ref{prop_loglink_convex}: Each $\ell_{t,i}$ is linear in $\bm U$, the selected summand is linear, and $x\mapsto-\log(1-\e^{-x})$ has second derivative $\e^x/(\e^x-1)^2>0$ on $x>0$.
\ref{prop_loglink_target}: Expectation and the finite sum interchange, and $\E[\mathbbm{1}_{\{i\in S\}}]=p_i$ gives the closed form, which separates across $i$. In the mean $q=\e^{-\ell}$ each summand is the cross-entropy $-p\log q-(1-p)\log(1-q)$, strictly convex on $(0,1)$ for $p\in(0,1)$, divergent at both endpoints, with derivative $-p/q+(1-p)/(1-q)$ vanishing uniquely at $q=p$. The map $\ell\mapsto\e^{-\ell}$ is a strictly monotone bijection of $(0,\infty)$ onto $(0,1)$, so the minimizer is $-\log p$, and strict convexity in $\ell$ is the second derivative $(1-p)\e^\ell/(\e^\ell-1)^2>0$ of the summand.
\ref{prop_loglink_proj}: On $\Sy_+^m$ every term of~\eqref{lossloglink} is nonnegative, so conditioning on $(n_0,\bm z_0)$ as in the proof of Proposition~\ref{prop_target}\ref{prop_target_proj} is legitimate in $[0,\infty]$ and gives $\E[\Lcal^{\rm log}(\bm U,\xi_0)]=-\E\big[\sum_i\big(p_{0,i}\log\e^{-\ell_{0,i}}+(1-p_{0,i})\log(1-\e^{-\ell_{0,i}})\big)\big]$, and adding and subtracting $\E\big[\sum_ih(p_{0,i})\big]$ produces~\eqref{bklident}. The constant is finite because $0\le h\le\log2$. The projection statement follows from $d_{\rm Be}\ge0$ with equality only at $q=p$. For~\eqref{compsplitlog}, $d_{\rm Be}(p\|q)=\big(p\log p+(1-p)\log(1-p)\big)-\big(p\log q+(1-p)\log(1-q)\big)$, and the second bracket is affine in $p$ for every fixed $q=\e^{-\ell_{\bm U}(z)}$, so the truncation argument in the proof of Proposition~\ref{prop_target}\ref{prop_target_comp} applies, with $u\coloneqq\log(\bar p/q)-\log((1-\bar p)/(1-q))$, $v\coloneqq\log((1-\bar p)/(1-q))$, $\bar p\in(0,1)$ $\mu_1$-a.e.\ when $p_{0,i}\in(0,1)$ $\Pa$-a.s., and both sides infinite when $q=1$ on a set of positive $\mu_1$-measure. The composition part is $\epsilon^{\rm log}_{\rm comp}=\int h(\bar p)\,d\mu_1-\E\big[\sum_ih(p_{0,i})\big]\in[0,\bar n\log2]$, the gap in Jensen's inequality for the concave entropy, because $0\le h\le\log2$ and $\mu_1(\Zcal)\le\bar n$. For existence, $d_{\rm Be}(p\|\e^{-\ell})=p\ell-(1-p)\log(1-\e^{-\ell})-h(p)$, whose middle term is nonnegative, so on the cone $G(\bm U)\coloneqq\int d_{\rm Be}(\bar p\|\e^{-\ell_{\bm U}})\,d\mu_1\ge\langle\int\bar p\,\bm C\,d\mu_1,\bm U\rangle-\E[n_0]\log2\ge c_2\|\bm U\|_F-\E[n_0]\log2$ with $c_2\coloneqq\lambda_{\min}(\int\bar p\,\bm C\,d\mu_1)>0$. Lower semicontinuity follows by Fatou's lemma as for the Poisson risk, $G$ is finite somewhere on the cone when $\epsilon^{\log}_{\rm class}<\infty$, by~\eqref{compsplitlog}, and by~\eqref{bklident} and~\eqref{compsplitlog} the minimizers of $G$ are those of the population risk.
\ref{prop_loglink_score}: $\frac{d}{d\ell}\big[y\ell-(1-y)\log(1-\e^{-\ell})\big]=y-(1-y)\e^{-\ell}/(1-\e^{-\ell})=(y-\e^{-\ell})/(1-\e^{-\ell})=w_{\rm log}(\ell)(\e^{-\ell}-y)$, and $\partial\ell_{t,i}/\partial\bm U=\bm C_{t,i}$. The first-order conditions over the cone are those of the proof of Proposition~\ref{prop_target}\ref{prop_target_score}, which uses only convexity and differentiability at the minimizer, and give~\eqref{calibeqlog} after the sign of $w_{\rm log}$ is absorbed.
\end{proof}

\begin{proof}[Proof of Proposition~\ref{prop_scalediag}]
\ref{prop_scalediag_bd}: Both targets are functions of the own covariate, so $D=(\ell^{\rm Po}_\star-\bar p)-(p^{\rm log}_\star-\bar p)$ and the triangle inequality in $L^2(\mu_1)$ give $\|D\|\le\|\ell^{\rm Po}_\star-\bar p\|+\|p^{\rm log}_\star-\bar p\|\le2\max\big(\|\ell^{\rm Po}_\star-\bar p\|,\|p^{\rm log}_\star-\bar p\|\big)$. By~\eqref{klident} and~\eqref{compsplit} the Poisson risk equals $c^{\rm Po}_0+\epsilon^{\rm Po}_{\rm comp}+\int d_{\rm KL}(\bar p\|\ell_{\bm U})\,d\mu_1$ with finite constants, so its minimizer attains $\int d_{\rm KL}(\bar p\|\ell^{\rm Po}_\star)\,d\mu_1=\epsilon^{\rm Po}_{\rm dict}$. The bound $d_{\rm KL}(p\|\ell)\ge(p-\ell)^2/(2\max(p,\ell))$ of the proof of Corollary~\ref{cor_calib}, with $\max(\bar p,\ell^{\rm Po}_\star)\le\max(1,\|\ell^{\rm Po}_\star\|_\infty)$, then gives $\|\ell^{\rm Po}_\star-\bar p\|^2\le2\max(1,\|\ell^{\rm Po}_\star\|_\infty)\,\epsilon^{\rm Po}_{\rm dict}/\E[n_0]$. In the same way~\eqref{bklident}, \eqref{compsplitlog} and Pinsker's inequality $d_{\rm Be}(p\|q)\ge2(p-q)^2$ give $\|p^{\rm log}_\star-\bar p\|^2\le\epsilon^{\rm log}_{\rm dict}/(2\E[n_0])$.
\ref{prop_scalediag_null}: Both classes contain $\bar p$, so both targets equal it $\mu_1$-a.e.\ on their own scale, by Propositions~\ref{prop_target}\ref{prop_target_comp} and~\ref{prop_loglink}\ref{prop_loglink_proj}, and $D=0$ $\mu_1$-a.e. Write $\Mcal\coloneqq\E[\sum_{i\le n_0}\bm C_{0,i}\otimes\bm C_{0,i}]$. At $p\equiv c$ the Poisson score at its target is $\bm\Psi^{\rm Po}_t=\sum_i(1-y_{t,i}/c)\bm C_{t,i}$ and, by Proposition~\ref{prop_loglink}\ref{prop_loglink_score} with $w_{\rm log}=-1/(1-c)$, the log-link score is $\bm\Psi^{\rm log}_t=-\frac{c}{1-c}\bm\Psi^{\rm Po}_t$. The population Hessians are $\Mcal/c$ and $\E[\sum_i(g_{\log}')^2/(g_{\log}(1-g_{\log}))\,\bm C\otimes\bm C]=\frac{c}{1-c}\Mcal$ with $g_{\log}(\ell)=\e^{-\ell}$. The two linearizations are therefore $-c\Mcal^{-1}\bar{\bm\Psi}^{\rm Po}_T$ and $+\Mcal^{-1}\bar{\bm\Psi}^{\rm Po}_T$, and the linearization of $\hat D(z)$, $\langle\bm C(z),\cdot\rangle$ of the first plus $c\langle\bm C(z),\cdot\rangle$ of the second, vanishes identically. The remainder is $O_\Pa(T^{-1})$, because each benchmark fit satisfies $\hat{\bm U}-\bm U_\star=-\Hcal_\star^{-1}[\bar{\bm\Psi}_T]+O_\Pa(T^{-1})$, since in the expansion of the proof of Theorem~\ref{thm_clt}\ref{thm_clt_par} the pseudo-self-concordant sandwich and the central limit theorem for the bounded Hessians give $\widetilde\Hcal_T-\Hcal_\star=O_\Pa(T^{-1/2})$, and $\ell\mapsto\e^{-\ell}$ adds a second-order term of order $\|\hat{\bm U}^{\rm log}_{\rm b}-\bm U^{\rm log}_\star\|_F^2=O_\Pa(T^{-1})$.
\end{proof}

\begin{proof}[Proof of Corollary~\ref{cor_scalediag}]
\ref{prop_scalediag_alt}: The proof of Theorem~\ref{thm_clt} applies blockwise. The log-link block is convex by Proposition~\ref{prop_loglink}\ref{prop_loglink_convex}, the floor $\ell_{\bm U}\ge\underline c$ at the unselected items bounds its derivatives on the neighborhood as~\eqref{specbnd} does for the identity-link block, the selected summand being linear, and it makes the block pseudo-self-concordant in the sense of~\eqref{pseudoSC}. For the function $g$ of Lemma~\ref{lem_SC_log}, $g(x)=-\log(1-\e^{-x})$, $g''(x)=\e^x/(\e^x-1)^2$ and $g'''(x)=-\e^x(\e^x+1)/(\e^x-1)^3$, so $|g'''|/g''=\coth(x/2)\le\coth(\underline c/2)$ on $x\ge\underline c$, and~\eqref{pseudoSC} holds with $M^{\log}_\psi=B^2\coth(\underline c/2)$. Localization by Lemma~\ref{lem_loc}, the expansion, and the central limit theorem for the bounded strongly mixing stacked score \citep[Theorem~18.5.4]{ibr_lin_71} with the Cram\'er--Wold device go through. $\hat D(z)$ is a differentiable function of the pair of fits with gradient $\bm d(z)$, the second block by $\frac{d}{d\ell}(-\e^{-\ell})=\e^{-\ell}$, and the delta method gives the limit. The interval follows from the consistency of the plug-in as in Theorem~\ref{thm_clt}\ref{thm_clt_ci}. If the samples $\xi_t$ are independent across $t$, the centered stacked score at the targets is a function of $\xi_t$ alone, so it is an independent centered sequence and its autocovariances vanish.
\ref{prop_scalediag_face}: For a singular block ${\rm L}\in\{{\rm Po},\log\}$ write $\bm\Lambda^{\rm L}=\bm V^{\rm L}\bm\Gamma^{\rm L}{\bm V^{\rm L}}^\top$ for the gradient of its population risk at $\bm U^{\rm L}_\star$, $\Pcal^{\rm L}$ for the orthogonal projection of $\Sy^m$ onto $\Tcal^{\rm L}_\star\coloneqq\{\bm A:{\bm V^{\rm L}}^\top\bm A\bm V^{\rm L}=\bm 0\}$, and $\Scal^{\rm L}_\star[\bm A]\coloneqq\bm\Lambda^{\rm L}\bm A{\bm U^{\rm L}_\star}^{+}+{\bm U^{\rm L}_\star}^{+}\bm A\bm\Lambda^{\rm L}$, so that $\Pcal_{\rm j}=\Pcal^{\rm Po}\oplus\Pcal^{\log}$ and $\Scal_{\rm j}=\Scal^{\rm Po}_\star\oplus\Scal^{\log}_\star$. Treat each block as in the proof of Theorem~\ref{thm_clt}, the log-link block through Theorem~\ref{thm_log}. For a singular block the proof of part~\ref{thm_clt_sc} gives, eventually $\Pa$-a.s., $\hat{\bm U}^{\rm L}_T=\iota^{\rm L}(\hat{\bm\theta}^{\rm L}_T)$ with the chart $\iota^{\rm L}$ at $\bm U^{\rm L}_\star$, whose differential $J^{\rm L}$ is an isometry onto $\Tcal^{\rm L}_\star$, and
\[ \sqrt T\big(\hat{\bm U}^{\rm L}_T-\bm U^{\rm L}_\star\big)=-J^{\rm L}(\Kcal_\theta^{\rm L})^{-1}J^{{\rm L}*}\big[\sqrt T\bar{\bm\Psi}^{\rm L}_T\big]+o_\Pa(1),\qquad \Kcal_\theta^{\rm L}\coloneqq J^{{\rm L}*}(\Hcal^{\rm L}+\Scal^{\rm L}_\star)J^{\rm L}, \]
with $\Hcal^{\rm L}$ the Hessian of the block at its target and $\bar{\bm\Psi}^{\rm L}_T$ its average centered score. For a positive definite block the fit converges to an interior point of the cone, so it is eventually a stationary point of the empirical risk, and the expansion of the proof of part~\ref{thm_clt_par} gives the same display with $J^{\rm L}=\Ical$ and $\Kcal_\theta^{\rm L}=\Hcal^{\rm L}$. The two objectives are minimized separately, so the two expansions hold jointly. The stacked centered score is bounded, and its mixing coefficients are dominated by those of $\{\xi_t\}$, so the central limit theorem of the proof of part~\ref{prop_scalediag_alt} gives $\sqrt T(\bar{\bm\Psi}^{\rm Po}_T,\bar{\bm\Psi}^{\log}_T)\to_d\bm\Psi_{\infty,{\rm j}}\sim N(\bm 0,\Qcal_{{\rm j},{\rm LR}})$. With $J\coloneqq J^{\rm Po}\oplus J^{\log}$, an isometry onto $\Tcal^{\rm Po}_\star\times\Tcal^{\log}_\star$, and $\Kcal_\theta\coloneqq \Kcal_\theta^{\rm Po}\oplus \Kcal_\theta^{\log}$, positive definite by the curvature conditions and $\Scal_{\rm j}\succeq\bm 0$, the identity $J\Kcal_\theta^{-1}J^*=\mathcal W_{\rm j}$ holds as in the proof of Theorem~\ref{thm_clt}\ref{thm_clt_sc}, and $J^*=J^*\Pcal_{\rm j}$. The pair of fits therefore satisfies $\sqrt T(\hat{\bm U}-\bm U_\star)\to_d-\mathcal W_{\rm j}\Pcal_{\rm j}\bm\Psi_{\infty,{\rm j}}$, and the delta method with the gradient $\bm d(z)$ of part~\ref{prop_scalediag_alt}, together with the self-adjointness of $\mathcal W_{\rm j}$, gives the limit and $\sigma^2_{D,\star}(z)$. For the interval, the rank of each fit equals that of its target eventually, by Theorem~\ref{thm_clt}\ref{thm_clt_sc}, so the projections onto the range and the kernel of each fit converge to those of its target, the estimated multipliers converge as in that proof, the Hessians at the fits converge as in part~\ref{thm_clt_par}, and the stacked kernel estimator is consistent as in part~\ref{thm_clt_ci}. The centering matters because the scores at the fits have the sample means $\nabla\Rcal^{\rm L}_T(\hat{\bm U}^{\rm L}_T)\to\bm\Lambda^{\rm L}\ne\bm 0$. Compressed to the face of the fit they have sample mean zero, since $\nabla\Rcal^{\rm L}_T(\hat{\bm U}^{\rm L}_T)\hat{\bm U}^{\rm L}_T=\bm 0$ makes its range-range and range-kernel blocks vanish. The estimated operator $\hat\Pcal_{\rm j}(\hat\Hcal_{\rm j}+\hat\Scal_{\rm j})\hat\Pcal_{\rm j}$ then has the constant rank $\dim\Tcal^{\rm Po}_\star+\dim\Tcal^{\log}_\star$ for $T$ large, so its pseudo-inverse converges to $\mathcal W_{\rm j}$, and Slutsky's theorem gives the interval. The last statement is Theorem~\ref{thm_clt}\ref{thm_clt_bdy} applied to the block with vanishing multiplier.
\end{proof}

%%%%%%%%%%%%%%%%%%%%%%%%%%%%%%%%%%%%%%%%%%%%%%%%%%%%%%%%%%%%%%%%%%%%%%%%
\section{Bands at a boundary target}\label{app_faceband}

\begin{remark}[Bands at a boundary target]\label{rem_faceband}
Theorem~\ref{thm_clt}\ref{thm_clt_par}, \ref{thm_clt_curve} and~\ref{thm_clt_ci} describe the benchmark estimator, which drops the cone constraint, and Theorem~\ref{thm_clt}\ref{thm_clt_sc} the cone-constrained fit at a boundary target of rank $r<m$ under strict complementarity, $\bm\Gamma\succ\bm 0$. There the critical cone, the part of the tangent cone~\eqref{tangentcone} on which the multiplier vanishes, is the linear space $\Tcal_\star$, of dimension $r(r+1)/2+r(m-r)$, the tangent space of the manifold of positive semidefinite matrices of rank $r$, and $\Scal_\star$ is the curvature of that manifold charged by the multiplier, $\langle\bm A,\Scal_\star[\bm A]\rangle=2\tr(\bm\Gamma\bm N^\top\bm\Theta_r^{-1}\bm N)$, with $\bm N$ the range-kernel block of $\bm A$ in the eigenbasis of the target and $\bm\Theta_r$ its nonzero eigenvalues (proof of Theorem~\ref{thm_clt}\ref{thm_clt_sc}, and \citealp{shapiro97,bonnansshapiro00} and \citealp[Theorem~1]{shapiro19} for semidefinite programs). The term is positive semidefinite and vanishes when $\bm\Gamma=\bm 0$. Leaving it out misstates the variance of range-kernel directions whose eigenvalue is small, and overstates it when the score covariance is close to the Hessian. We call the band of Theorem~\ref{thm_clt}\ref{thm_clt_ci} with $\sigma^2_\star(z)=\langle\bm C(z),\Vcal_\star[\bm C(z)]\rangle$ the face band, after the face of the cone that the fit occupies. The Hessian and the score covariance are compressed to the tangent space at the fit of the manifold of its rank, which identifies $r$ for $T$ large under strict complementarity. The multiplier is estimated by $\hat{\bm Q}_\perp^\top\nabla\Rcal_T(\hat{\bm U}^{\rm Po}_T)\hat{\bm Q}_\perp$, with $\hat{\bm Q}=[\hat{\bm Q}_r\ \hat{\bm Q}_\perp]$ the eigenbasis of the fit. At an interior target the two bands coincide. At a boundary target the full-space band of the benchmark estimator is the wrong one for the reported fit. It also carries the $(m-r)(m-r+1)/2$ kernel-block directions outside $\Tcal_\star$, along which the cone-constrained fit does not move under strict complementarity (Theorem~\ref{thm_clt}\ref{thm_clt_sc}), so their contribution to the variance of $\langle\bm C(z),\bm\eta\rangle$ is spurious for the reported fit and the band is too wide (Table~\ref{tab_e8}). The diagonal specification is treated alike. Its face is the coordinates with $u_k>0$, and the curvature term is absent because the orthant is polyhedral.
\end{remark}

\begin{remark}[Penalties of order $T^{-1/2}$]\label{rem_penclt}
The limit theory is stated at $\lambda=0$. If the fits use $\lambda_T$ with $\sqrt T\lambda_T\to\kappa\ge0$, the proof of Theorem~\ref{thm_clt} carries over with the average score $\bar{\bm\Psi}^{\rm Po}_T$ replaced by $\bar{\bm\Psi}^{\rm Po}_T+\lambda_T\bm G$, centered at the unpenalized target $\bm U^{\rm Po}_\star$. The limits in parts~\ref{thm_clt_par} and~\ref{thm_clt_bdy} become $\bm\eta-\kappa{\Hcal^{\rm Po}_\star}^{-1}[\bm G]$ and $\Pi_\star\big(\bm\eta-\kappa{\Hcal^{\rm Po}_\star}^{-1}[\bm G]\big)$, and that of part~\ref{thm_clt_sc} becomes $\bm\eta_\star-\kappa\big(\Pcal(\Hcal^{\rm Po}_\star+\Scal_\star)\Pcal\big)^+\Pcal[\bm G]$, the shifted limits of \citet[Theorem~2]{knightfu00} for a linear penalty over a cone. A penalty of order $o(T^{-1/2})$ therefore leaves every limit unchanged. At $\lambda=\kappa T^{-1/2}$ and an interior target the band of part~\ref{thm_clt_ci} is off center by a fixed number $\kappa|\langle\bm C(z),{\Hcal^{\rm Po}_\star}^{-1}[\bm G]\rangle|/\sigma(z)$ of its standard errors whatever $T$. By Theorem~\ref{thm_log} the same holds for the log-link fit, and Table~\ref{tab_e14pen} shows the bias in standard errors nearly flat in $T$ at fixed $\kappa$.
\end{remark}

\clearpage
\section{Internet appendix}\label{app_internet}

\subsection{The projection band at a zero multiplier}\label{app_projband}

\begin{proposition}[Plug-in projection band]\label{prop_plug}
Under the hypotheses of Theorem~\ref{thm_clt}\ref{thm_clt_bdy} let $\Vcal\succ\bm 0$, let $d_1\ge\dots\ge d_m$ and $\hat{\bm q}_1,\dots,\hat{\bm q}_m$ be the eigenpairs of $\hat{\bm U}^{\rm Po}_T$, $\hat\Hcal_T\coloneqq\nabla^2\Rcal_T(\hat{\bm U}^{\rm Po}_T)$, the Hessian at the cone-constrained fit rather than at the benchmark of Theorem~\ref{thm_clt}\ref{thm_clt_ci}, $\hat\Qcal_T$ as in Theorem~\ref{thm_clt}\ref{thm_clt_ci} from the scores at $\hat{\bm U}^{\rm Po}_T$ centered at their mean, and $\hat\Vcal_T\coloneqq\hat\Hcal_T^{-1}\hat\Qcal_T\hat\Hcal_T^{-1}$. For $c_T\to\infty$ with $c_T/\sqrt T\to0$ let
\[ \hat s_T\coloneqq\big\langle\hat{\bm q}_m\hat{\bm q}_m^\top,\hat\Vcal_T[\hat{\bm q}_m\hat{\bm q}_m^\top]\big\rangle^{1/2},\qquad \hat r_T\coloneqq\#\big\{j\le m:\ d_j>c_T\,\hat s_T/\sqrt T\big\}, \]
$\hat{\bm V}_T\coloneqq[\hat{\bm q}_{\hat r_T+1}\ \cdots\ \hat{\bm q}_m]$, $\hat\Tcal_T$ the cone~\eqref{tangentcone} with $\hat{\bm V}_T$ in place of $\bm V$, and $\hat\Pi_T$ the $\hat\Hcal_T$-projection onto $\hat\Tcal_T$. Then $\Pa(\hat r_T=r)\to1$, and for $\bm\eta_T$ centered Gaussian with covariance $\hat\Vcal_T$ conditionally on $\Dcal_T$, the conditional law of $\hat\Pi_T(\bm\eta_T)$ given $\Dcal_T$ converges weakly in probability to that of $\Pi_\star(\bm\eta)$. With $\hat q_\beta(z)$ the conditional $\beta$-quantile of $\langle\bm C(z),\hat\Pi_T(\bm\eta_T)\rangle$, the interval $\big[\hat\ell^{\rm Po}_T(z)-\hat q_{1-\beta/2}(z)/\sqrt T,\ \hat\ell^{\rm Po}_T(z)-\hat q_{\beta/2}(z)/\sqrt T\big]$ has asymptotic coverage $1-\beta$ wherever the distribution function of $\langle\bm C(z),\Pi_\star(\bm\eta)\rangle$ is continuous and strictly increasing at both quantiles, for $m-r=1$ at every $z$ with $\ell^{\rm Po}_\star(z)>0$. Under the hypotheses of Theorem~\ref{thm_clt}\ref{thm_clt_sc}, $\Pa(\hat r_T=r)\to1$ as well.
\end{proposition}

\begin{proof}[Proof of Proposition~\ref{prop_plug}]
\emph{Rank.} Write $\hat{\bm U}$ for $\hat{\bm U}^{\rm Po}_T$, $\lambda_1\ge\dots\ge\lambda_r>0$ for the nonzero eigenvalues of $\bm U^{\rm Po}_\star$, and $\lambda_j\coloneqq0$ for $j>r$. By Weyl's inequality, $|d_j-\lambda_j|\le\opn{\hat{\bm U}-\bm U^{\rm Po}_\star}$ for every $j$, and $\opn{\hat{\bm U}-\bm U^{\rm Po}_\star}=O_\Pa(T^{-1/2})$ by the localization of Theorem~\ref{thm_main}\ref{thm_main_det}, in force eventually, and the central limit theorem of Theorem~\ref{thm_clt}\ref{thm_clt_par}. Since $\hat{\bm U}\to\bm U^{\rm Po}_\star$ and the pseudo-self-concordant sandwich gives $\hat\Hcal_T\to\Hcal^{\rm Po}_\star\succeq\alpha\Ical$ almost surely, $\hat\Hcal_T$ is invertible with probability tending to one. The sample mean of the scores at the fit converges to $\nabla\Ecal(\bm U^{\rm Po}_\star)$, which is $\bm 0$ under Theorem~\ref{thm_clt}\ref{thm_clt_bdy} and $\bm\Lambda$ under Theorem~\ref{thm_clt}\ref{thm_clt_sc}, so the estimator built from the centered scores satisfies $\hat\Qcal_T\to\Qcal_{\rm LR}$ in both cases, as in Theorem~\ref{thm_clt}\ref{thm_clt_ci}, and $\hat\Vcal_T\to\Vcal\succ\bm 0$ in probability. Since $\|\hat{\bm q}_m\hat{\bm q}_m^\top\|_F=1$, $\lambda_{\min}(\hat\Vcal_T)\le\hat s_T^2\le\lambda_{\max}(\hat\Vcal_T)$, so $\hat s_T$ is bounded away from zero and infinity in probability, which is all the rank step uses. Hence $\tau_T\coloneqq c_T\hat s_T/\sqrt T$ satisfies $\tau_T\to0$ and $\sqrt T\tau_T\to\infty$ in probability. For $j\le r$, $d_j\to\lambda_j>0$ gives $\Pa(d_j>\tau_T)\to1$, and for $j>r$, $\Pa(d_j>\tau_T)\le\Pa\big(\sqrt T\opn{\hat{\bm U}-\bm U^{\rm Po}_\star}>\sqrt T\tau_T\big)\to0$. So $\Pa(\hat r_T=r)\to1$. Under the hypotheses of Theorem~\ref{thm_clt}\ref{thm_clt_sc}, $d_j=0$ for $j>r$ eventually, $\Pa$-a.s., by that part, and the argument for $j\le r$ is unchanged.

Let $\hat{\bm P}_V\coloneqq\hat{\bm V}_T\hat{\bm V}_T^\top$ and $\bm P_V\coloneqq\bm V\bm V^\top$. On $\{\hat r_T=r\}$, $\hat{\bm V}_T$ spans the eigenspace of the $m-r$ smallest eigenvalues of $\hat{\bm U}$ and $\bm V$ that of the $m-r$ smallest eigenvalues of $\bm U^{\rm Po}_\star$, which are separated from the others by the gap $\lambda_r$, so the variant of the Davis--Kahan theorem of \citet[Theorem~2]{yu_wan_sam_15} gives $\opn{\hat{\bm P}_V-\bm P_V}\le2\|\hat{\bm U}-\bm U^{\rm Po}_\star\|_F/\lambda_r\to0$ in probability.

For an orthogonal projection $\bm P'$ of rank $m-r$ write $\Tcal_{\bm P'}\coloneqq\{\bm A\in\Sy^m:\bm P'\bm A\bm P'\succeq\bm 0\}$, so that $\hat\Tcal_T=\Tcal_{\hat{\bm P}_V}$ on $\{\hat r_T=r\}$ and $\Tcal(\bm U^{\rm Po}_\star)=\Tcal_{\bm P_V}$. If $\bm P_n\to\bm P_V$ then $\Tcal_{\bm P_n}\to \Tcal_{\bm P_V}$ in the sense of Painlev\'e--Kuratowski. A limit of $\bm A_n\in \Tcal_{\bm P_n}$ lies in $\Tcal_{\bm P_V}$ by continuity of $\bm A\mapsto\bm P_V\bm A\bm P_V$, and for $\bm A\in \Tcal_{\bm P_V}$ the matrices $\bm A_n\coloneqq\bm A+2\opn{\bm A}\opn{\bm P_n-\bm P_V}\bm I_m$ lie in $\Tcal_{\bm P_n}$ and converge to $\bm A$, since for a unit vector $\bm w$ in the range of $\bm P_n$, $\|\bm w-\bm P_V\bm w\|_2=\|(\bm P_n-\bm P_V)\bm w\|_2\le\opn{\bm P_n-\bm P_V}$ and $\bm w^\top\bm A\bm w\ge(\bm P_V\bm w)^\top\bm A(\bm P_V\bm w)-2\opn{\bm A}\|\bm w-\bm P_V\bm w\|_2\ge-2\opn{\bm A}\opn{\bm P_n-\bm P_V}$. Now let $(\bm P_n,\Hcal_n,\bm\eta_n)\to(\bm P_V,\Hcal,\bm\eta)$ with $\Hcal_n,\Hcal\succ\bm 0$. The functions $f_n(\bm A)\coloneqq\langle\bm A-\bm\eta_n,\Hcal_n[\bm A-\bm\eta_n]\rangle+\delta_{\Tcal_{\bm P_n}}(\bm A)$ epi-converge to $f(\bm A)\coloneqq\langle\bm A-\bm\eta,\Hcal[\bm A-\bm\eta]\rangle+\delta_{\Tcal_{\bm P_V}}(\bm A)$, as the sum of a continuously convergent sequence and an epi-convergent one \citep[Theorem~7.46]{roc_wet_98}, they are eventually level-bounded, since $f_n(\bm A)\ge\lambda_{\min}(\Hcal_n)\|\bm A-\bm\eta_n\|_F^2$, and $f$ has the unique minimizer $\Pi(\bm P_V,\Hcal,\bm\eta)\coloneqq\argmin_{\bm A\in \Tcal_{\bm P_V}}\langle\bm A-\bm\eta,\Hcal[\bm A-\bm\eta]\rangle$, so $\argmin f_n\to\argmin f$ \citep[Theorem~7.33]{roc_wet_98}. The map $\Pi$ is therefore continuous, and $\Pi(\bm P_V,\Hcal^{\rm Po}_\star,\bm\eta)=\Pi_\star(\bm\eta)$.

Let $f$ be bounded and continuous on $\Sy^m$, let $\bm g$ be a standard Gaussian element of $\Sy^m$ independent of $\Dcal_T$, and let $\hat\Vcal_T^{1/2}$ be the positive square root, so that, conditionally on $\Dcal_T$, $\bm\eta_T$ has the law of $\hat\Vcal_T^{1/2}[\bm g]$ and $\E[f(\hat\Pi_T(\bm\eta_T))\mid\Dcal_T]=\E_{\bm g}\big[f\big(\Pi(\hat{\bm P}_V,\hat\Hcal_T,\hat\Vcal_T^{1/2}[\bm g])\big)\big]$ on $\{\hat r_T=r\}$. Every subsequence has a further subsequence along which $\hat r_T=r$ eventually and $(\hat{\bm P}_V,\hat\Hcal_T,\hat\Vcal_T)\to(\bm P_V,\Hcal^{\rm Po}_\star,\Vcal)$ almost surely, and along it the continuity of $\Pi$ and of the square root and dominated convergence give $\E[f(\hat\Pi_T(\bm\eta_T))\mid\Dcal_T]\to\E_{\bm g}[f(\Pi(\bm P_V,\Hcal^{\rm Po}_\star,\Vcal^{1/2}[\bm g]))]=\E[f(\Pi_\star(\bm\eta))]$ almost surely, which is the convergence in probability claimed. For the interval, $\sqrt T(\hat\ell^{\rm Po}_T(z)-\ell^{\rm Po}_\star(z))\to_dW\coloneqq\langle\bm C(z),\Pi_\star(\bm\eta)\rangle$ by Theorem~\ref{thm_clt}\ref{thm_clt_bdy}, the conditional law of $\langle\bm C(z),\hat\Pi_T(\bm\eta_T)\rangle$ converges weakly in probability to that of $W$ by the continuous mapping theorem, so that $\hat q_\beta(z)\to q_\beta$ in probability whenever the distribution function of $W$ is continuous and strictly increasing at its $\beta$-quantile $q_\beta$, and $\Pa\big(\hat q_{\beta/2}(z)\le\sqrt T(\hat\ell^{\rm Po}_T(z)-\ell^{\rm Po}_\star(z))\le\hat q_{1-\beta/2}(z)\big)\to\Pa(q_{\beta/2}\le W\le q_{1-\beta/2})=1-\beta$ by Slutsky's theorem. For $m-r=1$ let $\bm v$ span the kernel and $\bm n\coloneqq\bm v\bm v^\top$. The cone is the half-space $\{\langle\bm n,\bm A\rangle\ge0\}$, so $\Pi_\star(\bm\eta)=\bm\eta$ on $\{\langle\bm n,\bm\eta\rangle\ge0\}$ and $\Pi_\star(\bm\eta)=\bm\eta-\langle\bm n,\bm\eta\rangle\,\Hcal^{\rm Po\,-1}_\star[\bm n]/\langle\bm n,\Hcal^{\rm Po\,-1}_\star[\bm n]\rangle$ otherwise, and $W=\langle\bm C(z),\bm\eta\rangle$ on the first event and $W=\langle\bm C(z)-\gamma\bm n,\bm\eta\rangle$ on the second, with $\gamma\coloneqq\langle\bm C(z),\Hcal^{\rm Po\,-1}_\star[\bm n]\rangle/\langle\bm n,\Hcal^{\rm Po\,-1}_\star[\bm n]\rangle$. Since $\bm C(z)=\bm\phi(z)^\top\bm\phi(z)$ has rank one, $\bm C(z)$ is proportional to $\bm n$ if and only if $\bm\phi(z)^\top$ is proportional to $\bm v$, that is, if and only if $\ell^{\rm Po}_\star(z)=\bm\phi(z)\bm U^{\rm Po}_\star\bm\phi(z)^\top=0$. If $\ell^{\rm Po}_\star(z)>0$, then $\bm C(z)$ and $\bm n$ are linearly independent, as are $\bm C(z)-\gamma\bm n$ and $\bm n$, so on each event the pair $(W,\langle\bm n,\bm\eta\rangle)$ is a nondegenerate Gaussian vector of $\R^2$ restricted to a half-plane, by $\Vcal\succ\bm 0$, and $W$ has a density that is positive on $\R$. If $\bm\phi(z)^\top=a\bm v$, then $\bm C(z)=a^2\bm n$, $\gamma=a^2$ and $W=a^2\max(\langle\bm n,\bm\eta\rangle,0)$, which has the atom $\Pa(\langle\bm n,\bm\eta\rangle\le0)=1/2$ at zero.
\end{proof}

The fits of the application sit on the boundary of the cone (Section~\ref{sec_app_data}). Two designs on the dictionary $\{1,\phi_1,\dots,\phi_4\}$, with independent inclusions in the baseline design, put a target of rank $4$ there. In the first, $p^\circ$ is a sum of squares with no weight on $\phi_4$, so the model is correct and the multiplier vanishes, $\bm\Gamma=\bm 0$. In the second, $p^\circ$ is that curve less $0.026\,\phi_4^2$, not a sum of squares, and its projection has multiplier $\bm\Gamma=0.91\succ\bm 0$, so strict complementarity holds. Table~\ref{tab_e8} compares four intervals. With $\bm\Gamma\succ\bm 0$ the face band is right. From $T=1000$ on its standard error matches the realized dispersion and it covers at $0.93$ to $0.96$, while dropping the curvature makes it $6$ percent too wide and the full-space band $14$ to $15$ percent, against $16$ and $29$ percent at $T=250$, as Remark~\ref{rem_faceband} anticipates. The benchmark estimates the indefinite matrix that represents $p^\circ$ and misses the cone-constrained target. With $\bm\Gamma=\bm 0$ the face band covers from $T=1000$ on and falls to $0.92$ at one point at $T=250$, where the estimated multiplier adds curvature the target lacks. The rank of the fit does not settle at that of the target. Under $\bm\Gamma=\bm 0$, in the separate run of Table~\ref{tab_e8thr}, it is $4$ in $69$, $77$ and $60$ percent of the replications at $T=250$, $1000$ and $4000$, $3$ in $24$, $3$ and $0$ percent, and $5$ in $7$, $19$ and $40$ percent, the last share rising toward the limit probability one half of a positive kernel block of $\Pi_\star(\bm\eta)$, as Theorem~\ref{thm_clt}\ref{thm_clt_bdy} predicts.

The rule of Proposition~\ref{prop_plug} with $c_T=\sqrt{\log T}$ selects the target's rank $4$ in $0.4$, $20$ and $92$ percent of the replications at $T=250$, $1000$ and $4000$, and its band covers at $0.89$ to $0.91$, $0.92$ to $0.94$ and $0.93$ to $0.95$, the last in agreement with an oracle that projects with the true kernel (Table~\ref{tab_e8thr}). The design is a hard case for any rank rule, since the fourth eigenvalue of the target, $0.016$, is within the sampling scale of the kernel block up to $T=1000$, so at small $T$ the rule treats range directions as kernel directions and imposes a constraint the target does not satisfy. For a second zero-multiplier target with the wider spectrum $0.179$, $0.108$, $0.066$, $0.052$ and $0$, the band covers at $0.80$ to $0.86$ at $T=250$ and at $0.93$ to $0.96$ from $T=1000$ on, while under strict complementarity the face band of Table~\ref{tab_e8} is the band to use. In the first zero-multiplier design at $T=250$ the face band covers at $0.92$ to $0.97$ (Table~\ref{tab_e8}), the projection band at $0.89$ to $0.91$, and the benchmark band, valid at a zero multiplier for its own estimator, at $0.93$ to $0.96$, so the evidence at $T=250$ is design-specific.

\begin{table}[!tbp]
\centering
\caption{Projection bands at a zero-multiplier target, $500$ replications, baseline design, a separate run of the designs of Table~\ref{tab_e8}. Entries are the range over $z\in\{-0.5,0,0.5\}$ of the coverage of nominal $95$ percent intervals. Fit rank gives the shares of fits of rank $3$, $4$ and $5$ at the relative tolerance $10^{-4}$. Threshold is the band of Proposition~\ref{prop_plug} with $c_T=\sqrt{\log T}$, with the shares of selected ranks at most $2$, exactly $3$ and exactly $4$. Oracle projects with the true kernel and the Hessian and score covariance at the target. The rows for $\bm\Gamma\succ\bm 0$ are the second design of Table~\ref{tab_e8}, in this run, where the limit is Gaussian on the face and the projection band is not the one to use. Wide is a second zero-multiplier target with eigenvalues $0.179$, $0.108$, $0.066$, $0.052$ and $0$, obtained from the first by raising the ridge from $0.6$ to $4$.}
\label{tab_e8thr}
\begin{footnotesize}
\setlength{\tabcolsep}{3.5pt}
\begin{tabular}{lrcccc}
\toprule
 & & fit rank & \multicolumn{2}{c}{threshold} & oracle \\
\cmidrule(lr){4-5}
 & $T$ & $3/4/5$ & coverage & rank $\le2/3/4$ & coverage \\
\midrule
$\bm\Gamma=\bm 0$ & $250$ & $0.24/0.69/0.07$ & $0.89$--$0.91$ & $0.50$$/0.50/0.00$ & $0.97$--$0.98$ \\
 & $1000$ & $0.03/0.77/0.19$ & $0.92$--$0.94$ & $0.08$$/0.72/0.20$ & $0.96$--$0.97$ \\
 & $4000$ & $0.00/0.60/0.40$ & $0.93$--$0.95$ & $0.00/0.08/0.92$ & $0.93$--$0.95$ \\
\addlinespace
$\bm\Gamma\succ\bm 0$ & $250$ & $0.55/0.45/0.00$ & $0.92$--$0.94$ & $0.30/0.70$$/0.00$ & $0.96$--$0.98$ \\
 & $1000$ & $0.25/0.75/0.00$ & $0.94$--$0.95$ & $0.01$$/0.95/0.05$ & $0.94$--$0.99$ \\
 & $4000$ & $0.07/0.93/0.00$ & $0.93$--$0.96$ & $0.00/0.39/0.61$ & $0.95$--$0.99$ \\
\addlinespace
wide, $\bm\Gamma=\bm 0$ & $250$ & $0.18/0.73/0.09$ & $0.80$--$0.86$ & $0.46/0.41$$/0.13$ & $0.96$--$0.98$ \\
 & $1000$ & $0.04/0.69/0.26$ & $0.93$--$0.96$ & $0.07/$$0.19$$/0.74$ & $0.94$--$0.97$ \\
 & $4000$ & $0.00/0.60/0.40$ & $0.93$--$0.96$ & $0.00/0.01/0.99$ & $0.94$--$0.96$ \\
\bottomrule
\end{tabular}
\end{footnotesize}
\end{table}

\begin{table}[!tbp]
\centering
\caption{Coverage at a boundary target, $500$ replications, baseline design. Entries are the range over $z\in\{-0.5,0,0.5\}$ of the coverage of nominal $95$ percent intervals, the mean standard error over the Monte Carlo standard deviation, averaged over the points, and the share of fits of the target's rank $4$, the rank read off the fit with a relative eigenvalue tolerance of $10^{-4}$. Face is the band of Remark~\ref{rem_faceband}, no curv.\ is that band without the curvature of the face, full is the full-space band, and benchmark is the unconstrained estimator.}
\label{tab_e8}
\begin{footnotesize}
\setlength{\tabcolsep}{2.5pt}
\begin{tabular}{lrcccccccc}
\toprule
 & & \multicolumn{4}{c}{coverage} & \multicolumn{3}{c}{s.e./MC s.d.} & rank \\
\cmidrule(lr){3-6}\cmidrule(lr){7-9}
 & $T$ & face & no curv. & full & benchmark & face & no curv. & full & $4$ \\
\midrule
$\bm\Gamma\succ\bm 0$ & $250$ & $0.94$--$0.96$ & $0.96$--$0.98$ & $0.98$--$0.99$ & $0.75$--$0.95$ & $1.05$ & $1.16$ & $1.29$ & $0.46$ \\
 & $1000$ & $0.93$--$0.95$ & $0.94$--$0.97$ & $0.97$--$0.98$ & $0.26$--$0.94$ & $0.98$ & $1.06$ & $1.15$ & $0.73$ \\
 & $4000$ & $0.95$--$0.96$ & $0.96$--$0.96$ & $0.97$--$0.99$ & $0.00$--$0.94$ & $0.99$ & $1.06$ & $1.14$ & $0.94$ \\
\addlinespace
$\bm\Gamma=\bm 0$ & $250$ & $0.92$--$0.97$ & $0.94$--$0.97$ & $0.96$--$0.98$ & $0.93$--$0.96$ & $1.04$ & $1.13$ & $1.21$ & $0.72$ \\
 & $1000$ & $0.93$--$0.95$ & $0.94$--$0.96$ & $0.95$--$0.97$ & $0.94$--$0.96$ & $1.00$ & $1.04$ & $1.09$ & $0.73$ \\
 & $4000$ & $0.95$--$0.96$ & $0.95$--$0.96$ & $0.95$--$0.97$ & $0.94$--$0.96$ & $1.00$ & $1.01$ & $1.06$ & $0.65$ \\
\bottomrule
\end{tabular}
\end{footnotesize}
\end{table}

%%%%%%%%%%%%%%%%%%%%%%%%%%%%%%%%%%%%%%%%%%%%%%%%%%%%%%%%%%%%%%%%%%%%%%%%

\subsection{Identification of the matrix and of the curve}\label{app_ident}

The curvature condition~\eqref{curv} fails whenever a nonzero $\bm A\in\Sy^m$ generates the zero curve, $\bm\phi(z)\bm A\bm\phi(z)^\top=0$ for almost every covariate value, since every term of $\langle\bm A,\Hcal^{\rm Po}_\star[\bm A]\rangle$, an expectation of terms $(\bm\phi(z_{0,i})\bm A\bm\phi(z_{0,i})^\top)^2$ with positive weights as in~\eqref{curvexplicit}, then vanishes. Linear independence of the dictionary does not exclude this, because the $M$ products $\phi_k\phi_l$ may be linearly dependent on the support of the covariates. When they are, $\bm U$ is not identified but the curve is, since the population risk is strictly convex in the curve. Let $\Ncal\coloneqq\{\bm A\in\Sy^m:\bm\phi(z)\bm A\bm\phi(z)^\top=0\text{ for }\mu_1\text{-a.e.\ }z\}$. The risk is invariant under $\bm U\mapsto\bm U+\bm A$, $\bm A\in\Ncal$, wherever both points lie in $\Ocal$, and for continuous $\bm\phi$, $\bm C(z)\perp\Ncal$ on the support of $\mu_1$. If~\eqref{curv} holds on $\Ncal^\perp$, parts~\ref{thm_clt_par}, \ref{thm_clt_curve} and~\ref{thm_clt_ci} of Theorem~\ref{thm_clt} hold for the projection of the benchmark estimator onto $\Ncal^\perp$ and for its curve, by the proof of Theorem~\ref{thm_clt}\ref{thm_clt_par} applied in the quotient space $\Sy^m/\Ncal$, for any measurable selection of the benchmark minimizer, with ${\Hcal^{\rm Po}_\star}^{-1}$ replaced by the inverse of $\Hcal^{\rm Po}_\star$ restricted to $\Ncal^\perp$, equivalently its Moore--Penrose inverse. The cone-constrained fit is not covered by this argument, since the cone does not decompose along $\Ncal$. The fitted curve need not be unique, since the empirical loss is strictly convex in the curve only at the events. The design of Section~\ref{sec_sim_design} at $m=8$ is an instance to working precision. Its products are linearly independent, but its design operator is singular in double precision.

For the dictionaries of Section~\ref{sec_app} the question is empirical, but the answer has a structure. For the grid and lin the kernel part of the dictionary consists of affine functions of the covariate vector itself, $\bm\phi(z)=(1,z^\top/\sqrt d)\bm W$ with $d=121$ and $\bm W\in\R^{(d+1)\times m}$ of full column rank, so that $\bm\phi\bm A\bm\phi^\top$ is the quadratic form $\bm W\bm A\bm W^\top$ evaluated at $(1,z/\sqrt d)$, and a zero curve is a quadric through every covariate value of the window. The data force none, but covariates common to the names of a month can come close to one. Their products are functions of the month alone, evaluated at the $97$ months of a window, and nearly collinear when they move slowly. The characteristics, rank-transformed within the month, contain no such covariate. The three state summaries of Section~\ref{sec_app_data}, appended to the dictionary rather than to $z$, are common to the names of a month, but the ten products of the constant and the summaries are far from collinear. At the median window their singular values fall only to $5\cdot10^{-3}$ of the largest, and to $1.6\cdot10^{-3}$ at worst. A direction $\bm A$ of the design operator with eigenvalue $\epsilon$ moves the curve by $\sqrt\epsilon\,\|\bm A\|_F$ in root mean square over the window, so where the operator is poorly conditioned, what is poorly determined is the matrix along those directions, not the curve. The dictionary lin spans all affine functions of $z$, so its zero curves are the quadrics through the characteristic vectors of a window. We have not formed its design operator, of order $8{,}128$ with the summaries, and read lin for the curve only.

The dictionaries arc and half are not affine, and they have exact null directions. In arc the features $z_j/\sqrt d$ and $\sqrt{1-z_j^2}/\sqrt d$ satisfy $(z_j/\sqrt d)^2+(\sqrt{1-z_j^2}/\sqrt d)^2=1/d$, a multiple of the square of the constant, and half and the lifted summaries satisfy the analogous identities. This gives $d$ null directions plus one per summary on complete records, and under the masked dictionary at least those of the coordinates observed for every name and of the summaries. Along them the fitted curve is unchanged at every covariate value, so predictions do not depend on the representation the solver returns, while the matrix, its rank and the bands of Section~\ref{sec_theory} are not identified. We use these dictionaries, and lg and lgz, whose Gaussian sections are not affine either, for prediction only.

The pivoted Cholesky dictionaries, $m\le40$, escape this structure exactly. A quadratic form $\bm W\bm A\bm W^\top$ that is a function of the month alone has a vanishing block $\bm W_c\bm A\bm W_c^\top$ on the characteristics, with $\bm W_c$ the $121$ rows of $\bm W$ that multiply them, provided the projections $\bm W_c^\top z$ of the characteristic vectors of the names of a month do not lie on a common quadric of $\R^m$. A quadric of $\R^m$ has at most $m(m+3)/2+1$ coefficients, $861$ at $m=40$, against about $4{,}600$ names a month, so this is a general-position condition on the data, which the positive definiteness of the design operator reported in Section~\ref{sec_app_data} checks at $m\le20$. Up to an invertible change of basis, the columns of $\bm W_c$ are the characteristic parts of the $m$ selected training covariates, $m\le40$ vectors in $\R^{121}$ that are linearly independent in general position. Hence $\bm W_c$ has full column rank, $\bm W_c\bm A\bm W_c^\top=0$ forces $\bm A=0$, and no nonzero direction of the kernel part generates a curve that depends on the month alone. Of the five summaries appended to it, the two name-level ones vary across the names of a month and are not affine in $z$, and the three state summaries contribute the month-level products just described. The positive definiteness of the design operator reported in Section~\ref{sec_app_data}, computed with all five at $m\le20$, confirms that these grid dictionaries are identified at $m\le20$. At $m=40$ identification rests on the general-position argument alone.

\subsection{Simulation in the regime of the application}\label{app_simapp}\label{sec_sim_app}

This experiment mimics the application with a known truth. Each period has $3{,}450$ to $5{,}750$ items with $d=10$ correlated covariates, and $q_A$ is a sum of two squares of affine functions of $z$. In design~A the margins are $p_t=\min\{c_t\,g_A,0.95\}$, with $g_A\propto q_A$ of mean $0.0025$, the failure rate of Section~\ref{sec_app_data}, fitted by the identity-link fit. In design~B they are $p_t=\min\{c_t\,\e^{-c_0-6q_A},0.95\}$ with the same mean, fitted by the log-link fit. A $97$-period window holds about $1{,}150$ events, against $940$ in the application. The rate is stationary, $c_t=1$, or drifts as an AR(1) in $\log c_t$. The dictionaries are the linear kernel with $m=11$, exact and unpenalized, and the linear plus a Gaussian kernel with $m=40$ at $\lambda=T^{-1/2}$. Each of $50$ replications runs $100$ test periods, each fitted on the preceding $97$, and bins the test items by the true margin and, as in the application, by the fitted probability. The covariates are $z_j=2\Phi(\sqrt{0.3}\,f+\sqrt{0.7}\,\varepsilon_j)-1$ with a factor $f$ common to the coordinates of an item, so that they are uniform with correlation $0.29$, and under drift $\log c_t$ is an AR(1) with coefficient $0.95$ and innovation standard deviation $0.1$. The cap at $0.95$ never binds, since the largest margin is $0.028$ in the stationary design~A and $0.0052$ in design~B, a far narrower range than the three orders of magnitude of the fitted probabilities in the application (Section~\ref{sec_theory_log}). The bands are the face bands of Remark~\ref{rem_faceband} at the rank of each fit. At $m=11$ the targets are correctly specified and unpenalized, so the multiplier vanishes. For the log-link fit of design~B, Theorem~\ref{thm_clt}\ref{thm_clt_bdy} applies through Theorem~\ref{thm_log}. For the identity-link fit of design~A the margins approach zero near the null set of $g_A$, (P3) fails and the theorem does not apply. The face band is an approximation in both designs. At $m=40$ the fits are penalized at $\lambda=T^{-1/2}$, so by Remark~\ref{rem_penclt} the penalty shifts their limit off the target, and the dictionary contains no nonzero affine function, so the class approximates the margins rather than representing them.

Table~\ref{tab_e14} reports the results. With a stationary rate and no penalty the bands are close to right in bins of the true margin, except for the identity-link fit in the lowest three bins. The log-link fit covers at $88$ to $96$ percent, with standardized errors whose standard deviation is $1.02$ to $1.20$, so its bands are about ten percent too narrow. There, near the null set of $g_A$, the margins approach zero, (P3) fails, and a fit that cannot fall below zero errs upward by up to $2.9$ standard errors. At $m=40$ both fits miss in bins of the true margin, by up to $7.2$ standard errors for the identity-link fit and $9.4$ for the log-link fit. Binning by the fitted probability adds selection, which moves the coverage of the log-link bands to $2$ to $90$ percent. A drifting rate leaves the bands centered but too narrow for the margins of the test period, which the drift moves away from the target of the window.

The penalty $\lambda=T^{-1/2}$ biases the log-link fit, by $14$ standard errors at $m=11$, the shift of Remark~\ref{rem_penclt}, while the identity-link fit barely moves. Table~\ref{tab_e14pen} confirms the scale $\lambda T/\sqrt N$ of Section~\ref{sec_theory_log}, with a constant of $49$ to $65$. At $m=11$, holding the bias at a quarter of a standard error calls for
\begin{equation}\label{lamlog}
   \lambda=\frac{\sqrt N}{200\,T},
\end{equation}
about sixty times smaller than $T^{-1/2}$ in the application. The misses of the identity-link grid fit in Figure~\ref{fig_reliability_grid}, factors of $2.1$ to $3.4$ in the fourth to the ninth decile, exceed by far what the selection, the null set or a drifting rate produce in the runs of Table~\ref{tab_e14}, where the identity-link fit of design~A stays within eight percent of the realized frequency in every bin of its fit, drifting rate included, apart from a miss of $22$ percent in the lowest bin at $m=40$. They also lie in other bins, which points to the class gap.

\begin{table}[!tp]
\centering
\caption{The regime of the application: $d=10$, $3{,}450$ to $5{,}750$ items per period, event rate $0.25$ percent, rolling $97$-period windows, $100$ test periods, $50$ replications per design, each fitted by the link of its class. Dictionaries: linear kernel, $m=11$, unpenalized, and linear plus Gaussian kernel, $m=40$, at $\lambda=T^{-1/2}$. With the test items binned by the true margin, the columns report, over the ten pooled bins, the range of the coverage of the nominal $95$ percent band for the fitted bin mean, the largest absolute mean over the replications of its standardized difference from the true bin mean, and the range of the standard deviation of that difference. With the items binned by the fitted probability, as in the application, they report the ranges of the coverage of the band for the realized frequency (total) and of the band for the fitted bin mean (fit). RMSE of the fitted probability in units of $10^{-3}$, and seconds per fit. Bin means of different test periods are pooled with covariance equal to the product of their standard errors times the overlap share of their windows, $(97-|a-c|)/97$ for periods $a$ and $c$, an approximation, since the sandwich of one window does not carry the covariance across windows. The Monte Carlo standard error of a coverage near $0.95$ is $0.03$.}
\label{tab_e14}
\begin{footnotesize}
\setlength{\tabcolsep}{2.5pt}
\begin{tabular}{lccccccc}
\toprule
 & \multicolumn{3}{c}{bins of the true margin} & \multicolumn{2}{c}{bins of the fit} & & \\
\cmidrule(lr){2-4}\cmidrule(lr){5-6}
 & coverage & max $|$mean$|$ & s.d. & total & fit & RMSE & seconds \\
\midrule
\multicolumn{8}{l}{\emph{A.\ Margins in the identity class, identity-link fit}} \\
\addlinespace
stationary, $m=11$ & $0.08$--$0.98$ & $2.85$ & $0.65$--$1.02$ & $0.92$--$1.00$ & $0.88$--$1.00$ & $0.32$ & $0.020$ \\
stationary, $m=40$ & $0.00$--$0.94$ & $7.20$ & $0.59$--$1.24$ & $0.86$--$1.00$ & $0.56$--$1.00$ & $0.74$ & $0.147$ \\
drifting, $m=11$ & $0.10$--$0.64$ & $2.83$ & $0.67$--$2.84$ & $0.70$--$0.96$ & $0.48$--$0.88$ & $1.05$ & $0.021$ \\
drifting, $m=40$ & $0.00$--$0.58$ & $6.84$ & $0.71$--$2.97$ & $0.64$--$0.96$ & $0.38$--$0.78$ & $1.29$ & $0.151$ \\
\midrule
\multicolumn{8}{l}{\emph{B.\ Margins in the log class, log-link fit}} \\
\addlinespace
stationary, $m=11$ & $0.88$--$0.96$ & $0.24$ & $1.02$--$1.20$ & $0.50$--$0.96$ & $0.02$--$0.90$ & $0.60$ & $4.2$ \\
stationary, $m=40$ & $0.00$--$0.90$ & $9.41$ & $0.77$--$1.49$ & $0.14$--$0.98$ & $0.00$--$0.86$ & $0.94$ & $13.2$ \\
drifting, $m=11$ & $0.54$--$0.86$ & $0.40$ & $1.33$--$2.87$ & $0.56$--$0.98$ & $0.28$--$0.74$ & $1.04$ & $4.1$ \\
drifting, $m=40$ & $0.00$--$0.46$ & $8.45$ & $1.11$--$2.92$ & $0.32$--$0.86$ & $0.00$--$0.60$ & $1.28$ & $14.2$ \\
\bottomrule
\end{tabular}
\end{footnotesize}
\end{table}

\begin{table}[!tp]
\centering
\caption{Penalty bias of the log-link fit in single stationary windows of design~B at $m=11$, $16$ replications per row, at $\lambda=\kappa T^{-1/2}$, measured as the largest absolute difference over the ten bins of the true margin between the mean standardized error of the fitted bin mean at $\lambda$ and at $\lambda=0$, in standard errors of the band. Items per period on average, and $N$ the expected number of events per window. In the first five rows $\lambda T/\sqrt N=\kappa/\sqrt{np}$, with $n$ items per period and rate $p$, does not change with $T$, so the last three rows test the scale. The last column divides the bias at $\kappa=1/25$ by $\lambda T/\sqrt N$. The first five rows vary the window length, the next two the items per period, the last the failure rate.}
\label{tab_e14pen}
\begin{footnotesize}
\begin{tabular}{rrrrccc}
\toprule
 & items per & rate & & \multicolumn{2}{c}{bias at $\kappa$} & \\
\cmidrule(lr){5-6}
$T$ & period & (percent) & $N$ & $1$ & $1/25$ & ratio \\
\midrule
$24$ & $4{,}600$ & $0.25$ & $276$ & $12.6$ & $0.77$ & $65$ \\
$48$ & $4{,}600$ & $0.25$ & $552$ & $14.1$ & $0.69$ & $58$ \\
$97$ & $4{,}600$ & $0.25$ & $1{,}116$ & $14.3$ & $0.61$ & $52$ \\
$194$ & $4{,}600$ & $0.25$ & $2{,}231$ & $14.6$ & $0.61$ & $52$ \\
$388$ & $4{,}600$ & $0.25$ & $4{,}462$ & $14.2$ & $0.58$ & $49$ \\
\addlinespace
$97$ & $2{,}300$ & $0.25$ & $558$ & $17.0$ & $0.96$ & $58$ \\
$97$ & $9{,}200$ & $0.25$ & $2{,}231$ & $10.8$ & $0.43$ & $52$ \\
\addlinespace
$97$ & $4{,}600$ & $1.00$ & $4{,}462$ & $7.3$ & $0.30$ & $51$ \\
\bottomrule
\end{tabular}
\end{footnotesize}
\end{table}

\subsubsection{The penalty for prediction}\label{sec_sim_lampred}

\begin{table}[!tp]
\centering
\caption{The penalty of the log-link fit in simulation. Design~B in stationary windows, $128$ replications per row, about $12$ events per period. Excess risk is the mean Bernoulli divergence $d_{\rm Be}(p\,\|\,\hat p)$ of the fit from the true margins per item-period over five new test periods, in units of $10^{-4}$. Rule: $\lambda=T^{-1/2}$, $\sqrt T$ phantom events per direction (Section~\ref{sec_theory_log}). Column~\eqref{lamlog}: rule~\eqref{lamlog}, calibrated at $m=11$, $\sqrt N/200$ phantom events, $0.12$ to $0.24$. $\lambda T=1$, $2$, $4$: that many phantom events. TIC: in each window, the penalty with the smallest Takeuchi information criterion, in the form that \citet{shibata89} derives for penalized fits. The criterion is the in-sample loss plus $\tr(\hat\Hcal^+\hat\Qcal)/T$ on the face of the fit, where $\hat\Hcal$ is the face operator of Remark~\ref{rem_faceband}, curvature term included, and $\hat\Qcal$ is the Bartlett covariance of the face scores at lag three, among all penalties fitted, the multiples $xT^{-1/2}$ for $x$ from $4$ down to $1/256$ in factors of four, $\lambda T\in\{0.5,1,2,4\}$, rule~\eqref{lamlog} and three further data-driven penalties. The dictionaries at $m=11$ and $m=40$ are those above, and at $m=250$ the linear plus Gaussian kernel at $250$ pivoted Cholesky sections. Fits that stopped before convergence, at most $3$ of the $128$ of a row at rule~\eqref{lamlog}, $6$ at the smallest multiple and none at the rule or at $\lambda T\in\{1,2,4\}$, are scored as they stand.}
\label{tab_lamstudy}
\begin{footnotesize}
\begin{tabular}{lcccccc}
\toprule
 & rule & \eqref{lamlog} & $\lambda T=1$ & $\lambda T=2$ & $\lambda T=4$ & TIC \\
\midrule
$m=11$, $T=97$ & $1.77$ & $0.54$ & $0.42$ & $0.42$ & $0.68$ & $0.43$ \\
$m=40$, $T=48$ & $2.25$ & $4.13$ & $1.54$ & $1.33$ & $1.58$ & $1.44$ \\
$m=40$, $T=97$ & $1.56$ & $2.06$ & $0.92$ & $0.77$ & $0.85$ & $0.78$ \\
$m=40$, $T=194$ & $1.06$ & $1.08$ & $0.58$ & $0.46$ & $0.46$ & $0.46$ \\
$m=250$, $T=97$ & $1.02$ & $5.33$ & $1.42$ & $0.82$ & $0.70$ & $0.72$ \\
\bottomrule
\end{tabular}
\end{footnotesize}
\end{table}

Table~\ref{tab_lamstudy} reports the excess risk of the log-link fit of design~B in stationary windows up to $m=250$. One to four phantom events per direction are best in every row, the risk-optimal count rising with the dictionary from about two at $m=40$ to four at $m=250$. The rule $\lambda=T^{-1/2}$ has up to $4.2$ times the smallest risk at $m=11$ and up to $2.3$ times in the richer dictionaries, and rule~\eqref{lamlog} overfits the richer dictionaries. Two phantom events, however, bias the bin means by $1.6$ to $3.0$ standard errors, so prediction and bands call for different penalties.

\subsection{Missing covariates}\label{app_missing}

\begin{figure}[!tbp]
\centering
\figorbox[width=\textwidth]{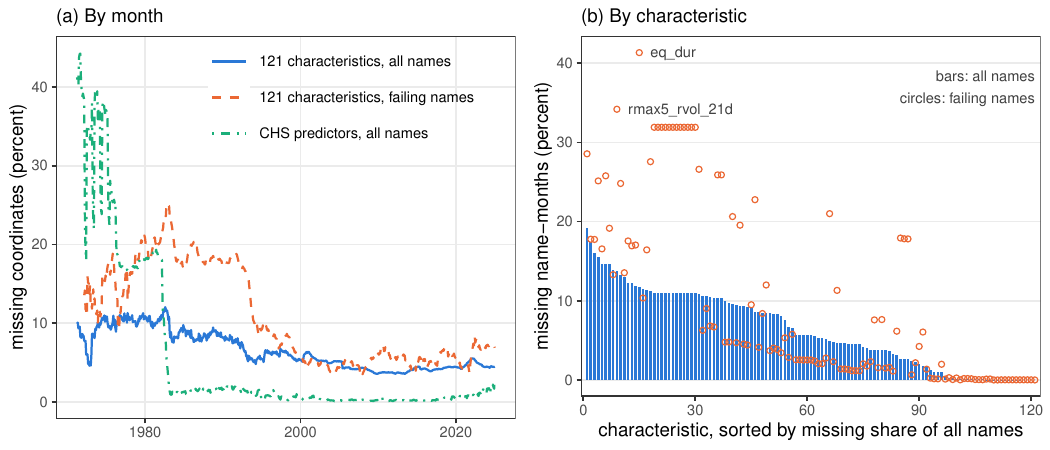}
\caption{Missing covariates in the $646$ test months. (a) Share of missing coordinates by month, for all names and for failing names on the $121$ characteristics (failing names pooled over the trailing twelve months), and for all names on the eight predictors of \citet{campbellhilscherszilagyi08}. (b) Share of name-months in which each of the $121$ characteristics is missing, sorted, for all names (bars) and for failing names (circles).}
\label{fig_missing}
\end{figure}

The screen of $121$ characteristics uses the full sample. A screen at the same threshold within each training window would keep $47$ to $93$ characteristics in the windows of $1971$--$1987$ and drop none of the $121$ after $2006$. Over $2007$--$2024$, where the two screens agree, the log-link fit of arc leads the ordinary logit by $0.036$ ($t=-8.3$). Figure~\ref{fig_missing} shows how much is missing. Over the test months a name misses $6.4$ percent of its characteristics on average and $1.7$ percent at the median, and $32$ percent of the name-months observe all $121$. Failing names miss more, $9.9$ percent on average, and only $14$ percent of them are complete. Missingness falls over time, from about nine percent of the coordinates in the $1970$s and $1980$s to under five percent after $2000$. Size, price and the last month's return are always observed. The characteristics missing most often need a long history, such as the five-year beta, or several years of accounts. For failing names the measures built from the daily returns of the last month, such as volatility, are missing in a third of the name-months against eleven percent overall, since a stock close to delisting trades thinly. Before $1983$ up to $44$ percent of the CHS coordinates are missing, and about one percent after. The accounting ratios are missing in $5$ to $9$ percent of all name-months, the market-based predictors in about one percent or less, and $10.5$ percent of the name-months miss at least one predictor.

\subsection{The masked dictionary}\label{app_masked}

Imputation is standard practice in the literature on firm characteristics. The most common rule sets a missing coordinate to the cross-sectional median or the median rank and evaluates the model there as if it had been observed \citep[e.g.,][]{shihaokellyxiu20}. \citet{bryzgalovalernerlettaupelger25} document that this rule is the most widely used one, that missingness in firm fundamentals follows systematic patterns, and that the median biases the imputed values. They impute from a latent factor model of the characteristics instead, and \citet{freybergeretal25} impute by regression on the observed characteristics. The masked dictionary~\eqref{maskeddict} of Section~\ref{sec_dict} instead evaluates a kernel dictionary at the conditional law of the missing coordinates given the observed ones.

Write $o=o(z)\subseteq\{1,\dots,d\}$ for the observed coordinates of a covariate vector $z$ in $\Zcal\subseteq\R^d$ and $\bar o$ for the missing ones, let $k$ be a kernel on complete covariate vectors with reproducing kernel Hilbert space $\HH_k$, and let $x_1,\dots,x_{N_{\rm ref}}$ be complete reference records, the name-months of the training window that observe every coordinate. With $\bm K\coloneqq[k(x_i,x_j)]_{i,j\le N_{\rm ref}}$, $m$ steps of the pivoted Cholesky decomposition select pivots $\Jcal=\{j_1,\dots,j_m\}$, each where the residual diagonal of the current approximation is largest, and return $\bm F\in\R^{N_{\rm ref}\times m}$ and an upper triangular $\bm R\in\R^{m\times m}$ with
\begin{equation}\label{pivchol}
    \bm K_{:,\Jcal}\bm R=\bm F,\qquad \bm R^\top\bm F_{\Jcal,:}=\bm I_m,\qquad \bm R\bm R^\top=\bm K_{\Jcal\Jcal}^{-1},
\end{equation}
at cost $O(N_{\rm ref}m^2)$ from the diagonal of $\bm K$ and the $m$ pivot columns \citep{HPS12,filipovicschneider24mom}. The Newton basis $\bm e(z)\coloneqq[k(z,x_{j_1}),\dots,k(z,x_{j_m})]\bm R$ is orthonormal in $\HH_k$, since its Gram matrix is $\bm R^\top\bm K_{\Jcal\Jcal}\bm R=\bm I_m$, and spans $\HH_\Jcal\coloneqq{\rm span}\{k(\cdot,x_j):j\in\Jcal\}$. With the linear kernel it is the grid dictionary of Section~\ref{sec_app}. For a covariate vector with missing coordinates, let $w_1(z),\dots,w_{N_{\rm ref}}(z)\ge0$ be weights that sum to one, an estimate of the conditional law of the missing coordinates $X_{\bar o}$ of a complete covariate vector $X$ given $X_o=z_o$, so that~\eqref{maskeddict} is the masked dictionary. We use the nearest-neighbor weights $w_i(z)=1/k_{\rm nn}$ on the $k_{\rm nn}$ reference records closest to $z$ in its observed coordinates \citep{troyanskaya01}. Their distances to a block of names are three matrix products, and a name costs $k_{\rm nn}$ evaluations of the basis. Averaging the index instead of the basis would give an index of rank up to $k_{\rm nn}$ per name and multiply the cost of every evaluation by $k_{\rm nn}$. With $\bm U=\bm L\bm L^\top$, it is the sum of squares of the conditional means of the functions $\bm e\bm L$.

\begin{proposition}[Masked dictionary]\label{prop_masked}
Let $P_\Jcal$ be the orthogonal projection of $\HH_k$ onto $\HH_\Jcal$ and $k_w(z)\coloneqq\sum_iw_i(z)\,k\big((z_o,x_{i,\bar o}),\cdot\big)$ the embedding of $z$ under the weights.
\begin{enumerate}
\item\label{prop_masked_proj} $\bm\phi(z)$ is the coordinate vector of $P_\Jcal k_w(z)$ in the basis $\bm e$, so $\bm\phi(z)\bm\phi(z')^\top=\langle P_\Jcal k_w(z),P_\Jcal k_w(z')\rangle_{\HH_k}$.
\item\label{prop_masked_lin} If $k(z,z')=1+\langle z,z'\rangle/d$, then $\bm e$ is affine and $\bm\phi(z)=\bm e\big(z_o,\sum_iw_i(z)x_{i,\bar o}\big)$, the basis at the conditional mean of the missing coordinates.
\end{enumerate}
\end{proposition}

\begin{proof}[Proof of Proposition~\ref{prop_masked}]
Write $\tilde x_i\coloneqq(z_o,x_{i,\bar o})$, so that $k_w(z)=\sum_iw_ik(\tilde x_i,\cdot)$.
\ref{prop_masked_proj}: The Newton basis is orthonormal in $\HH_k$ and spans $\HH_\Jcal$, so $P_\Jcal f=\sum_le_l\langle e_l,f\rangle$. By the reproducing property $\langle e_l,k(x,\cdot)\rangle=e_l(x)$, so $\langle\bm e,k_w(z)\rangle=\sum_iw_i\bm e(\tilde x_i)=\bm\phi(z)$ by linearity, and the inner product of two projections is that of their coordinate vectors.
\ref{prop_masked_lin}: The sections $k(\cdot,x_j)=1+\langle\cdot,x_j\rangle/d$, $j\in\Jcal$, are affine, hence so is $\bm e$, and an affine map commutes with the average in~\eqref{maskeddict}.
\end{proof}

Part~\ref{prop_masked_proj} identifies the masked dictionary with the Nystr\"om approximation \citep{chatalicetal22} of a conditional mean embedding \citep{songetal09}. By part~\ref{prop_masked_lin}, imputation at a fixed value is the masked dictionary of the linear kernel under weights that put all mass on that value, which is the conditional mean only where the conditional mean of the coordinate given the observed ones equals that value. For any other kernel the basis is nonlinear, and by Jensen's inequality its value at an imputed point differs in general from its conditional mean. The weights estimate the law among complete records. Under missingness at random \citep{rubin76} the law of the missing coordinates of a name given its observed ones and its pattern is their law given the observed ones, and the weights estimate it when, in addition, the indicator of a complete record is independent of the missing coordinates given the observed ones, a condition that missingness at random does not imply, and, for the nearest-neighbor weights, when the observed coordinates of incomplete names lie in the support of those of the complete records. The share of missing characteristics, which predicts distress in its own right, remains a summary of the name (Section~\ref{sec_app_data}). Proposition~\ref{prop_masked_app} below bounds the error of the projection and relates the masked dictionaries of nested observation patterns.

\begin{proposition}\label{prop_masked_app}
In the setting of Proposition~\ref{prop_masked}, let $\varrho_\Jcal(x)^2\coloneqq k(x,x)-\|\bm e(x)\|_2^2$ be the residual of the decomposition at $x$, which equals the residual diagonal at the reference records.
\begin{enumerate}
\item\label{prop_masked_err} $0\le\|k_w(z)\|_{\HH_k}^2-\|\bm\phi(z)\|_2^2\le\max\{\varrho_\Jcal(z_o,x_{i,\bar o})^2:w_i(z)>0\}$.
\item\label{prop_masked_tower} If $w(z)$ is the conditional law of $X_{\bar o}$ given $X_o=z_o$, then $\bm\phi(z)=\E[\bm e(X)\mid X_o=z_o]$, and for observed sets $o\subseteq o'$ the dictionary of the finer pattern averages to that of the coarser one, $\E[\bm\phi(X_{o'})\mid X_o]=\bm\phi(X_o)$, with $\bm\phi(X_o)$ the masked dictionary of a record that reveals $X_o$.
\end{enumerate}
\end{proposition}

\begin{proof}
Write $\tilde x_i\coloneqq(z_o,x_{i,\bar o})$.
\ref{prop_masked_err}: By Proposition~\ref{prop_masked}\ref{prop_masked_proj} and Pythagoras, $\|k_w(z)\|^2-\|\bm\phi(z)\|_2^2=\|(I-P_\Jcal)k_w(z)\|^2\ge0$, and $\|(I-P_\Jcal)k(x,\cdot)\|^2=k(x,x)-\|\bm e(x)\|_2^2=\varrho_\Jcal(x)^2$. The triangle inequality gives $\|(I-P_\Jcal)k_w(z)\|\le\sum_iw_i\varrho_\Jcal(\tilde x_i)\le\max_{w_i>0}\varrho_\Jcal(\tilde x_i)$. At a reference record $x_i$, $\bm e(x_i)=\bm K_{i,\Jcal}\bm R=\bm F_{i,:}$ by~\eqref{pivchol}, so $\varrho_\Jcal(x_i)^2=\bm K_{ii}-\|\bm F_{i,:}\|_2^2$ is the residual diagonal.
\ref{prop_masked_tower}: The first statement is~\eqref{maskeddict} with the conditional law in place of the weights. The second is the tower property of conditional expectation, $\sigma(X_o)\subseteq\sigma(X_{o'})$.
\end{proof}

\emph{Implementation.} The weights use $k_{\rm nn}=32$ nearest neighbors. Characteristics observed in fewer than half of the window's name-months, a few balance-sheet items in the early windows, are set to the median rank for every name and are not completed. The summaries enter the grid and lin as raw columns, and arc and half lifted like a coordinate, the state summaries after the map $s\mapsto2s-1$ from their percentiles. Reference records are the window's name-months that observe every active characteristic, subsampled by a fixed stride to at most $20{,}000$, and nearest neighbors are by Euclidean distance over the observed coordinates. The ex post check of the first-order conditions uses the tolerances $10^{-6}$ for the gradient, $10^{-5}$ for complementarity, $10^{-9}$ for feasibility and $10^{-11}$ for the relative available decrease.

\subsection{Subsamples and penalties}\label{app_extra}

Table~\ref{tab_app_robust} reports the log scores of Section~\ref{sec_app_calib} by subsample, Table~\ref{tab_lamstudy_app} the log-link fits of Section~\ref{sec_app_dicts} at other penalties, and Table~\ref{tab_app_qlogit} the link and the cone of Section~\ref{sec_app_diag}.

Two choices of the scoring affect only the weaker comparisons. With the equal-weighted cosine estimator of Table~\ref{tab_e2}, $\nu=29$, the lead over the ordinary logit keeps $t=-7.8$ and the AUC gain $t=3.5$, while the leads over KLR in arc and in half, at $t=-1.5$ and $-1.4$, and the gap of the CHS logit, at $t=1.7$, are no longer significant at five percent, against the $t_{29}$ critical value of $2.05$. The clip at $10^{-6}$ binds mainly for the logistic models. KLR puts $33$ failures in arc and $34$ in half below it, the sum-of-squares fits at most one, and their scores move by at most $0.0015$ for clips from $10^{-4}$ to $10^{-10}$, against $0.019$ for KLR in arc. The log-link fit of arc has the best score at every such clip.

Table~\ref{tab_lamstudy_app} reports, ex post, the log-link fits of Section~\ref{sec_app_dicts} at other penalties. Two phantom events per direction, close to the risk-optimal count in the simulations (Appendix~\ref{sec_sim_lampred}), improve every log-link fit, by $0.005$ to $0.031$ points, and leave the ranking of the dictionaries unchanged. The rank-angle dictionaries are the best at every penalty of the table, while the order of arc and half depends on the penalty. Rule~\eqref{lamlog} loses in four of the five richer dictionaries, as in the simulation at $m=250$. The TIC choice, close to the best count in the simulations, also loses to the rule in four of the five richer dictionaries, by $0.017$ to $0.066$ points, and gains only in half and on the CHS predictors. In arc and half the face, and with it the TIC, is that of the representation the solver returns. The trace penalty selects among representations but need not do so uniquely, so the TIC of these dictionaries is descriptive.

Table~\ref{tab_app_qlogit} asks whether the link or the index class drives the gain. A logistic regression whose index is minus a sum of squares, which Section~\ref{sec_loglink} identifies with the log-link fit of a rare event to first order, trails the log-link fit by only $0.003$ points, and KLR in the squared kernel, which drops the cone, does no better. At this penalty the gain comes from the index class, a negative quadratic form in a rank-angle dictionary, more than from the link. KLR in $k^2$ changes the penalty along with the cone, from the trace to the squared Frobenius norm, and KLR in the dictionary is unpenalized, so these two rows do not separate the cone from the regularization.

\begin{table}[!tbp]
\centering
\caption{The link and the cone on the failure event. Log score as in Table~\ref{tab_app_main}, and its difference to the log-link fit of the same dictionary with the Newey--West $t$-statistic (six lags) of the monthly name-weighted differences in brackets. Further columns give the mean monthly AUC, the calibration RMSE in units of $10^{-3}$ and the share of all failures in the top percentile of each month. Logit, $-$SOS index: the logistic regression ${\rm expit}(-\bm\phi\bm U\bm\phi^\top)$, $\bm U\in\Sy_+^m$, at the trace penalty $\lambda=T^{-1/2}$ of the sum-of-squares fits. KLR in $k^2$: ${\rm expit}(\bm\phi\bm A\bm\phi^\top)$ with $\bm A$ symmetric under the penalty $\lambda\|\bm A\|_F^2$, kernel logistic regression in the squared kernel restricted to the dictionary, fitted on every eighth test month, since its $m(m+1)/2$ free parameters make it the costliest fit. KLR: kernel logistic regression in the dictionary, unpenalized.}
\label{tab_app_qlogit}
\begin{footnotesize}
\setlength{\tabcolsep}{4pt}
\begin{tabular}{llcr@{\,}lccc}
\toprule
 & & log score & \multicolumn{2}{c}{vs log link} & AUC & RMSE & top $1\%$ \\
\midrule
\multicolumn{8}{l}{\emph{all $646$ test months}} \\
arc, $m=243$ & log link & $1.137$ & & & $0.940$ & $0.65$ & $0.392$ \\
 & logit, $-$SOS index & $1.140$ & $+0.003$ & $(6.9)$ & $0.940$ & $0.75$ & $0.390$ \\
lg, $m=200$ & log link & $1.192$ & & & $0.933$ & $0.82$ & $0.344$ \\
 & logit, $-$SOS index & $1.195$ & $+0.003$ & $(7.7)$ & $0.933$ & $0.90$ & $0.342$ \\
\addlinespace
\multicolumn{8}{l}{\emph{every eighth test month, $81$ months}} \\
arc, $m=243$ & log link & $1.216$ & & & $0.940$ & $0.90$ & $0.360$ \\
 & logit, $-$SOS index & $1.219$ & $+0.003$ & $(3.2)$ & $0.940$ & $1.01$ & $0.359$ \\
 & KLR in $k^2$ & $1.223$ & $+0.007$ & $(2.2)$ & $0.937$ & $0.82$ & $0.362$ \\
 & KLR & $1.231$ & $+0.014$ & $(0.6)$ & $0.927$ & $0.33$ & $0.392$ \\
\bottomrule
\end{tabular}
\end{footnotesize}
\end{table}

\begin{table}[!tbp]
\centering
\caption{Subsamples of the failure event. Log score of the ordinary logit in units of $10^{-2}$, and the difference of each model to it with the Newey--West $t$-statistic (six lags) of the monthly name-weighted differences in brackets. Complete CHS: test names with all eight predictors of \citet{campbellhilscherszilagyi08} observed. Halves: event months March $1971$ to January $1998$ and February $1998$ to December $2024$, $323$ test months each. The arc fits are as in Table~\ref{tab_app_main}. With the EWC estimator of Table~\ref{tab_e2} every comparison keeps its sign, and only the identity-link lead in the first half ($t=-2.0$) and the two half-sample differences of the CHS logit with missingness indicators ($t=-2.0$ and $2.0$) lose significance at five percent. CHS, indicators: the variant of Table~\ref{tab_app_main}.}
\label{tab_app_robust}
\begin{footnotesize}
\setlength{\tabcolsep}{3pt}
\resizebox{\textwidth}{!}{%
\begin{tabular}{lccr@{\,}lr@{\,}lr@{\,}lr@{\,}lr@{\,}l}
\toprule
 & name-months & failures & \multicolumn{2}{c}{ordinary logit} & \multicolumn{2}{c}{CHS logit} & \multicolumn{2}{c}{CHS, indicators} & \multicolumn{2}{c}{arc, log link} & \multicolumn{2}{c}{arc, identity link} \\
\midrule
all names & $2{,}965{,}609$ & $7{,}371$ & $1.185$ & & $+0.015$ & $(2.4)$ & $-0.001$ & $(-0.3)$ & $-0.048$ & $(-11.5)$ & $-0.005$ & $(-1.0)$ \\
complete CHS & $2{,}653{,}972$ & $6{,}318$ & $1.115$ & & $+0.003$ & $(0.7)$ & $+0.004$ & $(0.8)$ & $-0.041$ & $(-9.8)$ & $+0.004$ & $(1.0)$ \\
first half & $1{,}407{,}429$ & $3{,}214$ & $1.135$ & & $+0.011$ & $(1.1)$ & $-0.019$ & $(-2.5)$ & $-0.051$ & $(-7.4)$ & $-0.020$ & $(-2.8)$ \\
second half & $1{,}558{,}180$ & $4{,}157$ & $1.231$ & & $+0.018$ & $(2.7)$ & $+0.015$ & $(2.3)$ & $-0.046$ & $(-9.1)$ & $+0.008$ & $(1.3)$ \\
\bottomrule
\end{tabular}}
\end{footnotesize}
\end{table}

\begin{table}[!tp]
\centering
\caption{The penalty of the log-link fit in the application, $646$ test months. Log score at the rule as in Table~\ref{tab_app_main}, and the difference of each other penalty to it, with the Newey--West $t$-statistic (six lags) of the monthly name-weighted differences in brackets. Columns as in Table~\ref{tab_lamstudy}, with rule~\eqref{lamlog} calibrated at $m=11$. Fits that stopped before convergence, two of $5{,}168$ at rule~\eqref{lamlog} and at most $0.3$ percent at any other penalty, are scored as they stand. With the EWC estimator of Table~\ref{tab_e2} every comparison keeps its sign, and only rule~\eqref{lamlog} on the CHS predictors in arc ($t=1.8$) loses significance at five percent.}
\label{tab_lamstudy_app}
\begin{footnotesize}
\setlength{\tabcolsep}{3pt}
\begin{tabular}{lc*{5}{r@{\,}l}}
\toprule
 & rule & \multicolumn{2}{c}{\eqref{lamlog}} & \multicolumn{2}{c}{$\lambda T=1$} & \multicolumn{2}{c}{$\lambda T=2$} & \multicolumn{2}{c}{$\lambda T=4$} & \multicolumn{2}{c}{TIC} \\
\midrule
\multicolumn{12}{l}{\quad\emph{$121$ characteristics}} \\
lin, $m=122$ & $1.202$ & $+0.068$ & $(11.5)$ & $-0.014$ & $(-6.1)$ & $-0.018$ & $(-9.8)$ & $-0.013$ & $(-10.5)$ & $+0.066$ & $(6.5)$ \\
lg, $m=200$ & $1.192$ & $+0.065$ & $(9.4)$ & $-0.013$ & $(-5.2)$ & $-0.022$ & $(-10.4)$ & $-0.017$ & $(-12.0)$ & $+0.025$ & $(4.1)$ \\
lgz, $m=200$ & $1.174$ & $+0.085$ & $(11.3)$ & $-0.004$ & $(-1.6)$ & $-0.015$ & $(-8.7)$ & $-0.013$ & $(-11.4)$ & $+0.038$ & $(5.8)$ \\
arc, $m=243$ & $1.137$ & $+0.106$ & $(10.1)$ & $-0.003$ & $(-0.9)$ & $-0.014$ & $(-7.0)$ & $-0.013$ & $(-10.4)$ & $+0.017$ & $(3.4)$ \\
half, $m=243$ & $1.172$ & $-0.010$ & $(-1.7)$ & $-0.041$ & $(-11.7)$ & $-0.031$ & $(-12.0)$ & $-0.019$ & $(-12.2)$ & $-0.015$ & $(-2.5)$ \\
\multicolumn{12}{l}{\quad\emph{CHS predictors}} \\
lin, $m=9$ & $1.209$ & $+0.002$ & $(0.8)$ & $-0.004$ & $(-2.6)$ & $-0.005$ & $(-3.9)$ & $-0.004$ & $(-4.9)$ & $-0.003$ & $(-2.1)$ \\
arc, $m=17$ & $1.128$ & $+0.011$ & $(2.1)$ & $-0.009$ & $(-4.3)$ & $-0.010$ & $(-6.6)$ & $-0.008$ & $(-8.3)$ & $0.000$ & $(0.0)$ \\
half, $m=17$ & $1.147$ & $-0.024$ & $(-7.9)$ & $-0.021$ & $(-12.0)$ & $-0.016$ & $(-12.9)$ & $-0.010$ & $(-12.9)$ & $-0.021$ & $(-5.5)$ \\
\bottomrule
\end{tabular}
\end{footnotesize}
\end{table}

\end{appendices}

\clearpage
\ifnum\joesubmission=1
\section*{Declaration of generative AI and AI-assisted technologies in the manuscript preparation process}
During the preparation of this work the author used Claude (Anthropic) in order to check proofs, search the literature and edit the text. Its use for the code of the simulations and the application is stated in Section~\ref{sec_sim_design}. After using this tool, the author reviewed and edited the content as needed and takes full responsibility for the content of the publication.
\fi
\bibliographystyle{plainnat}
\bibliography{master}

\begin{thebibliography}{101}
\providecommand{\natexlab}[1]{#1}
\providecommand{\url}[1]{\texttt{#1}}
\expandafter\ifx\csname urlstyle\endcsname\relax
  \providecommand{\doi}[1]{doi: #1}\else
  \providecommand{\doi}{doi: \begingroup \urlstyle{rm}\Url}\fi

\bibitem[Alizadeh et~al.(1997)Alizadeh, Haeberly, and
  Overton]{alizadehhaeberlyoverton97}
Farid Alizadeh, Jean-Pierre~A. Haeberly, and Michael~L. Overton.
\newblock Complementarity and nondegeneracy in semidefinite programming.
\newblock \emph{Mathematical Programming}, 77\penalty0 (2):\penalty0 111--128,
  1997.

\bibitem[Alizadeh et~al.(2008)Alizadeh, Eckstein, Noyan, and
  Rudolf]{alizadehetal08}
Farid Alizadeh, Jonathan Eckstein, Nilay Noyan, and G{\'a}bor Rudolf.
\newblock Arrival rate approximation by nonnegative cubic splines.
\newblock \emph{Operations Research}, 56\penalty0 (1):\penalty0 140--156, 2008.

\bibitem[Almeida et~al.(2025)Almeida, Masini, and
  Schneider]{almeidamasinischneider25}
Caio Almeida, Ricardo Masini, and Paul Schneider.
\newblock Constrained polynomial likelihood.
\newblock \emph{Journal of Business \& Economic Statistics}, 43\penalty0
  (2):\penalty0 482--493, 2025.
\newblock \doi{10.1080/07350015.2024.2394587}.

\bibitem[Amisano and Giacomini(2007)]{amisanogiacomini07}
Gianni Amisano and Raffaella Giacomini.
\newblock Comparing density forecasts via weighted likelihood ratio tests.
\newblock \emph{Journal of Business \& Economic Statistics}, 25\penalty0
  (2):\penalty0 177--190, 2007.

\bibitem[Andrews(1991)]{andrews91}
Donald W.~K. Andrews.
\newblock Heteroskedasticity and autocorrelation consistent covariance matrix
  estimation.
\newblock \emph{Econometrica}, 59\penalty0 (3):\penalty0 817--858, 1991.

\bibitem[Andrews(1999)]{andrews1999b}
Donald W.~K. Andrews.
\newblock Estimation when a parameter is on a boundary.
\newblock \emph{Econometrica}, 67\penalty0 (6):\penalty0 1341--1383, 1999.

\bibitem[Aubin-Frankowski and Szab{\'o}(2020)]{aubinfrankowskiszabo20}
Pierre-Cyril Aubin-Frankowski and Zolt{\'a}n Szab{\'o}.
\newblock Hard shape-constrained kernel machines.
\newblock In \emph{Advances in Neural Information Processing Systems 33}, 2020.

\bibitem[Azizpour et~al.(2018)Azizpour, Giesecke, and
  Schwenkler]{azizpourgieseckeschwenkler18}
Shahriar Azizpour, Kay Giesecke, and Gustavo Schwenkler.
\newblock Exploring the sources of default clustering.
\newblock \emph{Journal of Financial Economics}, 129\penalty0 (1):\penalty0
  154--183, 2018.
\newblock \doi{10.1016/j.jfineco.2018.04.008}.

\bibitem[Bach(2010)]{bach10}
Francis Bach.
\newblock Self-concordant analysis for logistic regression.
\newblock \emph{Electronic Journal of Statistics}, 4:\penalty0 384--414, 2010.
\newblock \doi{10.1214/09-EJS521}.

\bibitem[Banerjee et~al.(2005{\natexlab{a}})Banerjee, Guo, and
  Wang]{banerjeeguowang05}
Arindam Banerjee, Xin Guo, and Hui Wang.
\newblock On the optimality of conditional expectation as a {B}regman
  predictor.
\newblock \emph{IEEE Transactions on Information Theory}, 51\penalty0
  (7):\penalty0 2664--2669, 2005{\natexlab{a}}.

\bibitem[Banerjee et~al.(2005{\natexlab{b}})Banerjee, Merugu, Dhillon, and
  Ghosh]{banerjeeetal05jmlr}
Arindam Banerjee, Srujana Merugu, Inderjit~S. Dhillon, and Joydeep Ghosh.
\newblock Clustering with {B}regman divergences.
\newblock \emph{Journal of Machine Learning Research}, 6:\penalty0 1705--1749,
  2005{\natexlab{b}}.

\bibitem[Billingsley(1961)]{billingsley61}
Patrick Billingsley.
\newblock The {L}indeberg-{L}\'evy theorem for martingales.
\newblock \emph{Proceedings of the American Mathematical Society}, 12\penalty0
  (5):\penalty0 788--792, 1961.

\bibitem[Birgin et~al.(2000)Birgin, Mart{\'\i}nez, and
  Raydan]{birginmartinezraydan00}
Ernesto~G. Birgin, Jos{\'e}~Mario Mart{\'\i}nez, and Marcos Raydan.
\newblock Nonmonotone spectral projected gradient methods on convex sets.
\newblock \emph{SIAM Journal on Optimization}, 10\penalty0 (4):\penalty0
  1196--1211, 2000.

\bibitem[Bonnans and Shapiro(2000)]{bonnansshapiro00}
J.~Fr{\'e}d{\'e}ric Bonnans and Alexander Shapiro.
\newblock \emph{Perturbation Analysis of Optimization Problems}.
\newblock Springer, New York, 2000.

\bibitem[Bradley(2005)]{bra_05}
Richard~C. Bradley.
\newblock Basic properties of strong mixing conditions. {A} survey and some
  open questions.
\newblock \emph{Probability Surveys}, 2:\penalty0 107--144, 2005.
\newblock \doi{10.1214/154957805100000104}.

\bibitem[Bregman(1967)]{bregman67}
Lev~M. Bregman.
\newblock The relaxation method of finding the common points of convex sets and
  its application to the solution of problems in convex programming.
\newblock \emph{USSR Computational Mathematics and Mathematical Physics},
  7\penalty0 (3):\penalty0 200--217, 1967.

\bibitem[Bryzgalova et~al.(2025)Bryzgalova, Lerner, Lettau, and
  Pelger]{bryzgalovalernerlettaupelger25}
Svetlana Bryzgalova, Sven Lerner, Martin Lettau, and Markus Pelger.
\newblock Missing financial data.
\newblock \emph{Review of Financial Studies}, 38\penalty0 (3):\penalty0
  803--882, 2025.

\bibitem[Campbell et~al.(2008)Campbell, Hilscher, and
  Szilagyi]{campbellhilscherszilagyi08}
John~Y. Campbell, Jens Hilscher, and Jan Szilagyi.
\newblock In search of distress risk.
\newblock \emph{Journal of Finance}, 63\penalty0 (6):\penalty0 2899--2939,
  2008.

\bibitem[Chatalic et~al.(2022)Chatalic, Schreuder, Rosasco, and
  Rudi]{chatalicetal22}
Antoine Chatalic, Nicolas Schreuder, Lorenzo Rosasco, and Alessandro Rudi.
\newblock Nystr{\"o}m kernel mean embeddings.
\newblock In \emph{Proceedings of the 39th International Conference on Machine
  Learning}, volume 162 of \emph{Proceedings of Machine Learning Research},
  pages 3006--3024, 2022.

\bibitem[Chava and Jarrow(2004)]{chavajarrow04}
Sudheer Chava and Robert~A. Jarrow.
\newblock Bankruptcy prediction with industry effects.
\newblock \emph{Review of Finance}, 8\penalty0 (4):\penalty0 537--569, 2004.

\bibitem[Chen and Fang(2019)]{chenfang19}
Qihui Chen and Zheng Fang.
\newblock Improved inference on the rank of a matrix.
\newblock \emph{Quantitative Economics}, 10\penalty0 (4):\penalty0 1787--1824,
  2019.

\bibitem[Chernoff(1954)]{chernoff54}
Herman Chernoff.
\newblock On the distribution of the likelihood ratio.
\newblock \emph{Annals of Mathematical Statistics}, 25\penalty0 (3):\penalty0
  573--578, 1954.

\bibitem[Crammer and Globerson(2006)]{crammerglobersonuai06}
Koby Crammer and Amir Globerson.
\newblock Discriminative learning via semidefinite probabilistic models.
\newblock In \emph{Proceedings of the Twenty-Second Conference on Uncertainty
  in Artificial Intelligence (UAI)}, 2006.

\bibitem[Daley and Vere-Jones(2003)]{daleyverejones03}
Daryl~J. Daley and David Vere-Jones.
\newblock \emph{An Introduction to the Theory of Point Processes. {V}olume {I}:
  Elementary Theory and Methods}.
\newblock Springer, New York, second edition, 2003.

\bibitem[Das et~al.(2007)Das, Duffie, Kapadia, and Saita]{dasetal07}
Sanjiv~R. Das, Darrell Duffie, Nikunj Kapadia, and Leandro Saita.
\newblock Common failings: How corporate defaults are correlated.
\newblock \emph{Journal of Finance}, 62\penalty0 (1):\penalty0 93--117, 2007.

\bibitem[Diebold and Mariano(1995)]{dieboldmariano95}
Francis~X. Diebold and Roberto~S. Mariano.
\newblock Comparing predictive accuracy.
\newblock \emph{Journal of Business \& Economic Statistics}, 13\penalty0
  (3):\penalty0 253--263, 1995.

\bibitem[Donald et~al.(2007)Donald, Fortuna, and
  Pipiras]{donaldfortunapipiras07}
Stephen~G. Donald, Nat{\'e}rcia Fortuna, and Vladas Pipiras.
\newblock On rank estimation in symmetric matrices: The case of indefinite
  matrix estimators.
\newblock \emph{Econometric Theory}, 23\penalty0 (6):\penalty0 1217--1232,
  2007.

\bibitem[Doukhan(1994)]{dou_94}
Paul Doukhan.
\newblock \emph{Mixing: Properties and Examples}, volume~85 of \emph{Lecture
  Notes in Statistics}.
\newblock Springer, New York, 1994.

\bibitem[Duffie et~al.(2007)Duffie, Saita, and Wang]{DuffieSaitaWang2007}
Darrell Duffie, Leandro Saita, and Ke~Wang.
\newblock Multi-period corporate default prediction with stochastic covariates.
\newblock \emph{Journal of Financial Economics}, 83\penalty0 (3):\penalty0
  635--665, 2007.

\bibitem[Duffie et~al.(2009)Duffie, Eckner, Horel, and
  Saita]{duffieecknerhorelsaita09}
Darrell Duffie, Andreas Eckner, Guillaume Horel, and Leandro Saita.
\newblock Frailty correlated default.
\newblock \emph{Journal of Finance}, 64\penalty0 (5):\penalty0 2089--2123,
  2009.

\bibitem[Fang and Santos(2019)]{fangsantos19}
Zheng Fang and Andres Santos.
\newblock Inference on directionally differentiable functions.
\newblock \emph{Review of Economic Studies}, 86\penalty0 (1):\penalty0
  377--412, 2019.

\bibitem[Filipovic and Schneider(2026)]{filipovicschneider24mom}
Damir Filipovic and Paul Schneider.
\newblock Conditional moment kernels: Joint modeling of expected returns and
  covariances in unbalanced panels.
\newblock SFI working paper series; formerly circulated as ``Joint Estimation
  of Conditional Mean and Covariance for Unbalanced Panels'', 2026.

\bibitem[Flaxman et~al.(2017)Flaxman, Teh, and Sejdinovic]{flaxman17}
Seth Flaxman, Yee~Whye Teh, and Dino Sejdinovic.
\newblock {P}oisson intensity estimation with reproducing kernels.
\newblock \emph{Electronic Journal of Statistics}, 11\penalty0 (2):\penalty0
  5081--5104, 2017.

\bibitem[Freyberger et~al.(2025)Freyberger, H{\"o}ppner, Neuhierl, and
  Weber]{freybergeretal25}
Joachim Freyberger, Bj{\"o}rn H{\"o}ppner, Andreas Neuhierl, and Michael Weber.
\newblock Missing data in asset pricing panels.
\newblock \emph{Review of Financial Studies}, 38\penalty0 (3):\penalty0
  760--802, 2025.

\bibitem[Gallant and Nychka(1987)]{gallantnychka87}
A.~Ronald Gallant and Douglas~W. Nychka.
\newblock Semi-nonparametric maximum likelihood estimation.
\newblock \emph{Econometrica}, 55:\penalty0 363--390, 1987.

\bibitem[Geyer(1994)]{gey_94}
Charles~J. Geyer.
\newblock On the asymptotics of constrained {M}-estimation.
\newblock \emph{Annals of Statistics}, 22\penalty0 (4):\penalty0 1993--2010,
  1994.
\newblock \doi{10.1214/aos/1176325768}.

\bibitem[Giacomini and White(2006)]{giacominiwhite06}
Raffaella Giacomini and Halbert White.
\newblock Tests of conditional predictive ability.
\newblock \emph{Econometrica}, 74\penalty0 (6):\penalty0 1545--1578, 2006.

\bibitem[Giesecke and Schwenkler(2018)]{gieseckeschwenkler18}
Kay Giesecke and Gustavo Schwenkler.
\newblock Filtered likelihood for point processes.
\newblock \emph{Journal of Econometrics}, 204\penalty0 (1):\penalty0 33 -- 53,
  2018.
\newblock ISSN 0304-4076.
\newblock \doi{https://doi.org/10.1016/j.jeconom.2017.11.011}.
\newblock URL
  \url{http://www.sciencedirect.com/science/article/pii/S0304407618300058}.

\bibitem[Gneiting and Raftery(2007)]{gneitingraftery07}
Tilmann Gneiting and Adrian~E. Raftery.
\newblock Strictly proper scoring rules, prediction, and estimation.
\newblock \emph{Journal of the American Statistical Association}, 102\penalty0
  (477):\penalty0 359--378, 2007.

\bibitem[Gouri{\'e}roux et~al.(1984{\natexlab{a}})Gouri{\'e}roux, Monfort, and
  Trognon]{gmt84poisson}
Christian Gouri{\'e}roux, Alain Monfort, and Alain Trognon.
\newblock Pseudo maximum likelihood methods: Applications to {P}oisson models.
\newblock \emph{Econometrica}, 52\penalty0 (3):\penalty0 701--720,
  1984{\natexlab{a}}.

\bibitem[Gouri{\'e}roux et~al.(1984{\natexlab{b}})Gouri{\'e}roux, Monfort, and
  Trognon]{gmt84theory}
Christian Gouri{\'e}roux, Alain Monfort, and Alain Trognon.
\newblock Pseudo maximum likelihood methods: Theory.
\newblock \emph{Econometrica}, 52\penalty0 (3):\penalty0 681--700,
  1984{\natexlab{b}}.

\bibitem[Gu et~al.(2020)Gu, Kelly, and Xiu]{shihaokellyxiu20}
Shihao Gu, Bryan Kelly, and Dacheng Xiu.
\newblock Empirical asset pricing via machine learning.
\newblock \emph{Review of Financial Studies}, 33\penalty0 (5):\penalty0
  2223--2273, 2020.
\newblock \doi{10.1093/rfs/hhaa009}.

\bibitem[Harbrecht et~al.(2012)Harbrecht, Peters, and Schneider]{HPS12}
Helmut Harbrecht, Michael Peters, and Reinhold Schneider.
\newblock On the low-rank approximation by the pivoted {C}holesky
  decomposition.
\newblock \emph{Applied Numerical Mathematics}, 62:\penalty0 428--440, 2012.

\bibitem[Hausman(1978)]{hausman78}
Jerry~A. Hausman.
\newblock Specification tests in econometrics.
\newblock \emph{Econometrica}, 46\penalty0 (6):\penalty0 1251--1271, 1978.

\bibitem[Hjort and Pollard(1993)]{hjortpollard93}
Nils~Lid Hjort and David Pollard.
\newblock Asymptotics for minimisers of convex processes.
\newblock Statistical research report, University of Oslo, 1993.
\newblock arXiv:1107.3806.

\bibitem[Ibragimov and Linnik(1971)]{ibr_lin_71}
Ildar~A. Ibragimov and Yuri~V. Linnik.
\newblock \emph{Independent and Stationary Sequences of Random Variables}.
\newblock Wolters-Noordhoff, Groningen, 1971.

\bibitem[Jensen et~al.(2023)Jensen, Kelly, and Pedersen]{jensenkellypedersen23}
Theis~Ingerslev Jensen, Bryan Kelly, and Lasse~Heje Pedersen.
\newblock Is there a replication crisis in finance?
\newblock \emph{Journal of Finance}, 78\penalty0 (5):\penalty0 2465--2518,
  2023.

\bibitem[Knight and Fu(2000)]{knightfu00}
Keith Knight and Wenjiang Fu.
\newblock Asymptotics for lasso-type estimators.
\newblock \emph{Annals of Statistics}, 28\penalty0 (5):\penalty0 1356--1378,
  2000.

\bibitem[Krengel(1985)]{kre_85}
Ulrich Krengel.
\newblock \emph{Ergodic Theorems}, volume~6 of \emph{De Gruyter Studies in
  Mathematics}.
\newblock Walter de Gruyter, Berlin, 1985.

\bibitem[Kulesza and Taskar(2012)]{kulesza2012determinantal}
Alex Kulesza and Ben Taskar.
\newblock \emph{Determinantal Point Processes for Machine Learning}, volume
  5(2--3) of \emph{Foundations and Trends in Machine Learning}.
\newblock Now Publishers, 2012.

\bibitem[Lasserre(2010)]{lasserre10}
Jean~Bernard Lasserre.
\newblock \emph{Moments, Positive Polynomials and Their Applications}, volume~1
  of \emph{Imperial College Press Optimization Series}.
\newblock Imperial College Press, 2010.

\bibitem[Lazarus et~al.(2018)Lazarus, Lewis, Stock, and Watson]{lazarusetal18}
Eben Lazarus, Daniel~J. Lewis, James~H. Stock, and Mark~W. Watson.
\newblock {HAR} inference: Recommendations for practice.
\newblock \emph{Journal of Business \& Economic Statistics}, 36\penalty0
  (4):\penalty0 541--559, 2018.

\bibitem[Liang and Zeger(1986)]{liangzeger86}
Kung-Yee Liang and Scott~L. Zeger.
\newblock Longitudinal data analysis using generalized linear models.
\newblock \emph{Biometrika}, 73\penalty0 (1):\penalty0 13--22, 1986.

\bibitem[Macchi(1975)]{macchi1975coincidence}
Odile Macchi.
\newblock The coincidence approach to stochastic point processes.
\newblock \emph{Advances in Applied Probability}, 7\penalty0 (1):\penalty0
  83--122, 1975.

\bibitem[Marschner and Gillett(2012)]{marschnergillett12}
Ian~C. Marschner and Alexandra~C. Gillett.
\newblock Relative risk regression: Reliable and flexible methods for
  log-binomial models.
\newblock \emph{Biostatistics}, 13\penalty0 (1):\penalty0 179--192, 2012.

\bibitem[Marteau-Ferey et~al.(2019)Marteau-Ferey, Ostrovskii, Bach, and
  Rudi]{marteauferey19}
Ulysse Marteau-Ferey, Dmitrii Ostrovskii, Francis Bach, and Alessandro Rudi.
\newblock Beyond least-squares: Fast rates for regularized empirical risk
  minimization through self-concordance.
\newblock In \emph{Proceedings of the 32nd Conference on Learning Theory
  (COLT)}, volume~99 of \emph{PMLR}, pages 2294--2340, 2019.

\bibitem[Marteau-Ferey et~al.(2020)Marteau-Ferey, Bach, and
  Rudi]{marteauferey2020nonparametric}
Ulysse Marteau-Ferey, Francis Bach, and Alessandro Rudi.
\newblock Non-parametric models for non-negative functions.
\newblock In \emph{Advances in Neural Information Processing Systems 33}, 2020.

\bibitem[McCullagh and Nelder(1989)]{mccullaghnelder89}
Peter McCullagh and John~A. Nelder.
\newblock \emph{Generalized Linear Models}.
\newblock Chapman and Hall, London, second edition, 1989.

\bibitem[Merlev\`ede et~al.(2009)Merlev\`ede, Peligrad, and
  Rio]{mer_pel_rio_09}
Florence Merlev\`ede, Magda Peligrad, and Emmanuel Rio.
\newblock Bernstein inequality and moderate deviations under strong mixing
  conditions.
\newblock In \emph{High Dimensional Probability {V}: The {L}uminy Volume},
  volume~5 of \emph{Institute of Mathematical Statistics Collections}, pages
  273--292. Institute of Mathematical Statistics, Beachwood, OH, 2009.
\newblock \doi{10.1214/09-IMSCOLL518}.

\bibitem[Milz(2023)]{mil_23}
Johannes Milz.
\newblock Sample average approximations of strongly convex stochastic programs
  in {H}ilbert spaces.
\newblock \emph{Optimization Letters}, 17\penalty0 (2):\penalty0 471--492,
  2023.
\newblock \doi{10.1007/s11590-022-01888-4}.

\bibitem[Murphy and Winkler(1977)]{murphywinkler77}
Allan~H. Murphy and Robert~L. Winkler.
\newblock Reliability of subjective probability forecasts of precipitation and
  temperature.
\newblock \emph{Journal of the Royal Statistical Society: Series C (Applied
  Statistics)}, 26\penalty0 (1):\penalty0 41--47, 1977.

\bibitem[Mutn{\'y} and Krause(2022)]{mutnykrause22}
Mojm{\'i}r Mutn{\'y} and Andreas Krause.
\newblock Sensing {C}ox processes via posterior sampling and positive bases.
\newblock In \emph{Proceedings of the 25th International Conference on
  Artificial Intelligence and Statistics (AISTATS)}, volume 151 of \emph{PMLR},
  2022.

\bibitem[Muzellec et~al.(2021)Muzellec, Bach, and Rudi]{muzellec21psd}
Boris Muzellec, Francis Bach, and Alessandro Rudi.
\newblock Learning {PSD}-valued functions using kernel sums-of-squares.
\newblock arXiv:2111.11306, 2021.

\bibitem[Nesterov(2018)]{nesterov18}
Yurii Nesterov.
\newblock \emph{Lectures on Convex Optimization}, volume 137 of \emph{Springer
  Optimization and Its Applications}.
\newblock Springer, Cham, second edition, 2018.
\newblock \doi{10.1007/978-3-319-91578-4}.

\bibitem[Ostrovskii and Bach(2021)]{ost_bac_21}
Dmitrii~M. Ostrovskii and Francis Bach.
\newblock Finite-sample analysis of {M}-estimators using self-concordance.
\newblock \emph{Electronic Journal of Statistics}, 15\penalty0 (1):\penalty0
  326--391, 2021.
\newblock \doi{10.1214/20-EJS1780}.

\bibitem[Papp and Alizadeh(2014)]{pappalizadeh14}
D{\'a}vid Papp and Farid Alizadeh.
\newblock Shape-constrained estimation using nonnegative splines.
\newblock \emph{Journal of Computational and Graphical Statistics}, 23\penalty0
  (1):\penalty0 211--231, 2014.

\bibitem[Pinelis(1994)]{pinelis94}
Iosif Pinelis.
\newblock Optimum bounds for the distributions of martingales in {B}anach
  spaces.
\newblock \emph{Annals of Probability}, 22\penalty0 (4):\penalty0 1679--1706,
  1994.

\bibitem[Politis et~al.(1999)Politis, Romano, and Wolf]{politisromanowolf99}
Dimitris~N. Politis, Joseph~P. Romano, and Michael Wolf.
\newblock \emph{Subsampling}.
\newblock Springer, New York, 1999.

\bibitem[P{\"o}tscher and Prucha(1997)]{PoetscherPrucha1997}
Benedikt~M. P{\"o}tscher and Ingmar~R. Prucha.
\newblock \emph{Dynamic Nonlinear Econometric Models: Asymptotic Theory}.
\newblock Springer, New York, 1997.

\bibitem[Rockafellar(1997)]{roc_97}
R.~Tyrrell Rockafellar.
\newblock \emph{Convex Analysis}.
\newblock Princeton Landmarks in Mathematics. Princeton University Press,
  Princeton, NJ, 1997.
\newblock Reprint of the 1970 original.

\bibitem[Rockafellar and Wets(1998)]{roc_wet_98}
R.~Tyrrell Rockafellar and Roger J.-B. Wets.
\newblock \emph{Variational Analysis}, volume 317 of \emph{Grundlehren der
  mathematischen Wissenschaften}.
\newblock Springer, Berlin, 1998.

\bibitem[Rubin(1976)]{rubin76}
Donald~B. Rubin.
\newblock Inference and missing data.
\newblock \emph{Biometrika}, 63\penalty0 (3):\penalty0 581--592, 1976.

\bibitem[Rudi and Ciliberto(2021)]{rudiciliberto21}
Alessandro Rudi and Carlo Ciliberto.
\newblock {PSD} representations for effective probability models.
\newblock In \emph{Advances in Neural Information Processing Systems 34}, 2021.

\bibitem[Sadhwani et~al.(2021)Sadhwani, Giesecke, and
  Sirignano]{sadhwanigieseckesirignano21}
Apaar Sadhwani, Kay Giesecke, and Justin Sirignano.
\newblock Deep learning for mortgage risk.
\newblock \emph{Journal of Financial Econometrics}, 19\penalty0 (2):\penalty0
  313--368, 2021.
\newblock \doi{10.1093/jjfinec/nbaa025}.

\bibitem[Santos~Silva and Tenreyro(2006)]{silvatenreyro06}
J.~M.~C. Santos~Silva and Silvana Tenreyro.
\newblock The log of gravity.
\newblock \emph{Review of Economics and Statistics}, 88\penalty0 (4):\penalty0
  641--658, 2006.

\bibitem[Savage(1971)]{savage71}
Leonard~J. Savage.
\newblock Elicitation of personal probabilities and expectations.
\newblock \emph{Journal of the American Statistical Association}, 66\penalty0
  (336):\penalty0 783--801, 1971.

\bibitem[Self and Liang(1987)]{selfliang87}
Steven~G. Self and Kung-Yee Liang.
\newblock Asymptotic properties of maximum likelihood estimators and likelihood
  ratio tests under nonstandard conditions.
\newblock \emph{Journal of the American Statistical Association}, 82\penalty0
  (398):\penalty0 605--610, 1987.

\bibitem[Shapiro(1989)]{shapiro89}
Alexander Shapiro.
\newblock Asymptotic properties of statistical estimators in stochastic
  programming.
\newblock \emph{Annals of Statistics}, 17\penalty0 (2):\penalty0 841--858,
  1989.

\bibitem[Shapiro(1997)]{shapiro97}
Alexander Shapiro.
\newblock First and second order analysis of nonlinear semidefinite programs.
\newblock \emph{Mathematical Programming}, 77:\penalty0 301--320, 1997.

\bibitem[Shapiro(2019)]{shapiro19}
Alexander Shapiro.
\newblock Statistical inference of semidefinite programming.
\newblock \emph{Mathematical Programming}, 174:\penalty0 77--97, 2019.

\bibitem[Shapiro et~al.(2014)Shapiro, Dentcheva, and
  Ruszczy\'{n}ski]{sha_den_rus_14}
Alexander Shapiro, Darinka Dentcheva, and Andrzej Ruszczy\'{n}ski.
\newblock \emph{Lectures on Stochastic Programming: Modeling and Theory},
  volume~16 of \emph{MOS-SIAM Series on Optimization}.
\newblock Society for Industrial and Applied Mathematics, Philadelphia, PA,
  second edition, 2014.
\newblock \doi{10.1137/1.9781611973433}.

\bibitem[Shibata(1989)]{shibata89}
Ritei Shibata.
\newblock Statistical aspects of model selection.
\newblock In Jan~C. Willems, editor, \emph{From Data to Model}, pages 215--240.
  Springer, Berlin, Heidelberg, 1989.
\newblock \doi{10.1007/978-3-642-75007-6_5}.

\bibitem[Shumway(1997)]{shumway97}
Tyler Shumway.
\newblock The delisting bias in {CRSP} data.
\newblock \emph{Journal of Finance}, 52\penalty0 (1):\penalty0 327--340, 1997.

\bibitem[Shumway(2001)]{shumway01}
Tyler Shumway.
\newblock Forecasting bankruptcy more accurately: A simple hazard model.
\newblock \emph{Journal of Business}, 74\penalty0 (1):\penalty0 101--124, 2001.

\bibitem[{\'S}mieja et~al.(2018){\'S}mieja, Struski, Tabor, Zieli{\'n}ski, and
  Spurek]{smiejaetal18}
Marek {\'S}mieja, {\L}ukasz Struski, Jacek Tabor, Bartosz Zieli{\'n}ski, and
  Przemys{\l}aw Spurek.
\newblock Processing of missing data by neural networks.
\newblock In \emph{Advances in Neural Information Processing Systems},
  volume~31, pages 2719--2729, 2018.

\bibitem[Smola et~al.(2005)Smola, Vishwanathan, and
  Hofmann]{smolavishwanathanhofmann05}
Alexander~J. Smola, S.~V.~N. Vishwanathan, and Thomas Hofmann.
\newblock Kernel methods for missing variables.
\newblock In \emph{Proceedings of the Tenth International Workshop on
  Artificial Intelligence and Statistics}, volume~R5 of \emph{Proceedings of
  Machine Learning Research}, pages 325--332, 2005.

\bibitem[Song et~al.(2009)Song, Huang, Smola, and Fukumizu]{songetal09}
Le~Song, Jonathan Huang, Alexander~J. Smola, and Kenji Fukumizu.
\newblock Hilbert space embeddings of conditional distributions with
  applications to dynamical systems.
\newblock In \emph{Proceedings of the 26th Annual International Conference on
  Machine Learning (ICML)}, pages 961--968, New York, 2009. Association for
  Computing Machinery.
\newblock \doi{10.1145/1553374.1553497}.

\bibitem[Tropp(2015)]{tropp15}
Joel~A. Tropp.
\newblock An introduction to matrix concentration inequalities.
\newblock \emph{Foundations and Trends in Machine Learning}, 8\penalty0
  (1--2):\penalty0 1--230, 2015.
\newblock \doi{10.1561/2200000048}.

\bibitem[Troyanskaya et~al.(2001)Troyanskaya, Cantor, Sherlock, Brown, Hastie,
  Tibshirani, Botstein, and Altman]{troyanskaya01}
Olga Troyanskaya, Michael Cantor, Gavin Sherlock, Pat Brown, Trevor Hastie,
  Robert Tibshirani, David Botstein, and Russ~B. Altman.
\newblock Missing value estimation methods for {DNA} microarrays.
\newblock \emph{Bioinformatics}, 17\penalty0 (6):\penalty0 520--525, 2001.

\bibitem[Tsuchida et~al.(2024)Tsuchida, Ong, and Sejdinovic]{tsuchida24}
Russell Tsuchida, Cheng~Soon Ong, and Dino Sejdinovic.
\newblock Exact, fast and expressive {P}oisson point processes via squared
  neural families.
\newblock \emph{Proceedings of the AAAI Conference on Artificial Intelligence},
  38\penalty0 (18):\penalty0 20559--20566, 2024.

\bibitem[Wacholder(1986)]{wacholder86}
Sholom Wacholder.
\newblock Binomial regression in {GLIM}: Estimating risk ratios and risk
  differences.
\newblock \emph{American Journal of Epidemiology}, 123:\penalty0 174--184,
  1986.

\bibitem[Walder and Bishop(2017)]{walderbishop17}
Christian~J. Walder and Adrian~N. Bishop.
\newblock Fast {B}ayesian intensity estimation for the permanental process.
\newblock In \emph{Proceedings of the 34th International Conference on Machine
  Learning (ICML)}, volume~70 of \emph{PMLR}, pages 3579--3588, 2017.

\bibitem[Wedderburn(1974)]{wedderburn74}
R.~W.~M. Wedderburn.
\newblock Quasi-likelihood functions, generalized linear models, and the
  {G}auss--{N}ewton method.
\newblock \emph{Biometrika}, 61\penalty0 (3):\penalty0 439--447, 1974.

\bibitem[White(1994)]{white94}
Halbert White.
\newblock \emph{Estimation, Inference and Specification Analysis}.
\newblock Number~22 in Econometric Society Monographs. Cambridge University
  Press, Cambridge, 1994.

\bibitem[Wong and Shen(1995)]{wongshen95}
Wing~Hung Wong and Xiaotong Shen.
\newblock Probability inequalities for likelihood ratios and convergence rates
  of sieve {MLE}s.
\newblock \emph{Annals of Statistics}, 23\penalty0 (2):\penalty0 339--362,
  1995.

\bibitem[Wooldridge(1991)]{wooldridge91}
Jeffrey~M. Wooldridge.
\newblock Specification testing and quasi-maximum-likelihood estimation.
\newblock \emph{Journal of Econometrics}, 48\penalty0 (1--2):\penalty0 29--55,
  1991.

\bibitem[Wooldridge(1994)]{wooldridge94}
Jeffrey~M. Wooldridge.
\newblock Estimation and inference for dependent processes.
\newblock In Robert~F. Engle and Daniel~L. McFadden, editors, \emph{Handbook of
  Econometrics}, volume~4, chapter~45, pages 2639--2738. North-Holland,
  Amsterdam, 1994.

\bibitem[Wooldridge(2001)]{Wooldridgebook2001}
Jeffrey~M. Wooldridge.
\newblock \emph{Econometric Analysis of Cross Section and Panel Data}.
\newblock MIT Press, Cambridge, MA, 2001.

\bibitem[Yu et~al.(2015)Yu, Wang, and Samworth]{yu_wan_sam_15}
Yi~Yu, Tengyao Wang, and Richard~J. Samworth.
\newblock A useful variant of the {D}avis--{K}ahan theorem for statisticians.
\newblock \emph{Biometrika}, 102\penalty0 (2):\penalty0 315--323, 2015.
\newblock \doi{10.1093/biomet/asv008}.

\bibitem[Zhu and Hastie(2005)]{zhuhastie05}
Ji~Zhu and Trevor Hastie.
\newblock Kernel logistic regression and the import vector machine.
\newblock \emph{Journal of Computational and Graphical Statistics}, 14\penalty0
  (1):\penalty0 185--205, 2005.

\bibitem[Zou(2004)]{zou04}
Guangyong Zou.
\newblock A modified {P}oisson regression approach to prospective studies with
  binary data.
\newblock \emph{American Journal of Epidemiology}, 159\penalty0 (7):\penalty0
  702--706, 2004.

\end{thebibliography}

\end{document}